\documentclass[11pt]{article}
 
\usepackage{geometry}
\usepackage{setspace}
\usepackage{bm}
\usepackage{graphicx}
\graphicspath{ {.} }
\usepackage{hyperref}
\usepackage{xr}
\newcommand{\ra}[1]{\renewcommand{\arraystretch}{#1}}

\usepackage{amssymb, amsmath, amsthm, amsfonts, mathtools, bbm}
\usepackage{mathrsfs}
\usepackage{float}
\usepackage{color, colortbl}
\usepackage[dvipsnames]{xcolor}
\usepackage{url}
\usepackage{enumitem}
\usepackage{multirow}
\usepackage{siunitx}
\usepackage{placeins}
\usepackage{authblk}
\usepackage{appendix}
\usepackage{lscape}
\usepackage{booktabs}
\usepackage{longtable}
\usepackage{stackrel}

\usepackage{comment}

\usepackage{soul}       
\usepackage{ulem}   
\usepackage{nicefrac}
\usepackage{microtype} %this causes ``Command \showhyphens has changed'', as \showhyphens is also defined in ``\usepackage{babel}''; but we need it for our formatting

\allowdisplaybreaks[1]

\usepackage{xr-hyper}
\usepackage{hyperref}
\hypersetup{
    colorlinks=true,
    linkcolor=blue,
    citecolor=blue,
    filecolor=magenta,
    urlcolor=blue
}
\usepackage{cleveref}

\usepackage[round]{natbib}

\usepackage{amsthm} %for begin{proof}

\newtheorem{theorem}{Theorem}
\newtheorem{proposition}{Proposition}

\newtheorem{lemma}{Lemma}
\newtheorem{remark}{Remark}
\newtheorem{assumption}{Assumption}
\newtheorem{definition}{Definition}

\usepackage{algorithm}
\usepackage{algorithmic}

\usepackage{bm}
\usepackage{amsmath}
\usepackage{multirow}
\usepackage{booktabs}
\usepackage{tabularx}
\usepackage{array}
\usepackage[T1]{fontenc}
\usepackage[utf8]{inputenc}
\usepackage[ngerman,english]{babel}
\usepackage{subcaption}
\def\E{\mathbb{E}}
\def\P{\mathbb{P}}

\def\tr{\mathrm{tr}}

\def\R{\mathbb{R}}

\usepackage[makeroom]{cancel}

\usepackage{fullpage}

\begin{document}
		
		\title{Federated Estimation of Loading Spaces in Tensor Factor Models under Differential Privacy}
                \author{Zetai Cen\thanks{Equal contribution. School of Mathematics, University of Bristol.
                Email: \url{zetai.cen@bristol.ac.uk}.
                Supported by Engineering and Physical Sciences Research Council (EP/Z531327/1).}}
                \author{Fan Wang\thanks{Equal contribution. School of Mathematics and Statistics, University of Melbourne.
                Email: \url{fan.wang.2@unimelb.edu.au}.}}
                \author{Haeran Cho\thanks{Corresponding author. School of Mathematics, University of Bristol.
                Email: \url{haeran.cho@bristol.ac.uk}.
                Supported by Engineering and Physical Sciences Research Council (EP/Z531327/1).}}
        \affil{}
		\date{}
		\maketitle
        
\begin{abstract}
Tensor factor models provide low-dimensional representations of
tensor-valued observations while preserving their multiway structure. When these observations are held by different clients, privacy constraints may prevent them from being pooled for centralized analysis. We study estimation of common loading spaces in Tucker tensor factor models under federated differential privacy, allowing heterogeneous sample sizes, privacy budgets, and factor strengths. The proposed method combines private stability testing, Huber aggregation, and Gaussian perturbation at each client with weighted aggregation at the server. We establish non-asymptotic upper bounds governed by a scaled harmonic mean of the squared local error rates, and prove minimax optimality up to logarithmic factors under strong factors, fixed factor ranks, and fixed noise variance. We further develop a private calibration procedure that selects the local tuning parameters and aggregation weights without knowledge of the noise level or client-specific factor strengths. The resulting adaptive estimator preserves the privacy guarantee and attains the error rate obtained with oracle tuning. We illustrate the finite-sample performance of the proposed estimators through simulation experiments and an application to a Meditation Traits EEG dataset.
\end{abstract}

\bigskip

\noindent
{\sl Keywords and phrases:}
Adaptive estimation;
differential privacy;
minimax optimality;
tensor factor model;
Tucker decomposition.

\newpage

\section{Introduction}\label{sec: intro}

Modern applications increasingly generate complex data with multiway structures. Examples include facial images indexed by person, lighting condition, and pixel location \citep{kolda2009tensor}; 
recommendation data indexed by user, item, and context \citep{karatzoglou2010multiverse}; 
and gene-expression data indexed by individual mouse, tissue type, and gene \citep{allen2012sparse}. 
The resulting tensors are often high-dimensional, making it desirable to represent the data parsimoniously while preserving the multiway structure as each tensor mode may exhibit a distinct source of variation. The Tucker tensor factor model provides such a representation through decomposing the data into a low-rank component driven by a core factor tensor, and an idiosyncratic noise, under which it is of interest to learn the low-dimensional spaces spanned by mode-specific loading matrices that describe the main patterns of variation along each mode \citep{cichocki2009nonnegative, Chenetal2022, Babiietal2025, barigozzi2022statistical}.

In many applications, tensor-valued observations are held by multiple clients, and privacy regulations or institutional restrictions may prevent them from being pooled for centralized analysis; consider e.g., patient imaging data held by hospitals, gene-expression measurements stored by research centers, and user--item data retained by platforms in different jurisdictions.  
Federated learning allows these clients to collaborate by transmitting locally computed model updates or summary statistics to a central server for aggregation, while retaining the raw observations locally~\citep{mcmahan2017communication}. However, keeping the raw data in local clients does not by itself does not provide a formal privacy guarantee, since transmitted outputs may still reveal information about individual records. 
Federated differential privacy (FDP) addresses this risk by requiring that any output released by each client to satisfy record-level differential privacy (DP), so that its distribution remains stable under the replacement of any single local record \citep{Dworketal2006, DworkRoth2014}.

Despite the growing literature on FDP estimation, methods and statistical theory specific for tensor data are underdeveloped, particularly on the trade-off between privacy and statistical efficiency. We propose to close this gap by studying the problem of estimating mode-wise loading spaces in Tucker tensor factor models under FDP, which comes with worst-case record-level FDP under arbitrary record replacements.

\subsection{Related literature}\label{subsec:related}

Our work lies at the intersection of tensor factor models, differentially private principal component analysis (PCA), and FDP.

\medskip\noindent\textbf{Tensor factor models.} 
Low-rank tensor models are commonly formulated through Tucker or canonical polyadic (CP) decompositions \citep{cichocki2009nonnegative, kolda2009tensor}. Building on the Tucker decomposition, statistical tensor factor models represent tensor-valued observations through latent core tensors and mode-specific loading matrices.
Generalizing from matrix factor models by \citet{Wangetal2019}, 
the seminal work by \citet{Chenetal2022} considered a Tucker-type model for high-dimensional tensor time series.
\citet{barigozzi2022statistical} studied a Tucker tensor factor model and developed an iterative projection estimator, deriving convergence rates and asymptotic distributions for the loading matrices and common factors.
\citet{Babiietal2025} showed that tensor PCA is rate-optimal under strong factors, while alternating least squares improves the convergence rate under weak factors.
More recent literature has seen a wide range of developments \citep[e.g.][]{ChenLam2024, Hanetal2024, CenLam2025, barigozzi2025tail, ZhangCenCho2026};
see a recent review in \cite{Bolivaretal2026}.
All these studies consider tensor factor models without privacy constraints.

\medskip\noindent\textbf{Differentially private PCA.}
PCA provides a low-dimensional representation of multivariate data by identifying orthogonal directions that capture the dominant variation \citep{jolliffe2002principal, jolliffe2016principal} and has popularly been adopted for estimation of factor models and other low-rank structures.
Early work on differentially private PCA studied fixed, bounded datasets without making distributional assumptions \citep{blum2005practical, chaudhuri2013near, dwork2014analyze}. 
Subsequent work considered private PCA in distributed settings under boundedness assumptions. \citet{wang2020principal} studied non-interactive locally private PCA and obtained optimal or near-optimal rates in low-dimensional and sparse high-dimensional settings. \citet{grammenos2020federated} developed a differentially private federated PCA for bounded streaming data under homogeneous sample sizes and privacy budgets. More recent work has considered unbounded observations under distributional assumptions. \citet{liu2022dppca} established nearly optimal rates for rank-one PCA under sub-Gaussian distributions. 
For general-rank PCA, \citet{cai2024optimal} obtained nearly minimax-optimal rates under the spiked covariance model, but their privacy guarantee holds only with high probability and for independent-copy replacements rather than arbitrary record replacements. \citet{dungler2026iterative} subsequently developed an iterative algorithm for general-rank PCA with worst-case privacy guarantees and established nearly matching upper and lower error bounds.
In all above, the focus lies in the study of PCA-based methods for a sequence of random vectors, and their guarantees do not extend to tensor-valaued data without careful investigation into the multiarray structure.

\medskip\noindent\textbf{Federated differential privacy.} 
This framework requires each client's complete release to satisfy DP with respect to its local dataset, while the central server has access only to these privatized releases \citep{wei2020federated, li2024federated}. 
FDP has been studied in a range of statistical problems, including mean estimation and linear regression \citep{li2024federated}, nonparametric regression and classification \citep{cai2024nonparametric, auddy2026minimax}, functional data analysis \citep{cai2024functional, xue2024optimal}, and Cox regression \citep{hung2025optimal}. 
In these works, minimax framework has been adopted for the quantification of the cost of privacy, establishing minimax lower bounds on the rates of estimation accompanied by matching upper bounds up to logarithmic factors.
Minimax-optimal rates have been established in these settings, up to logarithmic factors. 
Most closely related to our setting, \citet{Lietal2024_fdp} studied federated estimation of a common principal subspace and covariance matrix under the spiked covariance model. They obtained an upper bound matching the federated minimax lower bound up to logarithmic factors but as in \citet{cai2024optimal}, their privacy guarantee is limited in its formulation.

\subsection{List of contributions}

The main contributions of this paper are summarized as follows.

Firstly, we develop a novel loading-space estimator for Tucker tensor factor models under FDP, allowing for heterogeneous sample sizes, privacy budgets, and factor strengths. 
The proposed algorithm combines a private stability test performed on blockwise loading-space projectors, a Huber loss-based aggregation strategy, and Gaussian perturbations for constructing a private local projector. Then, the central server aggregates these local projectors using client-specific weights.
Through these steps, the algorithm guarantees worst-case record-level FDP under arbitrary record replacements, in contrast to existing studies which only provide high probability guarantees and permit independent-copy replacements \citep{cai2024optimal, Lietal2024_fdp}.

Secondly, we establish the non-asymptotic upper bound for the oracle procedure whose Huber radii, noise scales, and aggregation weights depend on the unknown noise level and client-specific factor strengths.
Our results show that these rates are governed by a scaled harmonic mean of the squared local error rates, which is no larger than the smallest squared local error rate and can vanish as the number of clients increases even when none of the local rates does. We derive a federated minimax lower bound with the same harmonic structure, establishing minimax optimality up to logarithmic factors under strong factors with fixed ranks and noise variance. 
Our upper bound is attained under the same worst-case FDP requirement as the lower bound
and, to the best of our knowledge, this is the first minimax-optimality result for low-rank structure estimation for tensor-valued data
under the worst-case FDP.

Thirdly, we develop a private calibration procedure that adaptively selects the local tuning parameters and aggregation weights without knowledge of the noise level or client-specific factor strengths.
In \cite{li2024federated}, clients share a common covariance matrix, so the common signal and noise parameters cancel from the normalized PCA aggregation weights. In our setting, we allow factor strengths to vary across clients, so the optimal weights additionally depend on unknown client-specific signal parameters. In addition, both their procedure and our oracle procedure require knowledge of the noise level and signal strength to calibrate the Gaussian perturbation added for privacy.  To address both sources of dependence, each client uses a separate calibration sample and a private threshold search to estimate the dispersion of its blockwise loading-space projectors.
The resulting scale estimate determines the Huber radii, the scale of perturbations, and the estimated local error rate used to compute the aggregation weights.
Our adaptive estimator retains the same worst-case privacy guarantee and, under the additional signal condition, attains the same error rate as the oracle estimator.

\subsection{Notation and organization}
\label{subsec:notation}

We use bold uppercase letters (e.g.~$\mathbf{X}$) for tensors, plain uppercase letters (e.g.~$A$) for matrices, bold lowercase letters (e.g.~$\mathbf{x}$) for vectors, and plain lowercase letters (e.g.~$t,\lambda$) for scalars.
For any positive integer $p$, let $[p] = \{1, \ldots, p\}$. 
For real numbers $a,b$, write $a\vee b=\max\{a,b\}$ and $a\wedge b=\min\{a,b\}$.
For a matrix $A\in\R^{p_1\times p_2}$, denote by $\|A\|$ its spectral norm, $\|A\|_{\mathrm{F}}$ its Frobenius norm, and $\lambda_1(A) \geq \dots \geq  \lambda_{p_1 \wedge p_2} (A)\geq 0$ its singular values. 
For $j\in [p_1\wedge p_2]$, let $\mathrm{SVD}_{j}(A)$ denote the matrix whose columns are the leading $j$ left singular vectors of $A$;
we fix a measurable selection rule whenever there is a tie at the singular values.
For symmetric matrices $A,B$ of the same dimension,
write $A\preceq B$ if $B-A$ is positive semidefinite.
For integers $p\geq r$, define the Stiefel manifold $\mathbb{O}_{p\times r} = \{U \in \mathbb{R}^{p\times r}:U^\top U=I_r\}$. For a matrix $U\in \mathbb{O}_{p\times r}$, let $\Pi(U) = UU^\top$ denote the orthogonal projector onto its column space. For any order-$q$ tensor $\mathbf{M} \in \R^{p_1 \times \cdots \times p_q}$ and $s \in [q]$, let $\mathcal{M}_s(\mathbf{M}) \in \mathbb{R}^{p_s\times \prod_{s'\neq s} p_{s'}}$ denote its mode-$s$ matricization. 
For $B\in \mathbb{R}^{m\times p_s}$, the mode-$s$ product of
$\mathbf{M}$ and $B$ is denoted by $\mathbf{M} \times_s B$. We write $\otimes_{j=1}^\ell A_j$ for the Kronecker product of the matrices $A_1, \dots, A_\ell$. 
For probability measures $\mathscr{P},\mathscr{Q}$,
$\mathrm{TV}(\mathscr{P},\mathscr{Q})$ denotes the total variation distance.

The rest of this paper is organized as follows.
Section~\ref{sec:methodology} introduces the federated model,
formalizes the DP constraints, and presents the estimation procedure together with its privacy guarantee.
In Section~\ref{sec:thm_bound}, we state the model assumptions and derive the non-asymptotic upper bound as well as establishing the matching minimax lower bound. 
Section~\ref{sec:adaptive-tuning} develops the private calibration procedure and establishes that the adaptive estimator preserves the privacy guarantee and attains the oracle rate. 
Monte Carlo experiments and a real data application are presented in Section~\ref{sec: numerical}, followed by a discussion of possible extensions in Section~\ref{sec:conclusion}. Proofs and auxiliary results are deferred to the supplementary material.

\section{Methodology}\label{sec:methodology}

In this section, we develop a federated estimator of a common mode-wise loading space from tensor-valued observations and establish its FDP guarantee. 
We first formulate the federated Tucker tensor factor model and define the relevant notions of DP in Section~\ref{subsec:model}. 
The proposed estimator is presented in Section~\ref{subsec:estimator}, 
and its privacy guarantee is established through sensitivity analysis in Section~\ref{subsec:privacy-guarantee}.

\subsection{Federated model and differential privacy} \label{subsec:model}

We consider a federated setting with $m$ local clients. Client~$j\in[m]$ observes $n_j$ order-$K$ tensors, collected
in the dataset $\mathcal{D}_j = \big( \mathbf{X}_1^{(j)}, \ldots, \mathbf{X}_{n_j}^{(j)} \big) 
\in \mathcal{Z}^{n_j}$, where $\mathcal{Z} = \R^{p_1 \times \cdots \times p_K}$. For each $j\in[m]$ and $t\in[n_j]$, we assume that the data follow a Tucker tensor factor model
\begin{equation}\label{eq:fed-model}
\mathbf{X}_t^{(j)} = \mathbf{S}_t^{(j)} + \mathbf{E}_t^{(j)} ,
\quad \mbox{where \ }
\mathbf{S}_t^{(j)} = \mathbf{F}_t^{(j)} \times_1 A_1 \times_2 \cdots \times_K A_K.
\end{equation}
Here, $\mathbf{F}_t^{(j)} \in \R^{r_1 \times \cdots \times r_K}$ is the client-specific core factor tensor, $A_k \in \R^{p_k \times r_k}$ is the mode-$k$ loading matrix, and $\mathbf{E}_t^{(j)}$ is the idiosyncratic noise.
We assume that all clients share the same mode-wise loading spaces, while their latent factors and factor strengths may differ. 
Such an assumption
is a natural one when clients collect comparable tensor-valued observations with aligned modes, as in the multi-hospital, multi-center, and multi-platform applications discussed in Section~\ref{sec: intro}.
The common loading spaces capture shared mode-wise structures, whereas the client-specific latent factors accommodate heterogeneity across clients.

Even when the loading matrices differ across the clients as
$\mathbf{S}_t^{(j)} = \mathbf{F}_t^{(j)} \times_1 A_1^{(j)} \times_2 \cdots \times_K A_K^{(j)}$,
if there exist $A_k$ and $B_k^{(j)} \in \R^{r_k \times r_k}$ satisfying $A_k^{(j)} = A_k B_k^{(j)}$ for all $j \in [m]$, the cross-client heterogeneity can be absorbed into the client-specific core factor, by noting that
$$ \mathbf{S}_t^{(j)} = \Big( \mathbf{F}_t^{(j)} \times_1 B_1^{(j)} \times_2 \cdots \times_K B_K^{(j)} \Big) \times_1 A_1 \times_2 \cdots \times_K A_K. $$
It also implies that $\mathbf{F}_t^{(j)}$ and $A_k, \, k \in [K]$, are identifiable up to an invertible transform only, which is an inherent feature of the Tucker tensor factor models.
On the other hand, the mode-wise loading spaces denoted by $\mathrm{span}(A_k)$, are well-defined, and thus serve as our estimation target.

We require each client's complete release to satisfy record-level DP, with each tensor-valued observation treated as a single record. For a positive integer $N$, two datasets $\mathcal{D}, \mathcal{D}' \in \mathcal{Z}^N$ are neighbors, written as $\mathcal{D} \sim \mathcal{D}'$, if they differ by the replacement of at most one record. The replacement may be arbitrary and need not follow the identical model in~\eqref{eq:fed-model}. We first define central differential privacy (CDP) for a single client.

\begin{definition}[$(\varepsilon,\delta)$-central differential privacy]
\label{def:CDP}
Let $\varepsilon >0$ and $\delta \in [0,1)$. 
A randomized mechanism $\mathscr{Q}: \mathcal{Z}^N \to \mathcal{Y}$ is $(\varepsilon, \delta)$-centrally differentially private if, for every pair $\mathcal{D} \sim \mathcal{D}'$ and every measurable set $S\subseteq \mathcal{Y}$,
\[
\P_{\mathscr{Q}}\{\mathscr{Q}(\mathcal{D}) \in S\}
\leq
e^\varepsilon \P_{\mathscr{Q}}\{\mathscr{Q}(\mathcal{D}') \in S\} + \delta ,
\]
where the probability is taken over the internal randomness of $\mathscr{Q}$.
\end{definition}

We extend this requirement to the federated setting by imposing a separate privacy constraint on each client's dataset.

\begin{definition}[$(\boldsymbol\varepsilon, \boldsymbol\delta)$-federated differential privacy]
\label{def:FDP}
Let $\boldsymbol\varepsilon=(\varepsilon_1,\ldots,\varepsilon_m)^\top
\in (0,\infty)^m$ and
$\boldsymbol\delta=(\delta_1,\ldots,\delta_m)^\top \in [0,1)^m$. A randomized mechanism $\mathscr{M}: \mathcal{Z}^{n_1} \times \cdots \times \mathcal{Z}^{n_m} \to \mathcal{Y}$ is $(\boldsymbol\varepsilon, \boldsymbol\delta)$-federated differentially private if, for every $j\in[m]$, every pair $\mathcal{D}_j \sim \mathcal{D}_j'$, every collection $\mathcal{D}_{-j} = \{\mathcal{D}_\ell: \ell\neq j\}$, and every measurable set $S\subseteq \mathcal{Y}$,
\begin{align*}
\P_{\mathscr{M}}\{\mathscr{M}(\mathcal{D}_j, \mathcal{D}_{-j})\in S\}\leq
e^{\varepsilon_j} \P_{\mathscr{M}}\{\mathscr{M}(\mathcal{D}_j', \mathcal{D}_{-j}) \in S\} + \delta_j,
\end{align*}
where the probability is taken over the internal randomness of $\mathscr{M}$.
\end{definition}

Throughout, we consider noninteractive protocols in which each client $j$ sends the server a single message
$T_j= \mathscr{Q}_j(\mathcal{D}_j)$ computed from its local data. The local mechanisms use mutually independent randomness
and generate their messages without data-dependent feedback from the server or other clients. If each $\mathscr{Q}_j$ satisfies $(\varepsilon_j,\delta_j)$-CDP, then the joint release $(T_1,\ldots,T_m)$ satisfies $(\boldsymbol\varepsilon, \boldsymbol\delta)$-FDP. 
The server constructs its estimate from these messages alone,
so the estimate also satisfies $(\boldsymbol\varepsilon, \boldsymbol\delta)$-FDP by the post-processing property~\citep{DworkRoth2014}.
We focus on estimating $\mathrm{span}(A_k)$ under FDP for a fixed target mode $k\in[K]$. To estimate loading spaces for multiple modes, we allocate each client's privacy budget across modes and apply sequential composition; see Appendix~\ref{subsec: multimode} for a detailed discussion.

\subsection{Federated private loading-space estimation}
\label{subsec:estimator}

\begin{algorithm}[h!t!]
\caption{Federated private estimation of the mode-$k$ loading space}
\label{alg:fed-est}
\begin{algorithmic}
\STATE \textbf{Input:} Datasets $\{\mathcal{D}_j\}_{j=1}^m$; rank $r_k$;
privacy parameters $(\boldsymbol\varepsilon, \boldsymbol\delta)$; block sizes $\{b_j\}_{j=1}^m$; Huber radii $\{R_{j,k}\}_{j=1}^m$;
and weights $\{w_{j,k}\}_{j=1}^m$.
\FOR{$j=1$ to $m$}
    \STATE Set $B_j=\lfloor n_j/b_j \rfloor$. Randomly permute $[n_j]$, independently of the observations, and partition the resulting sequence into blocks $\mathcal{I}_{j,1},\ldots,\mathcal{I}_{j,B_j}$ with the first $B_j-1$ blocks containing $b_j$ indices each and the last block containing the remaining indices.
   \STATE For each $\ell\in[B_j]$, compute
        \[
        \widehat{U}_{j,k}^{[\ell]} :=        \mathrm{SVD}_{r_k}\Bigg\{
        \frac{1}{|\mathcal{I}_{j,\ell}|}
        \sum_{t\in \mathcal{I}_{j,\ell}}
        \mathcal{M}_k(\mathbf{X}_t^{(j)})
        \mathcal{M}_k(\mathbf{X}_t^{(j)})^\top
        \Bigg\}\quad 
        \mbox{and} \quad
        \widehat\Pi_{j,k}^{[\ell]}:= \Pi\big( \widehat{U}_{j,k}^{[\ell]} \big).
        \]
    \STATE Compute
    \[
    q_{j,k} := \frac{1}{B_j}
    \sum_{\ell=1}^{B_j} \sum_{h=1}^{B_j}
    \mathbbm{1}\Big\{
    \|\widehat\Pi_{j,k}^{[\ell]}-\widehat\Pi_{j,k}^{[h]}\|_{\mathrm{F}}
    \leq R_{j,k}/8 \Big\}.
    \]
    \STATE Independently draw $W_{j,k} \sim \mathrm{Lap}(4/\varepsilon_j)$.
    \IF{$q_{j,k} + W_{j,k}
    < (7 B_j)/8 + 2 + (4/\varepsilon_j) \log(2/\delta_j)$}
        \STATE Draw $\widetilde{U}_{j,k}$ independently of
        the data from the Haar distribution on
        $\mathbb{O}_{p_k\times r_k}$.
    \ELSE
        \STATE Form the Huber aggregate
        \[
        \widehat{H}_{j,k} :=
        \underset{H=H^\top}{\mathrm{argmin}}
        \sum_{\ell=1}^{B_j}
        \rho_{R_{j,k}}\big( \|H-\widehat\Pi_{j,k}^{[\ell]}\|_{\mathrm{F}} \big).  
        \tag*{See \eqref{eq:Huber-loss} for the definition of $\rho_{R_{j,k}}$.}
        \]
        \STATE Independently draw $G_{j,k}\in \mathbb{R}^{p_k\times p_k}$ with i.i.d.\ $\mathcal{N}(0,1)$ entries, and set
        \[
        \tau_{j,k}^2 :=
        \frac{64 R_{j,k}^2}{B_j^2 \varepsilon_j^2}
        \log\left( \frac{2.5}{\delta_j} \right) ,
        \quad 
        Z_{j,k} :=
        \frac{\tau_{j,k}}{\sqrt{2}}
        \big(G_{j,k} + G_{j,k}^\top \big) ,
        \quad \mbox{and} \quad
        \widetilde{U}_{j,k} := \mathrm{SVD}_{r_k}(\widehat{H}_{j,k} + Z_{j,k}) .
        \]
    \ENDIF
    \STATE Transmit $\widetilde{U}_{j,k}$ to the  server.
\ENDFOR
\STATE \textbf{Server:}  Output
$\widehat{A}_k^{\mathrm{fed}} := \mathrm{SVD}_{r_k}\big\{ \sum_{j=1}^m w_{j,k} \widetilde{U}_{j,k} \widetilde{U}_{j,k}^\top \big\}$.
\end{algorithmic}
\end{algorithm}

In this subsection, for a fixed mode $k\in[K]$, we construct a federated estimator of the common loading space $\mathrm{span}(A_k)$ under FDP. 
Each client uses blockwise estimation, a private stability test, and Gaussian perturbation of a Huber aggregate to compute a private local loading-space estimate, which it transmits to the server. The server then aggregates the corresponding projection matrices  using client-specific weights. 
The procedure is summarized in Algorithm~\ref{alg:fed-est}.

Specifically, client $j$ randomly partitions its observations into disjoint blocks and computes a loading-space estimate from each block. Blocking ensures that replacing one observation affects at most one block estimate. The client represents these estimates by their orthogonal projectors $\widehat\Pi_{j,k}^{[\ell]}$ and computes a consistency score $q_{j,k}$, 
a normalized count of projector pairs within Frobenius distance $R_{j,k}/8$ of each other. 
A sufficiently high score ensures that the Huber aggregate satisfies the sensitivity bound in Lemma~\ref{lem:Huber-sensitivity} below
under an arbitrary replacement of one observation. 
We refer to the following requirement on the score as the stability condition:
\begin{equation}
\label{eq:sensitivity:condition}
q_{j, k} \ge (7B_j)/8.
\end{equation}
The client performs a private stability test by adding Laplace noise to the consistency score and comparing the perturbed score with the prescribed threshold. The added noise protects the pass/fail decision, which could itself reveal information about individual observations.

If the test passes, the client forms the Huber aggregate $\widehat{H}_{j,k}$ using the Huber loss, defined for $R>0$ by
\begin{equation}\label{eq:Huber-loss}
\rho_R(u) =
\begin{cases}
u^2 /2, & 0\leq u\leq R ,\\
Ru - R^2/2, & u>R .
\end{cases}
\end{equation}
The quadratic part pools nearby block estimates, while the linear part limits the influence of extreme estimates on the aggregate. The client then adds symmetric Gaussian noise to $\widehat{H}_{j,k}$ and extracts the leading $r_k$ left singular vectors of the perturbed matrix to obtain $\widetilde{U}_{j,k}$. The choice of Gaussian noise scale involving $\tau_{j, k}$ is calibrated using the sensitivity bound in Lemma~\ref{lem:Huber-sensitivity}.
If the test fails, the client instead draws $\widetilde{U}_{j,k}$ independently of the data from the Haar distribution on $\mathbb{O}_{p_k\times r_k}$. In either case, it transmits $\widetilde{U}_{j,k}$ to the server. 
The server forms a weighted average of the corresponding projectors $\widetilde{U}_{j,k} \widetilde{U}_{j,k}^\top$ with client-specific weights, and takes its leading $r_k$-dimensional eigenspace as the federated estimate.

Several ingredients of Algorithm~\ref{alg:fed-est} have been studied in the literature on DP estimation.
Specifically, blockwise subspace estimation has been used for private PCA \citep{tsfadia2024differentially}. Our stability test follows the propose-test-release principle of \citet{dwork2009differential}. 
Its pairwise agreement score has the same form as that used by \citet{tsfadia2022friendlycore} for private mean estimation. 
Huber aggregation with Gaussian perturbation was studied for user-level private mean estimation~\citep{zhao2024huber}, 
while weighted aggregation of privatized projectors has been studied for federated PCA \citep{Lietal2024_fdp}. 
Our contribution is to combine a private stability test based on pairwise projector agreement with Huber aggregation.
The resulting mechanism  ensures DP for every pair of neighboring datasets in the full, unbounded sample space, which in turn carries over to the worst-case FDP of the final output, all the while permitting arbitrary-copy replacement.

\begin{remark}
We focus on estimating the mode-wise loading spaces, which are the main objects of interest in many applications. Using the estimated loading matrices $\widehat{A}_1, \ldots, \widehat{A}_K$ provided by the server, client~$j$ can estimate its common component by
\[
\widehat{\mathbf{S}}_t^{(j)} := \mathbf{X}_t^{(j)} \times_1 \Pi(\widehat{A}_1) \times_2 \cdots \times_K \Pi(\widehat{A}_K) , 
\quad \mbox{where } \Pi(\widehat{A}_k):= \widehat{A}_k \widehat{A}_k^\top.
\]
Likewise, the core factors can be estimated locally. Both reconstructions also use raw observations and therefore do not automatically inherit the privacy guarantee of the loading-space estimates. 
If these outputs are retained locally and not released, the FDP guarantee of the federated procedure remains unchanged.
\end{remark}

\subsection{Sensitivity analysis and privacy guarantees}
\label{subsec:privacy-guarantee}

In this subsection, we establish the privacy guarantee for Algorithm~\ref{alg:fed-est}.
A standard approach to releasing statistics under DP is output perturbation, 
i.e., adding noise calibrated to the sensitivity of the released statistics \citep{Dworketal2006, DworkRoth2014}. 

Recall that $\mathcal{D} \sim \mathcal{D}'$ denotes neighboring datasets.
Then for a vector-valued function $f$, we define its $\ell_1$- and $\ell_2$-sensitivities by
\[
\omega_{f,1} := \sup_{\mathcal{D} \sim \mathcal{D}'} \|f(\mathcal{D})-f(\mathcal{D}')\|_1, 
\quad \mbox{and} \quad  
\omega_{f,2} := \sup_{\mathcal{D} \sim \mathcal{D}'} \|f(\mathcal{D})-f(\mathcal{D}')\|_2,
\]
respectively.
For any $\varepsilon>0$,
the Laplace mechanism adds independent $\mathrm{Lap}(\omega_{f,1} /\varepsilon)$ noise to each coordinate, ensuring $(\varepsilon,0)$-CDP \citep[Theorem~3.6]{DworkRoth2014}.
For $\varepsilon, \delta \in (0, 1)$, the Gaussian mechanism adds independent centered Gaussian noise with variance $2\omega_{f,2}^2 \varepsilon^{-2} \log(1.25/\delta)$ to each coordinate, ensuring $(\varepsilon,\delta)$-CDP \citep[Theorem~A.1]{DworkRoth2014}.
For a matrix-valued function $f$, the identity $\|\mathrm{vec}(M)\|_2 =\|M\|_{\mathrm{F}}$ implies that its Frobenius-norm sensitivity,
\[
\omega_{f, \mathrm{F}} := \sup_{\mathcal{D} \sim \mathcal{D}'} \|f(\mathcal{D})-f(\mathcal{D}')\|_{\mathrm{F}},
\]
is equivalent to the $\ell_2$-sensitivity of the vectorized statistic. Thus, Gaussian perturbation can be calibrated directly to the Frobenius-norm sensitivity of matrix-valued statistics.

At each client~$j$, Algorithm~\ref{alg:fed-est} first applies a private stability test to the consistency score $q_{j,k}$.  
If the test fails, the client draws $\widetilde{U}_{j,k}$ independently of the data from the Haar distribution on~$\mathbb{O}_{p_k\times r_k}$; this draw itself satisfies $(0,0)$-CDP. 
If the test passes, the client forms the Huber aggregate $\widehat{H}_{j,k}$ and adds symmetric Gaussian noise before extracting the local loading-space estimate. 
The following lemma provides a uniform sensitivity bound required for the scalar-valued consistency score $q_{j,k}$, and for the matrix-valued Huber aggregate $\widehat{H}_{j,k}$. 

\begin{lemma}\label{lem:Huber-sensitivity}
Let $B\geq 8$, $R>0$, and let $P^{[1]}, \ldots, P^{[B]}$ be a list of symmetric matrices of the same dimension. Define
\begin{align}\label{eq:Huber}
    q := \frac{1}{B} \sum_{\ell=1}^B \sum_{h=1}^B
    \mathbbm{1}\big\{ \|P^{[\ell]}-P^{[h]}\|_{\mathrm{F}} \leq R/8 \big\} 
    \quad \mbox{and}  \quad
    \widehat{H} := \underset{H= H^\top}{\mathrm{argmin}} \sum_{\ell=1}^B \rho_R\big( \|H-P^{[\ell]}\|_{\mathrm{F}} \big),
\end{align}
where $\rho_R$ is the Huber loss defined in \eqref{eq:Huber-loss}. For a second list of symmetric matrices $P^{[1]\prime}, \ldots, P^{[B]\prime}$ differing from the original list in at most one index, define $q'$ and $\widehat{H}'$ analogously. Then
\begin{equation}\label{eq:Huber-score-sensitivity}
|q-q'| \leq 2 .
\end{equation}
Moreover, if $q\geq 7B/8$, then
\begin{equation}\label{eq:Huber-aggregate-sensitivity}
\|\widehat{H} - \widehat{H}'\|_{\mathrm{F}}
\leq \frac{4R}{B} .
\end{equation}
\end{lemma}

Lemma~\ref{lem:Huber-sensitivity} shows that the consistency score $q_{j,k}$ has sensitivity at most two uniformly over all neighboring datasets. Allocating privacy budget $\varepsilon_j/2$ to the private stability test therefore gives $W_{j,k}\sim \mathrm{Lap}(4/\varepsilon_j)$, so the test satisfies $(\varepsilon_j/2,0)$-CDP. 
Under the stability condition in~\eqref{eq:sensitivity:condition}, 
the Huber aggregate $\widehat{H}_{j,k}$ has Frobenius-norm sensitivity at most $4 R_{j,k}/B_j$. We allocate privacy budget of $(\varepsilon_j/2,\delta_j/2)$ to its Gaussian perturbation. The symmetric Gaussian noise in Algorithm~\ref{alg:fed-est} is obtained by symmetrizing a matrix with independent $\mathcal{N}(0, 2 \tau_{j,k}^2)$ entries. The Gaussian calibration therefore gives
\begin{equation}
\label{eq:taujk}
\tau_{j,k}^2 = 
\frac{(4R_{j,k}/B_j)^2}{(\varepsilon_j/2)^2} \log\bigg( \frac{1.25}{\delta_j/2} \bigg) 
= 
\frac{64 R_{j,k}^2}{B_j^2\varepsilon_j^2} \log\bigg( \frac{2.5}{\delta_j} \bigg).
\end{equation}
The remaining $\delta_j/2$ accounts for the possibility that the private stability test passes when the stability condition fails. The threshold is chosen to control the probability of such events. 
Combining this control with the calibrated Laplace and Gaussian perturbations and the data-independent Haar draw when the test fails yields the following privacy guarantee for Algorithm~\ref{alg:fed-est}.

\begin{proposition}\label{prop:standard-FDP}
Fix a mode $k\in[K]$. Suppose $B_j\geq 8$ and
$\varepsilon_j, \delta_j \in(0,1)$ for every $j\in[m]$. Assume that the partition law, block sizes $b_j$, Huber radii $R_{j,k}$, and server weights $w_{j,k}$ are fixed independently of the raw observations.
Then for every $j\in[m]$, the local mechanism in Algorithm~\ref{alg:fed-est} satisfies $(\varepsilon_j,\delta_j)$-CDP. 
Hence the federated estimator $\widehat{A}_k^{\mathrm{fed}}$ also satisfies $(\boldsymbol\varepsilon, \boldsymbol\delta)$-FDP.
\end{proposition}

Proposition~\ref{prop:standard-FDP} establishes a standard worst-case privacy guarantee for every pair of neighboring datasets in the full, unbounded sample space, without distributional assumptions on the observations. In contrast, the privacy guarantees in \citet{cai2024optimal} and \citet{Lietal2024_fdp} hold only with high probability under the sampling model and for independent-copy replacements. 
This is due to that their analyses calibrate the Gaussian mechanism using high-probability sensitivity bounds for spectral projectors computed from sample covariance matrices, which do not provide uniform sensitivity control under arbitrary observation replacement.

\section{Main theoretical results}\label{sec:thm_bound}

Having established the privacy guarantees in Section~\ref{subsec:privacy-guarantee}, we now derive the non-asymptotic upper bound and minimax lower bound for loading-space estimation in Sections~\ref{subsec:upper} and~\ref{subsec:lower}, respectively.
We first introduce the notation and assumptions used throughout this section.

Let $n:= \sum_{j=1}^m n_j$. For the fixed target mode~$k\in[K]$, define
\[
p:= \prod_{\ell=1}^K p_\ell, 
\quad
p_{-k}:= \prod_{\ell\in[K] \setminus\{k\}} p_\ell,
\quad \mbox{and} \quad
r_{-k}:= \prod_{\ell\in[K] \setminus\{k\}} r_\ell,
\]
with the convention $p_{-k}=r_{-k}=1$ when $K=1$. For each client $j\in[m]$, define the empirical mode-$k$ second-moment matrix of the core factors by 
\begin{equation}\label{def-lambda-1}
\Lambda_{j,(k)} := 
\frac{1}{n_j} \sum_{t=1}^{n_j} \mathcal{M}_k(\mathbf{F}_t^{(j)}) \mathcal{M}_k(\mathbf{F}_t^{(j)})^\top.
\end{equation}
Define the largest and smallest eigenvalues of $\Lambda_{j,(k)}$ and the maximum mode-$k$ factor strength over time by
\begin{equation}\label{def-lambda-2}
\lambda_{j,1}^{(k)}:= \|\Lambda_{j,(k)}\|, 
\quad
\lambda_{j,r_k}^{(k)}:= \lambda_{r_k}(\Lambda_{j,(k)}), 
\quad \mbox{and} \quad
\lambda_{j,\max}^{(k)}:= 
\max_{t\in[n_j]} \big\| \mathcal{M}_k(\mathbf{F}_t^{(j)}) \mathcal{M}_k(\mathbf{F}_t^{(j)})^\top\big\|,
\end{equation}
respectively.

We impose the following assumptions on the loading matrices,
latent factors, and idiosyncratic noise.

\begin{assumption}[Loading matrices]\label{ass:loading}
For every $\ell\in[K]$, $A_\ell\in \mathbb{O}_{p_\ell \times r_\ell}$ and $p_\ell > r_\ell$.
\end{assumption}

\begin{assumption}[Latent factors] \label{ass:latent}
The sequences $\{ \mathbf{F}_t^{(j)} \}_{t\in[n_j]}$, $j\in[m]$, are deterministic. 
There exists an absolute constant $C_\lambda >0$ such that 
$0 < \lambda_{j,r_k}^{(k)} \leq \lambda_{j,1}^{(k)} \leq  C_\lambda  \lambda_{j,r_k}^{(k)}$, 
and $\lambda_{j,\max}^{(k)} \leq C_\lambda \lambda_{j,1}^{(k)}$.
\end{assumption}

\begin{assumption}[Idiosyncratic noise]
\label{ass:noise}
The entries of $\{\mathbf{E}_t^{(j)}: j\in[m], t\in[n_j]\}$ are independent and identically distributed as $\mathcal{N}(0, \sigma^2)$.
\end{assumption}

Assumption~\ref{ass:loading} adopts the standard orthonormal normalization \citep{cai2024optimal, ZhangCenCho2026, Babiietal2025},
with loading scales absorbed into the core factors. 
Assumption~\ref{ass:latent} treats the core factors as deterministic,
as in \citet{BaiLi2012, Onatski2012, he2023robust, barigozzi2025tail}. 
The assumption requires comparable factor strengths within each client and prevents a single time point from dominating the empirical second moment of the core factor.  Assumption~\ref{ass:noise} imposes independent Gaussian idiosyncratic noise with a common variance, as in \citet{gu2020generalized, han2022tensor}.
Under Assumptions~\ref{ass:loading}--\ref{ass:noise}, we may define the population mode-$k$ covariance matrix of the observations as
\[
\Sigma_k^{(j)} :=
\E\Bigg\{ \frac{1}{n_j} \sum_{t=1}^{n_j}
\mathcal{M}_k(\mathbf{X}_t^{(j)}) \mathcal{M}_k(\mathbf{X}_t^{(j)})^\top \Bigg\} 
=
A_k \Lambda_{j,(k)} A_k^\top + \sigma^2 p_{-k} I_{p_k} .
\]
Its leading $r_k$-dimensional eigenspace is $\mathrm{span}(A_k)$, and the corresponding eigengap is $\lambda_{j,r_k}^{(k)}$.

\subsection{Non-asymptotic upper bound}\label{subsec:upper}

We first derive an error bound for Algorithm~\ref{alg:fed-est} under the oracle tuning, where the optimal choices of tuning parameters such as Huber radii and client-specific weights depend on unknown quantities.

\begin{theorem}
\label{thm:fed-upper}
Fix a mode~$k \in [K]$. Suppose that Assumptions~\ref{ass:loading}--\ref{ass:noise} hold. 
For each client $j \in [m]$, let
$\varepsilon_j, \delta_j \in (0,1)$.
In Algorithm~\ref{alg:fed-est}, choose $b_j = \left\lceil C_b \{ r_k + \log(n \vee p) \} \right\rceil$
and Huber radii
\begin{align}\label{eq:block-radius}
R_{j,k} =
C_R \sqrt{r_k} \Bigg[ \frac{\sigma}{\sqrt{\lambda_{j,1}^{(k)}}}
\sqrt{\frac{p_k + \log(n\vee p)}{b_j}}
+
\frac{\sigma^2}{\lambda_{j,1}^{(k)}}
\Bigg\{
\sqrt{\frac{p + p_{-k} \log(n\vee p)}{b_j}}
+
\frac{p_k + \log(n\vee p)}{b_j}
\Bigg\}\Bigg],
\end{align}
where $C_b, C_R >0$ are sufficiently large absolute constants.
Let $B_j=\lfloor n_j/b_j \rfloor$ and define the local error rate and aggregation weights by
\begin{equation}\label{eq:Psi-local}
\Psi_{j,k}^{\mathrm{loc}} := 
\frac{R_{j,k}}{\sqrt{r_k}} \sqrt{\frac{\log(n\vee p)}{B_j}} 
+ \tau_{j,k}\sqrt{p_k+\log(n\vee p)},
\quad \mbox{and} \quad 
w_{j,k} =
\frac{(\Psi_{j,k}^{\mathrm{loc}})^{-2}}{\sum_{\ell=1}^m (\Psi_{\ell,k}^{\mathrm{loc}})^{-2}},
\end{equation} 
respectively, where $\tau_{j,k}$ is defined in~\eqref{eq:taujk}. 
Suppose that, for every $j\in[m]$,
\begin{equation}\label{eq:upper-conditions-1}
\lambda_{j,1}^{(k)} / \sigma^{2} \geq C_0 ( p_k + \sqrt{p}),
\end{equation}
and
\begin{equation}\label{eq:upper-conditions-2}
n_j\varepsilon_j \geq C_0 \{r_k + \log(n\vee p)\} \big[ \log\{(n\vee p) /\delta_j\} + \sqrt{r_k\{p_k + \log(n\vee p)\}\{1 + \log(1/\delta_j)\}} \big],
\end{equation}
where $C_0 >0$ is a sufficiently large absolute constant. Then, there exist absolute constants $C,c  >0$ such that, with probability at least $1 -(n\vee p)^{-c}$,
\begin{equation}\label{eq:fed-upper}
\big\| \widehat{A}_k^{\mathrm{fed}} (\widehat{A}_k^{\mathrm{fed}})^\top - A_k A_k^\top \big\|_{\mathrm{F}}^2 
\leq 
\frac{C r_k \log(n\vee p)}{\sum_{j=1}^m (\Psi_{j,k}^{\mathrm{loc}})^{-2}} \wedge (2r_k).
\end{equation}
\end{theorem}

In Theorem~\ref{thm:fed-upper}, condition~\eqref{eq:upper-conditions-1} can be viewed as a signal-to-noise ratio condition, while \eqref{eq:upper-conditions-2} imposes a sample-size condition under privacy constraints. 
Similar conditions arise in the analysis of differentially private PCA \citep{cai2024optimal, Lietal2024_fdp}. The Huber radius $R_{j,k}$ is chosen to match the Frobenius-norm error rate of the blockwise projectors. The signal-to-noise ratio condition ensures $R_{j,k}\lesssim 1$, so that the blockwise projectors concentrate within a sufficiently small neighborhood of the true loading-space projector.  
Under this condition, the sample-size condition in \eqref{eq:upper-conditions-2} implies
\begin{equation}\label{eq:upper-conditions}
B_j\gtrsim 
\frac{1}{\varepsilon_j}
\log\{(n\vee p) /\delta_j\}
\quad\mbox{and}\quad
\sqrt{r_k} \Psi_{j,k}^{\mathrm{loc}}
\lesssim 1.
\end{equation}
The first bound provides enough blocks for the private stability test to pass with high probability, while the second controls the accuracy of the privatized local estimators. 

When there exists a single client ($m=1$),  Proposition~\ref{prop:standard-FDP} guarantees $(\varepsilon_1, \delta_1)$-CDP, 
while Theorem~\ref{thm:fed-upper} gives a Frobenius-norm error bound of order $\sqrt{r_k \log(n\vee p)} \Psi_{1,k}^{\mathrm{loc}} \wedge \sqrt{r_k}$. For fixed factor ranks and $p_k\gtrsim \log(n\vee p)$, this bound simplifies in the nontrivial regime to
\begin{equation}\label{eq:Psi-local-order}
\Bigg( \frac{\sigma \sqrt{p_k}}{\sqrt{\lambda_{1,1}^{(k)}}} + \frac{\sigma^2 \sqrt{p}}{\lambda_{1,1}^{(k)}} \Bigg) \frac{1}{\sqrt{n_1}} 
+
\Bigg( \frac{\sigma \sqrt{p_k}}{\sqrt{\lambda_{1,1}^{(k)}}} +\frac{\sigma^2 \sqrt{p}}{\lambda_{1,1}^{(k)}} \Bigg)  \frac{\sqrt{p_k \log(2.5/\delta_1)}}{n_1 \varepsilon_1} , 
\end{equation}
up to logarithmic factors. 
The two terms represent the non-private estimation error and the privacy cost, respectively. 
To compare this rate with existing differentially private PCA results, set $K=1$, so that $p_{-k}=1$. For $r_k=1$, 
the rate in \eqref{eq:Psi-local-order} matches the bounds derived for differentially private PCA methods proposed by \citet{liu2022dppca} and \citet{dungler2026iterative} in their signal-plus-noise setting, up to logarithmic factors. 
For fixed rank $r_k \geq 1$, it also matches the rate of \citet{cai2024optimal} under the spiked covariance model, up to logarithmic factors. 

For general $K\geq 1$, suppose further that $\sigma$ is fixed and that the factors are strong, with
$\lambda_{1,1}^{(k)}\asymp p$. 
The error bound in \eqref{eq:Psi-local-order} then simplifies, up to logarithmic factors, to
\begin{equation}\label{eq:cdp-upper-strong}
\frac{1}{\sqrt{n_1 p_{-k}}} + \frac{\sqrt{p_k \log(2.5/\delta_1)}}{n_1 \varepsilon_1 \sqrt{p_{-k}}}.
\end{equation}
The first term matches the minimax rate for non-private tensor loading-space estimation under fixed ranks and strong factors \citep[Theorem~3.5]{Hanetal2024}. 
For comparison, consider estimating the principal subspace of the vectorized observations in $\mathbb{R}^p$, where $p = p_k p_{-k}$, without exploiting the Tucker structure.
Using the PCA rate of \citet{cai2024optimal} as a benchmark,
the corresponding error bound under fixed ranks, fixed noise variance, and spikes of order~$p$ is
\[
\frac{1}{\sqrt{n_1}} 
+ \frac{\sqrt{p \log(2.5/\delta_1)}}{n_1 \varepsilon_1}
\]
up to logarithmic factors. 
Compared with this benchmark, the two terms in~\eqref{eq:cdp-upper-strong} are smaller by factors of $\sqrt{p_{-k}}$ and $p_{-k}$, respectively, illustrating the benefit of exploiting the Tucker structure.

For multiple clients ($m\geq 2$), the weights in \eqref{eq:Psi-local} minimize $\sum_{j=1}^m w_{j,k}^2 (\Psi_{j,k}^{\mathrm{loc}})^2$
over nonnegative weights summing to one.
The aggregation term in \eqref{eq:fed-upper}, i.e., $1/\sum_{j=1}^m (\Psi_{j,k}^{\mathrm{loc}})^{-2}$,
equals the harmonic mean of the squared local rates divided by $m$
and satisfies
\[
\frac{1}{\sum_{j=1}^m (\Psi_{j,k}^{\mathrm{loc}})^{-2}}
\leq 
\min\bigg\{ 
\min_{j\in[m]}(\Psi_{j,k}^{\mathrm{loc}})^2
,
\frac{1}{m^2} \sum_{j=1}^m (\Psi_{j,k}^{\mathrm{loc}})^2
,
\frac{2}{m} \operatorname*{median}_{j\in[m]}
(\Psi_{j,k}^{\mathrm{loc}})^2
\bigg\}.
\]
Under the conditions of Theorem~\ref{thm:fed-upper},
the federated bound therefore vanishes whenever any single-client bound vanishes. 
It can also vanish as $m$ increases, even when none of the local bounds does. The harmonic aggregation structure is consistent with the federated PCA results of \citet{Lietal2024_fdp}.
For $K = 1$, fixed ranks, and a common factor strength across clients, our bound recovers their rate up to logarithmic factors.

Suppose that the factor ranks and $\sigma$ are fixed, and consider homogeneous clients with $n_j\asymp n/m$,
$\lambda_{j,1}^{(k)}\asymp p$, 
$\varepsilon_j\asymp\varepsilon$, 
and $\delta_j\asymp\delta$.
In the nontrivial regime, the squared error bound in \eqref{eq:fed-upper} is of order
\[
\frac{1}{n p_{-k}} + \frac{m p_k \log(2.5/\delta)}{n^2 \varepsilon^2p_{-k}} ,
\]
up to logarithmic factors. 
The first term is the non-private rate based on all $n$ observations, and the second is the privacy cost. At fixed local sample sizes, aggregation reduces both terms by a factor of $m$,
whereas at fixed total sample size, the privacy cost increases with $m$.

The radii $R_{j,k}$, Gaussian noise scales $\tau_{j,k}$, and aggregation weights $w_{j,k}$ depend on the unknown noise level $\sigma$ and client-specific factor strengths $\lambda_{j,1}^{(k)}$. Under the common covariance model in \citet{Lietal2024_fdp}, the common signal and noise parameters get canceled out due to the normalization in~\eqref{eq:Psi-local} for $w_{j, k}$, 
but they are still needed to calibrate the Gaussian perturbation added for privacy.
In our heterogeneous setting, the normalized weights also depend on unknown client-specific signal parameters.
To address both sources of dependence, we develop in Section~\ref{sec:adaptive-tuning} a private calibration procedure that selects these quantities without knowing $\sigma$ or $\lambda_{j,1}^{(k)}$. Under the conditions stated there, the resulting adaptive estimator preserves the same privacy guarantee and attains the oracle high-probability error rate.

\subsection{Minimax lower bound}\label{subsec:lower}

We next establish a minimax lower bound under FDP to assess the optimality of the upper bound in Theorem~\ref{thm:fed-upper}.

\begin{theorem}\label{thm:fed-lower-target}
Fix a mode $k \in [K]$. 
For $\lambda, \sigma^2 >0$, let $\mathfrak{P}_k^{\mathrm{fed}}(\lambda,\sigma^2)$ denote the class of joint distributions generated by \eqref{eq:fed-model} under Assumptions~\ref{ass:loading}--\ref{ass:noise}, 
with $c\lambda I_{r_k} \preceq \Lambda_{j,(k)} \preceq C\lambda I_{r_k}$ for every $j\in[m]$, where $C \geq c >0$ are absolute constants. 
Let $\mathbf{n} :=(n_1,\ldots,n_m)^\top$, and  $\mathcal{U}_{\mathbf{n}, \boldsymbol\varepsilon, \boldsymbol\delta}^{\mathrm{fed},(k)}$ denote the class of all $\mathbb{O}_{p_k\times r_k}$-valued estimators satisfying $(\boldsymbol\varepsilon, \boldsymbol\delta)$-federated differential privacy within the noninteractive framework.
Suppose that the factor ranks are fixed and $p_k\geq 2r_k$. For every $j\in[m]$, assume that $\varepsilon_j\in (0,1)$, $\delta_j\in (0,e^{-1}]$,
\begin{equation*}
n_j r_{-k} \geq r_k,  
\quad \mbox{and} \quad
c_0 \sqrt{n_j} \varepsilon_j^2 
\geq 
r_k (p_k - r_k) \delta_j \log(1/\delta_j),
\end{equation*}
where $c_0 >0$ is a sufficiently small constant. Then there exists a constant $c_1 >0$ such that 
\begin{equation}
\label{eq:fed-lower} 
\inf_{\widehat{A}_k \in \mathcal{U}_{\mathbf{n}, \boldsymbol\varepsilon, \boldsymbol\delta}^{\mathrm{fed},(k)}} 
\sup_{\mathscr{P} \in \mathfrak{P}_k^{\mathrm{fed}}(\lambda,\sigma^2)} 
\mathbb{E}_{\mathscr{P}} \big\| \widehat{A}_k\widehat{A}_k^\top -A_kA_k^\top \big\|_{\mathrm{F}}^2 
\geq
c_1 \Bigg\{ 
\frac{r_k}{ \sum_{j=1}^m (\underline\Psi_{j,k}^{\mathrm{loc}})^{-2}} 
\wedge (2r_k) \Bigg\}, 
\end{equation} 
where
\begin{equation}\label{eq:Psi-lower} (\underline\Psi_{j,k}^{\mathrm{loc}})^2 := 
\frac{\sigma^2}{\lambda} 
\bigg( \frac{p_k}{n_j} 
+ \frac{p_k^2 r_k}{n_j^2\varepsilon_j^2} \bigg) , 
\quad j\in[m]. 
\end{equation} 
\end{theorem}

The lower bound in \eqref{eq:fed-lower} has the same harmonic aggregation structure as the upper bound in~\eqref{eq:fed-upper}. 
For fixed factor ranks and $\sigma$, with strong factors $\lambda\asymp p$, the local lower rates $\underline\Psi_{j,k}^{\mathrm{loc}}$ in \eqref{eq:Psi-lower} match the corresponding upper rates $\Psi_{j,k}^{\mathrm{loc}}$ in \eqref{eq:Psi-local}, up to logarithmic factors, and 
thus Algorithm~\ref{alg:fed-est} attains the near-optimal rate for estimating the common loading space under FDP in this regime. 
When $K=1$, the lower bound matches that derived in \citet{cai2024optimal} under CDP when $m=1$, and \citet{Lietal2024_fdp} under FDP when $m\geq 2$.
For $K>1$, both the non-private and privacy terms of the squared local rates contain the factor $p_{-k}^{-1}$, reflecting the contribution of the dimensions of the other modes.

\section{Adaptive tuning}
\label{sec:adaptive-tuning}

The oracle choices of the Huber radii $R_{j,k}$, Gaussian noise scales $\tau_{j,k}$ and aggregation weights $w_{j,k}$ in Theorem~\ref{thm:fed-upper} depend on the unknown noise level $\sigma$ and client-specific factor strengths $\lambda_{j,1}^{(k)}$. 
To address this, we develop a private calibration procedure that selects these quantities from the data. 
Section~\ref{subsec:estimator-adaptive} presents this procedure and the resulting adaptive version of Algorithm~\ref{alg:fed-est},
and Section~\ref{subsec:theory-adaptive} establishes its FDP guarantee and estimation error bound.

\subsection{Adaptive federated estimation}\label{subsec:estimator-adaptive}

\begin{algorithm}[!htbp]
\caption{Adaptive federated private loading-space estimation}
\label{alg:adaptive-fed-est}
\small
\begin{algorithmic}
\STATE \textbf{Input:} 
Datasets $\{\mathcal{D}_j\}_{j=1}^m$;
rank $r_k$; 
privacy parameters $(\boldsymbol\varepsilon,\boldsymbol\delta)$;
block sizes $\{b_j\}_{j=1}^m$; 
public grid $\mathcal{R}_k=\{d_s\}_{s=0}^{S_k}$;
constants $C_A,c_0$.
\FOR{$j=1$ to $m$}
    \STATE Randomly split $\mathcal{D}_j$, independently of the observations, into calibration and main samples of sizes $n_j^{(0)}=\lfloor n_j/2 \rfloor$ and $n_j^{(1)} = n_j-n_j^{(0)}$
    \STATE Form $B_j^{(a)}= \lfloor n_j^{(a)} /b_j \rfloor$ blocks in each sample, $a\in \{0,1\}$, as in Algorithm~\ref{alg:fed-est}. 
    \STATE For each calibration block $\ell \in [B_j^{(0)}]$, compute $\widehat{U}_{j,k}^{(0),[\ell]}$ as in Algorithm~\ref{alg:fed-est} and set $\widehat\Pi_{j,k}^{(0),[\ell]} = \Pi(\widehat{U}_{j,k}^{(0),[\ell]})$.
    \STATE For $s=0,\ldots,S_k$, compute
    \begin{equation*}
    q_{j,k}^{(0)}(d_s) := 
    \frac{1}{B_j^{(0)}} \sum_{\ell=1}^{B_j^{(0)}} \sum_{h=1}^{B_j^{(0)}}
    \mathbbm{1}\big\{ \|\widehat\Pi_{j,k}^{(0),[\ell]}
    - \widehat\Pi_{j,k}^{(0),[h]}\|_{\mathrm{F}}
    \leq d_s \big\}.
    \end{equation*}
    \STATE Independently draw $Y_j\sim \mathrm{Lap}(4/\varepsilon_j)$ and $Y_{j,s}\sim \mathrm{Lap}(8/\varepsilon_j)$, $s=0,\ldots,S_k$. 
    Set
    \begin{equation}
    \label{eq:adaptive-pilot-index}
    \widehat{s}_{j,k} := 
    \min\big\{s: q_{j,k}^{(0)}(d_s) + Y_{j,s} \geq (7 B_j^{(0)})/8 + Y_j \big\},  
    \end{equation}
    taking $\widehat{s}_{j,k}=S_k$ if the set is empty.
    \STATE Compute  
    \begin{equation}
    \label{eq:adaptive-radius}
    \widehat{R}_{j,k} := C_A d_{\widehat{s}_{j,k}} , 
    \quad 
    \widehat\tau_{j,k}^2 := \frac{64 \widehat{R}_{j,k}^2}{(B_j^{(1)})^2 \varepsilon_j^2} \log\bigg( \frac{2.5}{\delta_j} \bigg) 
    \quad \mbox{and} \quad \widehat\Psi_{j,k} := \frac{\widehat{R}_{j,k}}{\sqrt{r_k}} \sqrt{\frac{\log(n\vee p)}{B_j^{(1)}}}
    + \widehat\tau_{j,k}\sqrt{p_k}.
    \end{equation}
    \STATE Set $\widehat\iota_{j,k}=1$ if 
    \begin{equation}
    \label{eq:adaptive-screen}
    \widehat{R}_{j,k}\leq c_0 , 
    \quad 
    \sqrt{r_k} \widehat\Psi_{j,k} \leq c_0, 
    \quad\text{and}\quad 
    B_j^{(1)} \geq 16 + \frac{32}{\varepsilon_j} \log\bigg( \frac{8 m (n\vee p)^{c_0}}{\delta_j} \bigg),
    \end{equation}
    and set $\widehat\iota_{j,k}=0$ otherwise.
    \IF{$\widehat\iota_{j,k}=1$}
        \STATE Using the main-sample blocks formed above,
        apply the local estimation and release steps of
        Algorithm~\ref{alg:fed-est}, with Huber radius
        $\widehat{R}_{j,k}$ and privacy parameters
        $(\varepsilon_j,\delta_j)$, to obtain $\widetilde{U}_{j,k}$.
    \ELSE
        \STATE Draw $\widetilde{U}_{j,k}$ independently of the data from the Haar distribution on $\mathbb{O}_{p_k\times r_k}$.
    \ENDIF
    \STATE Transmit  $(\widehat\iota_{j,k}, \widehat\Psi_{j,k}, \widetilde{U}_{j,k})$ to the server.
\ENDFOR
\STATE \textbf{Server:} 
Set $\widehat{\mathcal{G}}_k := \{ j\in[m]: \widehat\iota_{j,k}=1 \}$.
\IF{$\widehat{\mathcal{G}}_k \neq \emptyset$}
    \STATE Set
    $\widehat{w}_{j,k} =
    \widehat\Psi_{j,k}^{-2} /
    \sum_{\ell\in \widehat{\mathcal{G}}_k}
    \widehat\Psi_{\ell,k}^{-2}$
    for $j\in \widehat{\mathcal{G}}_k$, and output
    \[
    \widehat{A}_k^{\mathrm{ad}} :=
    \mathrm{SVD}_{r_k}\bigg\{
    \sum_{j\in \widehat{\mathcal{G}}_k}
    \widehat{w}_{j,k}
    \widetilde{U}_{j,k} \widetilde{U}_{j,k}^\top
    \bigg\}.
    \]
\ELSE
    \STATE Output $\widehat{A}_k^{\mathrm{ad}}$ drawn independently of the data from the Haar distribution on $\mathbb{O}_{p_k\times r_k}$.
\ENDIF
\end{algorithmic}
\end{algorithm}

In this subsection, we develop an adaptive version of Algorithm~\ref{alg:fed-est}. 
In short, each client selects the Huber radius through a private threshold search using blockwise loading-space projectors from a calibration sample. 
The selected radius determines the Gaussian noise scale and the estimated local error rate. 
Clients that pass the screening conditions and use these calibrated parameters to compute private local loading-space estimates from separate main samples. 
The server aggregates the corresponding projectors using client-specific weights based on the estimated local error rates. 
The procedure is summarized in Algorithm~\ref{alg:adaptive-fed-est}.

Particularly, note that the oracle radius $R_{j,k}$ in \eqref{eq:block-radius} reflects the scale of  client $j$'s blockwise projector errors $\| \widehat\Pi_{j,k}^{[\ell]} - A_k A_k^\top \|_{\mathrm{F}}$. 
As $A_k A_k^\top$ is unknown, we instead use pairwise distances between the calibration projectors $\widehat\Pi_{j,k}^{(0),[\ell]}$. Under the conditions of Theorem~\ref{thm:adaptive-upper}, the empirical $7/8$ quantile of these distances is comparable to the oracle radius $R_{j,k}$ with high probability.
To approximate this quantile, given a constant $\gamma\geq 1$, let $d_{\min} := (n\vee p)^{-\gamma}$, and we define the public grid
\begin{equation}\label{eq:adaptive-grid}
\mathcal{R}_k := \{d_s=2^s d_{\min}: s=0,\ldots,S_k\},
\quad\text{where }
S_k := \big\lceil
\log_2( d_{\min}^{-1} \sqrt{2r_k}) 
\big\rceil.
\end{equation}
The grid contains $O\{\log(n\vee p)\}$ points, and its upper endpoint $d_{S_k}\geq \sqrt{2r_k}$ bounds all possible Frobenius distances between rank-$r_k$ orthogonal projectors. 
The ratio $q_{j,k}^{(0)}(d_s) /B_j^{(0)}$ is the proportion of the $(B_j^{(0)})^2$ ordered projector pairs whose distance does not exceed $d_s$. 
Therefore, the first grid point satisfying $q_{j,k}^{(0)}(d_s)\geq (7B_j^{(0)})/8$, which is akin to the stability condition in~\eqref{eq:sensitivity:condition}, rounds the empirical $7/8$ quantile up to the nearest grid point.

To protect this data-dependent selection, we use a private threshold search based on the AboveThreshold mechanism \citep[Algorithm~1]{DworkRoth2014}. 
In detail, the client adds independent Laplace noise to each consistency score and to the threshold, reusing the same noisy threshold across all grid points. It then selects the first index satisfying the comparison in \eqref{eq:adaptive-pilot-index}. The selected grid point determines $\widehat{R}_{j,k}$ and $\widehat\tau_{j,k}$ through \eqref{eq:adaptive-radius}, where $C_A >0$ is a sufficiently large absolute constant. 
These quantities in turn determine the estimated local error rate $\widehat\Psi_{j,k}$ in \eqref{eq:adaptive-radius},
constructed analogously to the oracle local error rate $\Psi_{j,k}^{\mathrm{loc}}$ in \eqref{eq:Psi-local}.

The client then uses the calibrated radius $\widehat{R}_{j,k}$, estimated local error rate $\widehat\Psi_{j,k}$, and main-sample block count $B_j^{(1)}$ to check the screening conditions in \eqref{eq:adaptive-screen}, where $c_0 >0$ is a sufficiently small absolute constant. These conditions are computable counterparts of the radius bound and the requirements in~\eqref{eq:upper-conditions}. 
If $\widehat\iota_{j,k}=1$, the client applies the local steps of Algorithm~\ref{alg:fed-est} to its main sample with Huber radius $\widehat{R}_{j,k}$ and Gaussian noise scale $\widehat\tau_{j,k}$. Otherwise, it draws a data-independent random subspace from the Haar distribution.
Each client transmits its screening indicator $\widehat\iota_{j,k}$, estimated local error rate $\widehat\Psi_{j,k}$, and private local loading-space estimate $\widetilde{U}_{j,k}$ to the server.
The server aggregates the corresponding projectors from clients with $\widehat\iota_{j,k}=1$, using weights proportional to $\widehat\Psi_{j,k}^{-2}$.
If all screening indicators are zero, it instead outputs a data-independent random subspace drawn from the Haar distribution.

\subsection{Theoretical guarantees}\label{subsec:theory-adaptive}

We next show that Algorithm~\ref{alg:adaptive-fed-est} satisfies $(\boldsymbol\varepsilon, \boldsymbol\delta)$-FDP and attains the oracle rate in Theorem~\ref{thm:fed-upper}.

\begin{proposition}\label{prop:adaptive-FDP}
Fix a mode $k\in[K]$. Suppose $B_j^{(0)}\geq 1$ and $\varepsilon_j, \delta_j \in(0,1)$ for every $j\in[m]$.
Assume that the split and partition laws, block sizes $b_j$, and public grid $\mathcal{R}_k$ are fixed independently of the raw observations.
Then each complete local release $(\widehat\iota_{j,k}, \widehat\Psi_{j,k}, \widetilde{U}_{j,k})$
produced by Algorithm~\ref{alg:adaptive-fed-est} is $(\varepsilon_j,\delta_j)$-centrally differentially private. Hence, the adaptive estimator $\widehat{A}_k^{\mathrm{ad}}$ satisfies $(\boldsymbol\varepsilon, \boldsymbol\delta)$-federated differential privacy.
\end{proposition}

Proposition~\ref{prop:adaptive-FDP} provides the same worst-case privacy guarantee as Proposition~\ref{prop:standard-FDP}, under arbitrary record replacements and without distributional assumptions on the observations.
For each client~$j\in[m]$, the private threshold search satisfies $(\varepsilon_j,0)$-CDP. 
The tuning parameters $\widehat{R}_{j,k}$ and $\widehat\tau_{j,k}$, screening indicator $\widehat\iota_{j,k}$, and estimated local error rate $\widehat\Psi_{j,k}$ are obtained from its output by post-processing. 
For every fixed calibration output, the main-sample procedure satisfies $(\varepsilon_j,\delta_j)$-CDP.
Since the split is independent of the observations and the two samples are disjoint, parallel composition preserves $(\varepsilon_j,\delta_j)$-CDP for the complete local release, 
even though each stage uses $\varepsilon_j$.
These CDP guarantees for the complete local releases imply $(\boldsymbol\varepsilon, \boldsymbol\delta)$-FDP for their joint release.
The server weights and final estimator inherit this guarantee by post-processing.

We next establish an error bound for the adaptive estimator.

\begin{theorem}\label{thm:adaptive-upper}
Fix a mode $k\in[K]$. Suppose that 
Assumptions~\ref{ass:loading}--\ref{ass:noise} hold.
Consider Algorithm~\ref{alg:adaptive-fed-est} with block sizes $b_j=\lceil C_b \{r_k + \log(n\vee p)\} \rceil$, $j\in[m]$, and the public grid in \eqref{eq:adaptive-grid} for a fixed constant $\gamma\geq 1$. 
Suppose that $p_k\geq 2r_k$, 
$p_k\geq C_p\log(n\vee p)$,
and for every $j\in[m]$, $\varepsilon_j, \delta_j \in(0,1)$,
\begin{equation}\label{eq:adaptive-upper-conditions-1}
C_0(p_k+\sqrt{p}) 
\leq  
\lambda_{j,1}^{(k)}/ \sigma^2 
\leq 
c_0 (n\vee p)^{2\gamma}
\frac{r_k p_k}{r_k + \log(n\vee p)},
\end{equation}
and
\begin{equation}
\label{eq:adaptive-upper-conditions-2}
n_j\varepsilon_j \geq 
C_0\{r_k + \log(n\vee p)\} \big[ \log\{( n\vee p) /\delta_j\}
+
\sqrt{r_k p_k \{1 + \log(1/\delta_j)\}} \big],
\end{equation}
where $C_b, C_p, C_0 >0$ are sufficiently large absolute constants, and $c_0 >0$
is the same constant as in the screening conditions \eqref{eq:adaptive-screen} of Algorithm~\ref{alg:adaptive-fed-est}.
Then there exist absolute constants $C, c >0$ such that, with probability at least $1 - (n\vee p)^{-c}$, we have
\begin{equation*}
\big\|\widehat{A}_k^{\mathrm{ad}}(\widehat{A}_k^{\mathrm{ad}})^\top - A_kA_k^\top \big\|_{\mathrm{F}}^2
\leq
\frac{C r_k\log(n\vee p)}{\sum_{j=1}^m (\Psi_{j,k}^{\mathrm{loc}})^{-2}}
\wedge (2r_k) ,
\end{equation*}
where $\Psi_{j,k}^{\mathrm{loc}}$ is defined in \eqref{eq:Psi-local}.
\end{theorem}

The lower bound of the signal-to-noise ratio in \eqref{eq:adaptive-upper-conditions-1} has the same form as \eqref{eq:upper-conditions-1} in Theorem~\ref{thm:fed-upper}.
Since $p_k\geq C_p \log(n\vee p)$, the sample-size condition in \eqref{eq:adaptive-upper-conditions-2} also has the same order as its oracle counterpart in~\eqref{eq:upper-conditions-2}.
The upper bound in \eqref{eq:adaptive-upper-conditions-1} ensures that the lower endpoint $d_{\min}=(n\vee p)^{-\gamma}$ of the public grid $\mathcal{R}_k$ lies below the oracle radius~$R_{j,k}$ in \eqref{eq:block-radius}.
Under these conditions, the proof of Theorem~\ref{thm:adaptive-upper} shows that, with high probability, the calibrated Huber radii, Gaussian noise scales and aggregation weights satisfy
\[
\widehat{R}_{j,k} \asymp R_{j,k} ,
\quad
\widehat\tau_{j,k} \asymp \tau_{j,k} ,
\quad \mbox{and} \quad 
\widehat{w}_{j,k} \asymp w_{j,k},
\quad 
\forall j\in[m] ,
\]
and all clients satisfy \eqref{eq:adaptive-screen}.
The adaptive estimator therefore attains the oracle rate in Theorem~\ref{thm:fed-upper} without the knowledge of the noise level or client-specific factor strengths.

\section{Numerical analysis}
\label{sec: numerical}

In this section, we investigate the finite-sample performance of the proposed estimators.
We first present simulation results in Section~\ref{subsec:simulation}, followed by an application to the meditation traits electroencephalography (EEG) dataset in Section~\ref{subsec:realdata}.

\subsection{Simulation studies}
\label{subsec:simulation}

We begin by describing the data generating process and the competing methods used throughout.

\medskip
\noindent
\textbf{Data generating process.}
We first consider model~\eqref{eq:fed-model} with order $K=2$, dimensions $p_1=p_2=15$, ranks $r_1=r_2=2$, and $m=3$ clients. 
Each client has $n_j \in \{2000, 4000, 6000, 8000\}$ observations, with equal sample sizes across the clients in each experiment. 
The target is the mode-$1$ loading space.
For each mode~$\ell$, we generate an orthonormal loading matrix $A_\ell \in \mathbb{O}_{p_\ell \times r_\ell}$ from a Gaussian matrix via QR decomposition and use the same loading matrices for all clients.
The entries of the core factors are independent $\mathcal{N}(0,\lambda)$
variables, with $\lambda = \prod_{\ell=1}^K p_\ell$, and the idiosyncratic noise entries are independent $\mathcal{N}(0,1)$ variables. 
Factors and noise
are independent across the observations and clients.
We assign heterogeneous privacy budgets $\boldsymbol\varepsilon= (0.6, 0.8, 1)^\top$
and $\boldsymbol\delta = (0.2, 0.2, 0.2)^\top$, to the three clients.
Unless otherwise specified, every experiment is repeated $100$ times.

\medskip
\noindent
\textbf{Competing methods.}
We compare the following estimators:
\begin{itemize}
    \item [(i)] \textsc{oracle}: 
    Algorithm~\ref{alg:fed-est} with the oracle radius $R_{j,k}$ in \eqref{eq:block-radius} and the oracle weights $w_{j,k}$ in \eqref{eq:Psi-local}, both computed from the true $\sigma$ and $\lambda_{j,1}^{(k)}$.

    \item [(ii)] \textsc{ada}:
    Algorithm~\ref{alg:adaptive-fed-est}, using adaptively chosen Huber Radii, Gaussian noise scale and weights. 
    
    \item [(iii)] \textsc{radius}:
    Algorithm~\ref{alg:fed-est} with the oracle weights $w_{j,k}$ and a randomized radius $R_{j,k} = U_{j,k} R_{j,k}^*$, where $R_{j,k}^*$ is the oracle radius in \eqref{eq:block-radius} and $U_{j,k} \sim \mathrm{Unif}(0.5, 1.5)$ is drawn independently for each client.

    \item [(iv)] \textsc{weight}:
    Algorithm~\ref{alg:fed-est} with the oracle radius $R_{j,k}$ in \eqref{eq:block-radius} and uniform weights $w_{j,k} = 1/m$.

    \item [(v)] \textsc{practical}:
    Algorithm~\ref{alg:fed-est} with the randomized radius as in \textsc{radius} and uniform weights as in \textsc{weight}.
    This method uses no oracle quantity other than the scale of $R_{j,k}^*$.

    \item [(vi)] \textsc{hp-vec}:
    the high-probability privatized algorithm by \cite{Lietal2024_fdp} applied to the vectorized observations $\mathrm{vec}(\mathbf{X}_t^{(j)})$ (detailed below).
\end{itemize}

The first two methods are those proposed in this paper. In particular, \textsc{oracle} with known $\sigma$ and $\lambda_{j,1}^{(k)}$ provides a benchmark, while \textsc{ada} does not require the knowledge of either quantity.
Methods~\textsc{radius}, \textsc{weight}, and \textsc{practical} are variants of \textsc{oracle}.
Note that none of these variants are genuinely feasible in practice, as even \textsc{practical} still relies partially on the oracle radius.
Methods~\textsc{oracle}, \textsc{ada}, \textsc{radius}, \textsc{weight}, and \textsc{practical} all directly consider the tensor-valued observations $\mathbf{X}_t^{(j)}$ and target the mode-$1$ loading space $\mathrm{span}(A_1)$.
In contrast, \textsc{hp-vec} considers the vectorized observations $\mathrm{vec}(\mathbf{X}_t^{(j)}) \in\R^p$, where $p=\prod_{\ell=1}^K p_\ell$, and targets the resulting loading space $\mathrm{span}(A_2\otimes A_1) \subset\R^p$ of the vectorized model.
For all methods, we set the target rank to $r_1$ on tensor-valued observations and to $r = \prod_{\ell=1}^K r_\ell$ on vectorized observations.
Additionally, we define methods~\textsc{oracle-m1} and \textsc{ada-m1} as applying \textsc{oracle} and \textsc{ada} to client~$1$ alone, to assess the benefit of federated aggregation.

\medskip
\noindent
\textbf{Evaluation and results.}
For evaluation, we report the empirical mean and standard deviation of the projector distance of the estimated projectors to their respective targets, i.e., 
$\big\|\widehat{A}_1 \widehat{A}_1^\top - A_1A_1^\top \big\|_{\mathrm{F}}$
for methods~(i)--(iv), and
$\big\|\widehat{A} \widehat{A}^\top - (A_2\otimes A_1) (A_2\otimes A_1)^\top \big\|_{\mathrm{F}}$
for method~(v).
To implement Algorithm~\ref{alg:fed-est},
we fix the block size at $b_j=4$, and $C_R=15$ in~\eqref{eq:block-radius}.
The same parameters are used in Algorithm~\ref{alg:adaptive-fed-est}, and we additionally set
$C_A=8$ and $d_{\min}=(n\vee p)^{-2}$ in \eqref{eq:adaptive-grid}, where $n=\sum_j n_j$ and $p=\prod_\ell p_\ell$.
We perform additional sensitivity analysis on $b_j$ and $C_R$ in Section~\ref{app:add:num} of the supplementary material.

\begin{table}[t]
\centering
\caption{Mean Frobenius projector distance over $100$ realizations for $K=2$, $p_\ell=15$, $r_\ell=2$, $m=3$, $\boldsymbol\varepsilon=(0.6,0.8,1)^\top$, and $\boldsymbol\delta=(0.2,0.2,0.2)^\top$, with standard deviations in parentheses.}
\label{tab:sim-K2-small-eps2-main}
\ra{1.3}
\begin{tabular}{lllll}
\toprule
Method & $n_j=2000$ & $n_j=4000$ & $n_j=6000$ & $n_j=8000$ \\
\midrule
\textsc{oracle}       & 1.410 (0.498) & 0.382 (0.341) & 0.177 (0.103) & 0.108 (0.015) \\
\textsc{ada}          & 1.347 (0.246) & 0.533 (0.075) & 0.313 (0.045) & 0.255 (0.037) \\
\midrule
\textsc{radius}       & 1.233 (0.489) & 0.913 (0.569) & 0.900 (0.605) & 0.799 (0.620) \\
\textsc{weight}       & 1.511 (0.426) & 0.425 (0.330) & 0.169 (0.070) & 0.117 (0.030) \\
\textsc{practical}    & 1.258 (0.438) & 0.927 (0.553) & 0.875 (0.587) & 0.790 (0.591) \\
\textsc{hp-vec}       & 0.442 (0.011) & 0.229 (0.005) & 0.157 (0.004) & 0.120 (0.003) \\
\midrule
\textsc{oracle-m1}    & 1.806 (0.216) & 1.501 (0.620) & 0.937 (0.767) & 0.437 (0.515) \\
\textsc{ada-m1}       & 1.863 (0.090) & 1.260 (0.243) & 0.802 (0.187) & 0.637 (0.098) \\
\bottomrule
\end{tabular}
\end{table}

Table~\ref{tab:sim-K2-small-eps2-main} reports the results for $K=2$. The method \textsc{weight} performs comparably to \textsc{oracle}, suggesting that uniform weighting incurs little loss in this setting with equal sample sizes and homogeneous factor strengths. Thus, we consider heterogeneous factor strengths at the end of this subsection to assess the benefit of the proposed aggregation weights. By contrast, \textsc{radius} and \textsc{practical} have larger mean errors and greater variability than \textsc{oracle} for $n_j\geq4000$, illustrating the sensitivity to the choice of radius. The reported errors of \textsc{hp-vec} are smaller than those of \textsc{oracle} at the first three sample sizes, but slightly larger at $n_j=8000$.
This may be accounted for by its weaker privacy requirement: \textsc{hp-vec} directly perturbs local empirical projectors with Gaussian noise and provides a high-probability privacy guarantee under independent-copy replacements.
Our procedure additionally uses Huber aggregation and private stability testing, with a data-independent fallback when the test fails. These safeguards ensure worst-case privacy under arbitrary record replacements, but can incur additional estimation error when local sample sizes are small. Moreover, \textsc{hp-vec} concerns the full vectorized loading space rather than the mode-$1$ loading space.

Among the procedures compared here, only \textsc{ada} requires no knowledge of the noise level or factor strengths for either local estimation or server aggregation. Nevertheless, \textsc{ada} performs reasonably well relative to \textsc{oracle} across the reported sample sizes, illustrating the effectiveness of private calibration. Both federated estimators also have smaller mean errors than their corresponding client-$1$ versions at every sample size, consistent with the aggregation benefit described in Theorem~\ref{thm:fed-upper}.

We next repeat the experiment with order $K=3$, dimensions $p_\ell=10$, and ranks $r_\ell=2$ for every mode. 
All other numerical choices are unchanged. Table~\ref{tab:sim-K3-small-eps2-main} shows that both \textsc{oracle} and \textsc{ada} outperform most of the other methods overall, with \textsc{oracle} attaining the smallest mean error at every sample size.  In particular, \textsc{hp-vec} has substantially larger reported errors than both of the proposed estimators. This difference is consistent with the benefit of exploiting the multiway structure: vectorization leads to a higher-dimensional estimation problem, whereas our procedures retain the tensor structure and directly estimate the mode-wise loading spaces.

\begin{table}[t]
\centering
\caption{Mean Frobenius projector distance over $100$ realizations for $K=3$, $p_\ell=10$, $r_\ell=2$, $m=3$, $\boldsymbol\varepsilon=(0.6,0.8,1)^\top$, and $\boldsymbol\delta=(0.2,0.2,0.2)^\top$, with standard deviations in parentheses.}
\label{tab:sim-K3-small-eps2-main}
\ra{1.3}
\begin{tabular}{lllll}
\toprule
Method & $n_j=2000$ & $n_j=4000$ & $n_j=6000$ & $n_j=8000$ \\
\midrule
\textsc{oracle}       & 0.086 (0.014) & 0.046 (0.008) & 0.031 (0.005) & 0.023 (0.004) \\
\textsc{ada}          & 0.266 (0.048) & 0.097 (0.018) & 0.058 (0.009) & 0.048 (0.008) \\
\midrule
\textsc{radius}       & 0.650 (0.588) & 0.591 (0.566) & 0.419 (0.453) & 0.501 (0.556) \\
\textsc{weight}       & 0.096 (0.018) & 0.049 (0.009) & 0.032 (0.006) & 0.025 (0.005) \\
\textsc{practical}    & 0.660 (0.576) & 0.592 (0.549) & 0.415 (0.439) & 0.500 (0.545) \\
\textsc{hp-vec}       & 1.986 (0.027) & 0.923 (0.008) & 0.614 (0.005) & 0.464 (0.003) \\
\midrule
\textsc{oracle-m1}    & 0.219 (0.166) & 0.100 (0.017) & 0.069 (0.011) & 0.053 (0.009) \\
\textsc{ada-m1}       & 1.794 (0.117) & 0.254 (0.049) & 0.150 (0.026) & 0.122 (0.021) \\
\bottomrule
\end{tabular}
\end{table}

\medskip
\noindent
\textbf{Heterogeneous factor strengths.}
We next examine how heterogeneity in factor strength affects estimation and aggregation across clients. We consider the setting with $K=2$, $p_1=p_2=15$, and $r_1=r_2=2$. There are three clients, each with sample size $n_j=12000$ and privacy parameters $(\varepsilon_j,\delta_j)=(0.8,0.2)$. The entries of the core factor at client~$j$ have variance $\lambda_j$.
With $p=p_1p_2=225$, we consider four configurations:
\begin{align*}
\mathrm{S1}: &\quad (\lambda_1,\lambda_2,\lambda_3)
= (p, p^{0.6}, p^{0.6}),
\quad
\mathrm{S2}: \quad (\lambda_1,\lambda_2,\lambda_3)
= (p^{0.6}, p, p) , 
\\
\mathrm{S3}: &\quad (\lambda_1,\lambda_2,\lambda_3)
= (p^{0.6}, p^{0.6}, p^{0.6}) ,
\quad
\mathrm{S4}: \quad (\lambda_1,\lambda_2,\lambda_3)
= (p,p,p).
\end{align*}
S1 and S2 represent heterogeneous configurations with one and two strong-factor clients, respectively, while S3 and S4 consider homogeneous weak- and strong-factor clients, respectively. All other implementation choices remain unchanged.

\begin{table}[t]
\centering
\caption{Mean Frobenius projector distance over $100$ realizations under factor-strength configurations~S1--S4 for $K=2$, $p_\ell=15$, $r_\ell=2$, $m=3$, $n_j=12000$, and $(\varepsilon_j,\delta_j)=(0.8,0.2)$, with standard deviations in parentheses.}
\label{tab:heterogeneous-strength}
\ra{1.3}
\begin{tabular}{lcccc}
\toprule
Method & S1 & S2 & S3 & S4\\
\midrule
\textsc{oracle}
& 0.126 (0.019) & 0.090 (0.014)
& 0.402 (0.067) & 0.074 (0.011)\\
\textsc{ada}
& 0.257 (0.032) & 0.179 (0.024)
& 1.864 (0.087) & 0.146 (0.020)\\
\midrule
\textsc{radius}
& 0.969 (0.836) & 0.935 (0.628)
& 0.720 (0.331) & 0.799 (0.578)\\
\textsc{weight}
& 0.312 (0.051) & 0.228 (0.035)
& 0.396 (0.067) & 0.073 (0.011)\\
\textsc{practical}
& 0.742 (0.401) & 0.811 (0.514)
& 0.720 (0.331) & 0.799 (0.578)\\
\textsc{hp-vec}
& 0.254 (0.006) & 0.187 (0.004)
& 0.308 (0.007) & 0.077 (0.002)\\
\midrule
\textsc{oracle-m1}
& 0.123 (0.017) & 0.629 (0.095)
& 0.623 (0.093) & 0.125 (0.017)\\
\textsc{ada-m1}
& 0.290 (0.043) & 1.858 (0.091)
& 1.858 (0.091) & 0.290 (0.043)\\
\bottomrule
\end{tabular}
\end{table}

Table~\ref{tab:heterogeneous-strength} reports the results for all the aforementioned methods. The results illustrate the importance of accounting for heterogeneous factor strengths when choosing aggregation weights. In the heterogeneous settings S1 and S2, \textsc{oracle} achieves substantially smaller errors than \textsc{weight}, which assigns equal weights to all clients,
whereas the two methods perform similarly in S3 and S4. In S1, where client~1 is the only client with strong factors, \textsc{oracle} and \textsc{oracle-m1} have similar mean errors. In S2, where client~1 has weak factors and the other two clients have strong factors, \textsc{oracle} substantially improves on \textsc{oracle-m1}. In both heterogeneous settings, \textsc{ada} also achieves smaller mean errors than \textsc{weight}, \textsc{radius}, and \textsc{practical}, despite requiring no knowledge of the noise level or factor strengths.
The large mean error of \textsc{ada} in S3 highlights the difficulty of adaptive estimation when all clients have weak factors.
Although the procedure does not require factor strengths as inputs, weaker factors can increase the dispersion of the blockwise calibration projectors, leading to larger estimated radii $\widehat{R}_{j,k}$ and local error rates $\widehat\Psi_{j,k}$ through \eqref{eq:adaptive-radius}.
Larger values make it harder to meet the screening conditions in \eqref{eq:adaptive-screen}; if no client passes screening, the server returns a Haar random subspace, which can lead to a large estimation error.

\subsection{Real data analysis}
\label{subsec:realdata}

We analyze the meditation traits dataset studied by \citet{braboszcz2017increased}. 
The dataset consists of tensor-valued EEG observations with three anatomical regions, six frequency bands, and three time frames per window,
giving a tensor $\mathbf{X}_t^{(j)} \in \mathbb{R}^{3 \times 6 \times 3}$ for each observation.
In total, it contains $13740$ windows from $29$ participants across four groups:
a control group and three practitioner groups (Vipassana, Himalayan Yoga, and Isha Shoonya).
We first demean each feature within each participant, and then form $m = 4$ clients by grouping the windows according to the group membership, giving local sample sizes $(n_1, n_2, n_3, n_4) = (2380, 3880, 3780, 3700)$.
The target parameter is the mode-$1$ loading space, corresponding to the anatomical regions. We fix the factor ranks at $(r_1,r_2,r_3) = (2,2,2)$, and set $\delta_j = 0.2$ for all clients.

We compare three estimators: \textsc{ada}, \textsc{plug}, and \textsc{hp-vec}. Specifically, \textsc{ada} applies Algorithm~\ref{alg:adaptive-fed-est}. Also for each $j\in[4]$, \textsc{ada-m$j$} applies the same adaptive procedure to client~$j$ alone, using its local sample and privacy parameters.
These single-client versions allow us to assess the benefit of federated aggregation. Second, we introduce \textsc{plug} which applies Algorithm~\ref{alg:fed-est} with the following plug-in estimates. For each client~$j$, let $\widehat\Sigma_k^{(j)}$ denote its sample mode-$k$ covariance matrix, with eigenvalues $\widehat{e}_{j,1} \geq \dots \geq \widehat{e}_{j,p_k}$.
Under the tensor factor model, the population covariance satisfies $\Sigma_k^{(j)} = A_k \Lambda_{j,(k)} A_k^\top + \sigma^2 p_{-k} I_{p_k}$. Its smallest $p_k-r_k$ eigenvalues therefore equal $\sigma^2p_{-k}$, so we set
\[
\widehat\sigma_{j,0}^2 := \operatorname{median}\big\{ \widehat{e}_{j,r_k+1}, \dots, \widehat{e}_{j,p_k} \big\}, \quad
\widehat\sigma_j := \sqrt{\widehat\sigma_{j,0}^2 / p_{-k}}, \quad \mbox{and} \quad
\widehat\lambda_{j,1}^{(k)} :=  \widehat{e}_{j,1} - \widehat\sigma_{j,0}^2 .
\]
The plug-in radius $\widehat{R}_{j,k}$ is obtained by replacing $\sigma$ and $\lambda_{j,1}^{(k)}$ in \eqref{eq:block-radius} with $\widehat\sigma_j$ and $\widehat\lambda_{j,1}^{(k)}$, respectively. The remaining tuning parameters and aggregation weights are computed using the same plug-in estimates. Since these plug-in estimates are not privatized, \textsc{plug} is not covered by the privacy guarantee in Proposition~\ref{prop:standard-FDP} and is therefore included only as an empirical benchmark.
Lastly, \textsc{hp-vec} applies the algorithm of \citet{Lietal2024_fdp} to the vectorized observations $\mathrm{vec}(\mathbf{X}_t^{(j)})$ with target rank $r=8$.
Its plug-in estimates are computed similarly to those of \textsc{plug}, but using the sample covariance matrix of the vectorized data. The original high-probability privacy guarantee does not directly cover this plug-in implementation.

We evaluate all tensor-based estimators, including the single-client versions, against the same non-private mode-1 reference constructed from all four clients. For \textsc{hp-vec}, we use a separate non-private reference for the full vectorized loading space. Both references are defined in Section~\ref{app:real} of the supplementary material.
Table~\ref{tab:Meditation-distances} reports the corresponding  Frobenius projector distance, which measures the degree of agreement with empirical loading spaces rather than error relative to unknown population loading spaces. The method \textsc{ada} has smaller distances than \textsc{plug} in every reported configuration and achieves the smallest distance among the tensor-based estimators in several configurations.
These results support the effectiveness of private adaptive calibration compared with direct plug-in tuning in this application.

We note that the performance varies across block sizes and privacy budgets. In particular, no method achieves the smallest distance in every configuration. The federated \textsc{ada} estimator does not uniformly improve on the best single-client result as it occasionally reports somewhat larger distances with $b_j=6$ or $8$.
This may partly reflect the trade-off associated with block size: increasing the block size provides more observations for each blockwise estimate but reduces the number of blocks available for calibration, testing, and local aggregation. 
This trade-off, together with randomness of the algorithm, can affect the final estimate.

To assess variability due to algorithmic randomness, we repeat the analysis $100$ times on the same dataset.
Specifically, we consider block sizes $b_j\in\{4,6,8\}$ and privacy parameters $\varepsilon_j= \varepsilon \in \{0.4, 0.6, 0.8, 1, 1.2, 1.5, 2\}$, with $\delta_j=0.2$ for all clients.
The block size applies only to the tensor methods.
We regenerate the random partitions and the randomness used in calibration, testing, and release in each run, while holding the dataset and deterministic tuning choices fixed.

\begin{table}[t]
\centering
\caption{Frobenius projector distance to the corresponding non-private reference on the meditation traits dataset for $K=3$, $(p_1,p_2,p_3)=(3,6,3)$, $r_\ell=2$, $m=4$, $n_j\in\{2380,3880,3780,3700\}$, $\varepsilon_j=\varepsilon$, and $\delta_j=0.2$.}
\label{tab:Meditation-distances}
\ra{1.3}
\begin{tabular}{cc ccccccc}
\toprule
$b_j$ & $\varepsilon$ & \textsc{ada} & \textsc{ada-m1} & \textsc{ada-m2} & \textsc{ada-m3} & \textsc{ada-m4} & \textsc{plug} & \textsc{hp-vec} \\
\cmidrule(lr){1-2} \cmidrule(lr){3-9}
4 & 0.6 & 0.228 & 1.398 & 0.303 & 0.342 & 0.250 & 0.785 & 0.750 \\
4 & 0.8 & 0.335 & 1.109 & 1.293 & 0.255 & 0.258 & 1.195 & 0.720 \\
4 & 1.0 & 0.147 & 0.545 & 1.382 & 0.194 & 0.261 & 1.070 & 0.694 \\
\cmidrule(lr){1-2} \cmidrule(lr){3-9}
6 & 0.6 & 0.391 & 1.055 & 0.503 & 0.680 & 0.842 & 1.311 & 0.746 \\
6 & 0.8 & 1.308 & 1.414 & 0.524 & 0.252 & 0.231 & 1.390 & 0.725 \\
6 & 1.0 & 0.853 & 0.693 & 0.226 & 0.098 & 0.995 & 1.402 & 0.684 \\
\cmidrule(lr){1-2} \cmidrule(lr){3-9}
8 & 0.6 & 0.396 & 0.544 & 1.391 & 0.226 & 0.865 & 1.214 & 0.762 \\
8 & 0.8 & 0.538 & 1.321 & 1.004 & 1.088 & 0.843 & 1.414 & 0.725 \\
8 & 1.0 & 0.837 & 1.386 & 0.430 & 0.420 & 0.833 & 1.262 & 0.669 \\
\bottomrule
\end{tabular}
\end{table}

\begin{figure}[t]
    \centering
    \includegraphics[width=\textwidth]{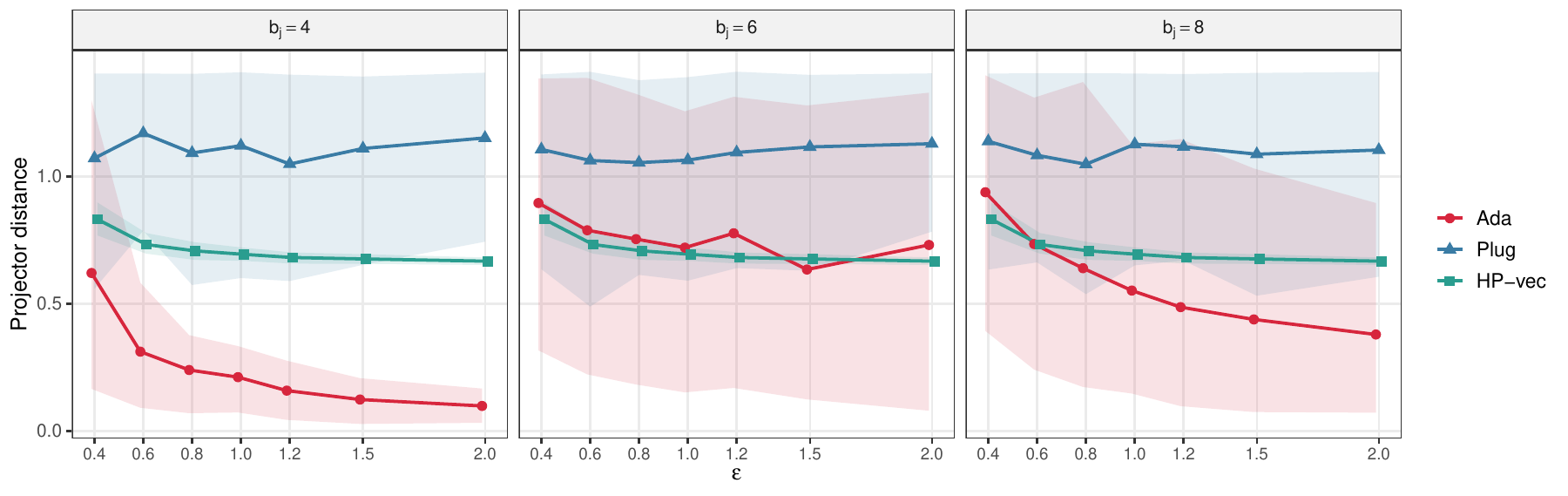}
    \caption{Mean Frobenius projector distance to the corresponding non-private reference over $100$ randomized runs on the meditation traits dataset for $K=3$, $(p_1,p_2,p_3)=(3,6,3)$, $r_\ell=2$, $m=4$, $n_j\in\{2380,3880,3780,3700\}$, $\varepsilon_j=\varepsilon$, and $\delta_j=0.2$, with shaded bands indicating the 10th--90th percentiles. Columns correspond to $b_j=4,6,8$.}
    \label{fig:Meditation-privacy-repeat}
\end{figure}

Figure~\ref{fig:Meditation-privacy-repeat} shows that, for $b_j=4$, \textsc{ada} outperforms the other methods across all examined privacy budgets. Its mean distance decreases steadily as $\varepsilon$ increases, and its variability is substantially smaller at larger values of $\varepsilon$. These results show that its favorable performance with this block size persists across repeated randomized runs. We additionally examine the frequency of complete fallback, defined as the event that no client passes the local stability test.
For \textsc{plug}, this event occurs in nearly all runs across the examined settings, which explains the large distances between its outputs and the non-private reference and limited improvement with increasing $\varepsilon$.
Complete fallback is rare for \textsc{ada}. Substantial variability nevertheless remains, particularly at $b_j=6$, suggesting that complete fallback alone does not explain the variation across the runs. The reported mean distances of \textsc{hp-vec} decrease modestly with $\varepsilon$ and then largely level off, with comparatively little variability at larger budgets.

\section{Discussion}\label{sec:conclusion}

We have developed a federated estimator of common loading spaces in Tucker tensor factor models under the worst-case, record-level FDP. 
Our method combines private stability testing, Huber aggregation and Gaussian perturbation at each client, with weighted aggregation at the server. 
We establish its non-asymptotic upper bound governed by a scaled harmonic mean of the squared local error rates, together with the minimax lower bound, and show that the former match the latter up to logarithmic terms under strong factors, fixed ranks and noise variance. 
We further propose a private calibration procedure that selects the Huber radii, Gaussian noise scales and aggregation weights without requiring any knowledge of the noise level or client-specific factor strengths. The resulting adaptive estimator preserves the privacy guarantee and, under the stated conditions, attains the oracle high-probability error rate.

Several extensions remain open for future research.
Firstly, we assume independent Gaussian idiosyncratic noise. Extending our results to heavy-tailed noise may require robust estimation within blocks and new error bounds for the resulting projector estimates. 
While it is possible to relax the independence assumption to allow for temporally dependent noise,
such temporal dependence raises additional privacy concerns in addition to requiring new estimation error bounds.
In particular, our record-level FDP guarantee remains valid under temporal dependence, but information about an observation may still be inferred from other correlated records. 
Addressing such inference calls for suitable privacy notions and the development of methods satisfying them, together with corresponding estimation error bounds.

Secondly, we assume that the factor ranks are known. 
For unknown ranks, \citet{han2022rank} and \citet{ChenLam2024}, 
among others, proposed information criteria or eigenvalue-ratio methods in the non-private setting.
Adapting these methods to our federated procedure would require establishing conditions for consistent rank selection and, more importantly, developing a private rank-selection mechanism to ensure FDP for the complete release.

Lastly, note that all clients are assumed to share the same loading spaces.
We could extend the model to include both shared and client-specific loading components. 
In the non-private setting,
\citet{yang2025estimating} first estimated each dataset's full signal space, including both shared and dataset-specific components. 
They then aggregated the corresponding projection matrices and took the leading eigenspace as an estimate of the shared space. We could adapt this approach under FDP by aggregating privatized local projectors to estimate the shared space, then estimating client-specific spaces by projecting out the estimated shared space. 
We leave the development of such an extension, including suitable identifiability conditions and estimation error bounds under FDP, to future work.

\bibliographystyle{apalike}
\bibliography{references}

\clearpage
\appendix

\renewcommand{\thesection}{\Alph{section}}
\renewcommand{\theHsection}{\Alph{section}} 
\renewcommand{\thesubsection}{\Alph{section}.\arabic{subsection}}
\renewcommand{\theHsubsection}{\Alph{section}.\arabic{subsection}} 

\setcounter{table}{0}
\renewcommand{\thetable}{A\arabic{table}}
\setcounter{figure}{0}
\renewcommand{\thefigure}{A\arabic{figure}}

\section{Additional details and discussion}\label{sec: additional_detail}

\subsection{Estimation of multiple loading matrices}
\label{subsec: multimode}

Suppose that our interest lies in estimating a subset of modes $\mathcal{K} \subseteq[K]$,
and each client $j$ has a total privacy budget $(\varepsilon_j,\delta_j)$.
Client $j$ allocates budgets $(\varepsilon_{j,k},\delta_{j,k})$ for $k\in\mathcal{K}$ such that
\begin{align*}
    \sum_{k\in\mathcal{K}} \varepsilon_{j,k}\leq \varepsilon_j,
    \quad
    \sum_{k\in\mathcal{K}} \delta_{j,k}\leq \delta_j .
\end{align*}
For each mode $k\in\mathcal{K}$, we run our procedure in the main text with budgets $(\varepsilon_{j,k}, \delta_{j,k})$ and with independently generated noise matrices.

By the composition theorem for DP
\citep[Theorem~3.16]{DworkRoth2014},
the vector of all mode-wise outputs released by client $j$
is $(\sum_{k\in\mathcal{K}} \varepsilon_{j,k}, \sum_{k\in\mathcal{K}} \delta_{j,k})$-CDP,
and hence $(\varepsilon_j,\delta_j)$-CDP.
The server then aggregates the privatized projectors separately for each mode.
This final aggregation step is post-processing of the released client outputs,
and hence it does not incur additional privacy loss.

More generally, the budgets $(\varepsilon_{j,k},\delta_{j,k})$ for every client-mode pair $(j,k)$ can be directly specified. This motivates us to extend the notion of FDP as follows.

\begin{definition}[Multi-mode federated differential privacy]
\label{def: multimode-FDP}
Let $E= (\varepsilon_{j,k})_{j\in[m], k\in[K]}$ and
$D=(\delta_{j,k})_{j\in[m], k\in[K]}$ be the matrices of privacy budgets.
A randomized mechanism
\begin{align*}
\mathcal{M}:
\mathcal{Z}^{n_1}\times \dots \times \mathcal{Z}^{n_m}
\to
\widetilde{\mathcal{Z}}_1 \times \dots \times \widetilde{\mathcal{Z}}_K
\end{align*}
is called $(E,D)$-multi-mode federated differentially private if, for every client $j\in[m]$, the combined output
$(\mathcal{M}_{j,1}, \dots, \mathcal{M}_{j,K})$ over all modes is $\big( \sum_{k=1}^K \varepsilon_{j,k}, \sum_{k=1}^K\delta_{j,k} \big)$-centrally differentially private with respect to client $j$'s local dataset.
\end{definition}

Applying our federated estimation mode-by-mode with independently generated noise yields an $(E,D)$-multi-mode FDP mechanism. In particular, if
\begin{align*}
    \sum_{k=1}^K \varepsilon_{j,k}\leq \varepsilon_j,
    \quad
    \sum_{k=1}^K \delta_{j,k}\leq \delta_j,
\end{align*}
for all $j\in[m]$,
then the resulting mechanism is also
$(\boldsymbol\varepsilon, \boldsymbol\delta)$-federated differentially private in the sense of Definition~\ref{def:FDP}, with
\begin{align*}
    \boldsymbol\varepsilon =
    \Big( \sum_{k=1}^K \varepsilon_{1,k}, \dots, \sum_{k=1}^K \varepsilon_{m,k} \Big)^\top, 
    \quad
    \boldsymbol\delta = \Big( \sum_{k=1}^K \delta_{1,k}, \dots, \sum_{k=1}^K \delta_{m,k} \Big)^\top .
\end{align*}

\subsection{Additional information  on the real data application}
\label{app:real}

The non-private estimators are computed on the same dataset as the private estimators.
To compute the non-private estimator for \textsc{ada} and \textsc{plug},
for each client~$j$, let 
\[
\widehat\Sigma_1^{(j), \mathrm{full}} := n_j^{-1} \sum_{t=1}^{n_j} \mathcal{M}_1(\mathbf{X}_t^{(j)}) \mathcal{M}_1(\mathbf{X}_t^{(j)})^\top
\]
denote the full-sample mode-$1$ covariance matrix, and $\widehat{U}_j^{\mathrm{np}} := \mathrm{SVD}_{r_1}(\widehat\Sigma_1^{(j), \mathrm{full}})$ the matrix of its top $r_1$ eigenvectors.
Then we compute
\[
\widehat{M}_1^{\mathrm{np}}
:=
\sum_{j=1}^m
\frac{n_j}{\sum_{\ell=1}^m n_\ell}
\widehat{U}_j^{\mathrm{np}}(\widehat{U}_j^{\mathrm{np}})^\top,
\]
and the non-private reference is set as $\widehat{A}_1^{\mathrm{np}} := \mathrm{SVD}_{r_1}(\widehat{M}_1^{\mathrm{np}})$.
The weights $\{n_j/\sum_\ell n_\ell\}_j$ are the federated analogue of the pooling, in that if the clients shared the same mode-$1$ loading space, the pooled sample covariance would be the same as the weighted average of the local sample covariances up to the isotropic noise.

The non-private reference for \textsc{hp-vec} is computed analogously as above, except that all estimations are based on the vectorized data.
In detail, its non-private reference consists of the leading $r=\prod_{k=1}^K r_k=8$ eigenvectors of
\[
\widehat\Sigma_{\mathrm{vec}}^{\mathrm{ref}} :=
\frac{1}{m}\sum_{j=1}^m
\left\{
\frac{1}{n_j}\sum_{t=1}^{n_j}
\mathrm{vec}(\mathbf{X}_t^{(j)}) \mathrm{vec}(\mathbf{X}_t^{(j)})^\top
\right\}.
\]

\subsection{Additional numerical results}
\label{app:add:num}

\begin{figure}[htp!]
\centering
\includegraphics[width=\linewidth]{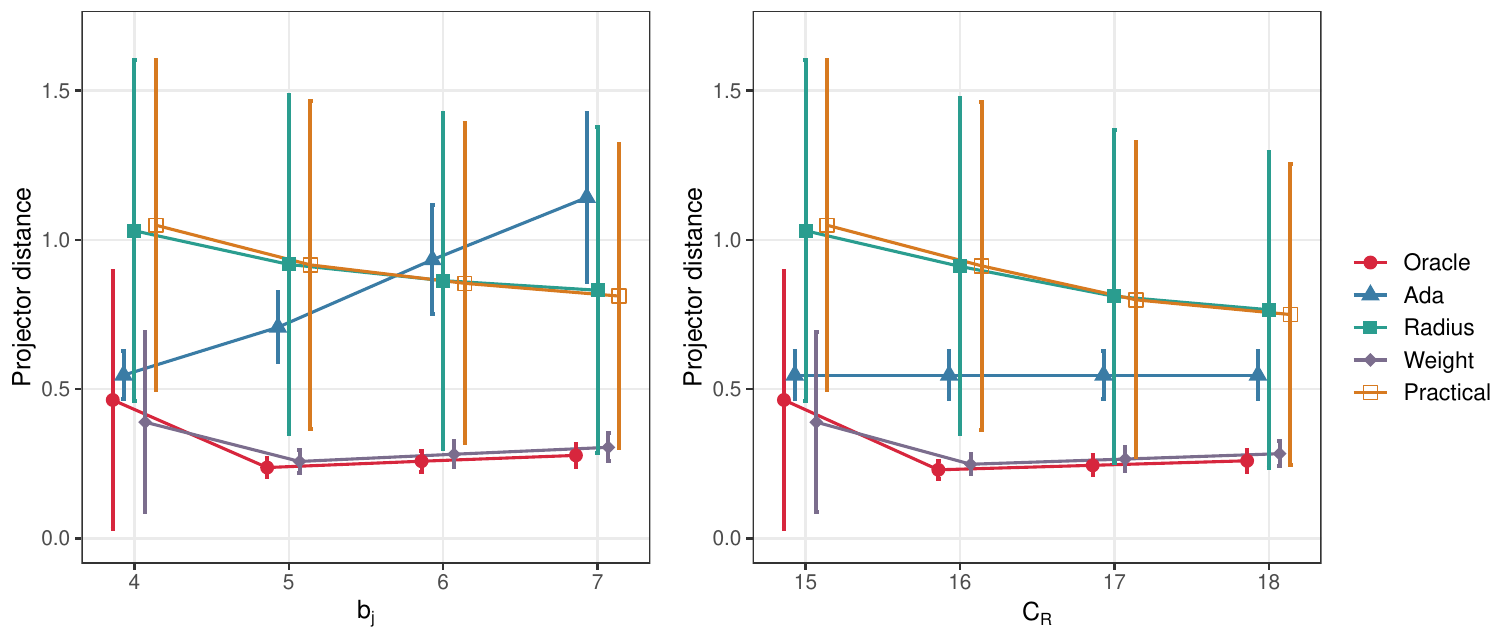}
{\setstretch{1}
\caption{Mean Frobenius projector distance for $K=2$, $p_\ell=15$, $r_\ell=2$, $m=3$, $n_j=4000$, $\boldsymbol\varepsilon=(0.6,0.8,1)^\top$, and $\boldsymbol\delta=(0.2,0.2,0.2)^\top$, with error bars indicating one standard deviation across replications.
The left panel varies $b_j$ with $C_R=15$, while the right panel varies $C_R$ with $b_j=4$.
Markers are slightly offset horizontally for clarity.
The method \textsc{ada} does not depend on $C_R$, so its results are unchanged in the right panel.}
\label{fig:sens_bj_CR}}
\end{figure}

\medskip
\noindent
\textbf{Sensitivity.}
We further examine the sensitivity of the estimators to the choice of $b_j$ and $C_R$. 
Continuing with the data generating process in Section~\ref{subsec:simulation}, we fix $K=2$, $p_1=p_2=15$, $r_1=r_2=2$, and $n_j=4000$,
retaining the privacy budgets and other settings of the preceding experiment. 
We vary $b_j\in \{4,5,6,7\}$ while holding $C_R=15$, and then vary $C_R\in \{15,16,17,18\}$ while holding $b_j=4$.
We experiment \textsc{oracle}, \textsc{ada}, \textsc{radius}, \textsc{weight}, and \textsc{practical} which involve $b_j$;
similarly for $C_R$ except that \textsc{ada} calibrates its radius and does not depend on $C_R$.
We use the same simulated datasets across all settings within each replication,
and the results are shown in Figure~\ref{fig:sens_bj_CR}.

For \textsc{oracle} and \textsc{weight}, increasing the
block size $b_j$ from $4$ to $5$ reduces the mean error and its variability, while further increasing its value produces modest changes. 
We note that the errors of \textsc{ada} increases with the block size,
which is mainly due to the screening rule in~\eqref{eq:adaptive-screen}, as client~1 (with stronger privacy constraint tends to fail the screening when $b_j=6$ or $7$.
Analogous patterns are observed when varying $C_R$, except that the performance of \textsc{ada} is unchanged because its calibrated radius and weights do not depend on $C_R$.
Overall, these results mark a range of relatively stable performance for \textsc{oracle}, while indicating considerable sensitivity of \textsc{ada} to the block size,
and for either, the privacy is always guaranteed as established in our theoretical results.

\section{Proofs for the main results}

\subsection{Proofs on privacy guarantee}

\begin{proof}[Proof of Lemma~\ref{lem:Huber-sensitivity}]
If the two lists are identical, both conclusions are immediate. Otherwise, let $i_0\in[B]$ be the unique index at which they differ, so that $P^{[\ell]\prime} =P^{[\ell]}$ for every $\ell\neq i_0$. Define the two objectives
\begin{equation}
\label{eq:two-Huber-objectives}
F(H):= \sum_{\ell=1}^B \rho_R\big(\|H-P^{[\ell]} \|_{\mathrm{F}} \big) ,
\quad
F'(H):= \sum_{\ell=1}^B \rho_R\big(\|H-P^{[\ell]\prime} \|_{\mathrm{F}} \big) .
\end{equation}
By \eqref{eq:Huber}, $\widehat{H}$ and $\widehat{H}'$ are the minimizers of $F$ and $F'$ over all compatible symmetric matrices, respectively.
Equation~\eqref{eq:Huber-score-sensitivity} is direct. To see this, in the double sum in \eqref{eq:Huber}, only the ordered pairs $(i_0,h)$ and $(\ell,i_0)$, where $h,\ell \neq i_0$, can change; the indicator corresponding to $(i_0,i_0)$ remains equal to one. 
Hence at most $2(B-1)$ indicators can change, and division by $B$ gives \eqref{eq:Huber-score-sensitivity}.

We show \eqref{eq:Huber-aggregate-sensitivity} in the remainder of this proof. Suppose $q\geq 7B/8$. Then note that for each $h\in[B]$, the sum over $\ell$ in \eqref{eq:Huber} is the number of matrices within the $R/8$-neighborhood (in terms of Frobenius norm) of $P^{[h]}$. Since the average of these $B$ neighborhood sizes is $q$, there must exist $h_0\in[B]$ such that
\[
\big| \{\ell\in[B]: \|P^{[\ell]}-P^{[h_0]} \|_{\mathrm{F}} \leq R/8\} \big|
\geq 7B/8 ,
\]
otherwise $q< 7B/8$.
Define the set of block indices
\begin{equation}
\label{eq:Huber-common-cluster}
\mathcal{I} :=
\{\ell\in[B] \setminus \{i_0\}: \|P^{[\ell]} - P^{[h_0]} \|_{\mathrm{F}}\leq R/8 \} .
\end{equation}
Then by the above and recall $B\geq 8$, we have
\begin{equation}
\label{eq:Huber-common-cluster-size}
|\mathcal{I}| \geq 7B/8 -1 \geq 3B/4 .
\end{equation}
Moreover, by construction, $P^{[\ell]\prime}=P^{[\ell]}$ and $\|P^{[\ell]}- P^{[h_0]} \|_{\mathrm{F}}\leq R/8$ for every $\ell\in\mathcal{I}$,
so that the two objectives in \eqref{eq:two-Huber-objectives} contain at least $3B/4$ common terms with $P^{[\ell]}, P^{[\ell]\prime}$ lie in the $R/8$-neighborhood of $P^{[h_0]}$.

We next show that the minimizers of both objectives are close to $P^{[h_0]}$. Consider first $F$ in \eqref{eq:two-Huber-objectives}. 
Let $H$ be any symmetric matrix satisfying $s:= \|H-P^{[h_0]}\|_{\mathrm{F}} \geq R/2$, and set $V:= (H-P^{[h_0]})/s$. Then $V=V^\top$, $\|V\|_{\mathrm{F}}=1$, and $H= P^{[h_0]} + sV$. By the Huber loss definition in \eqref{eq:Huber-loss}, we have 
\begin{equation}
\label{eq:Huber-gradient}
\nabla_H \rho_R\big(\|H-P\|_{\mathrm{F}}\big)
=
(H-P) \min\Big\{1,\frac{R}{\|H-P\|_{\mathrm{F}}}\Big\},
\end{equation}
where the gradient is defined to be zero when $H=P$. In particular, the gradient's Frobenius norm is at most $R$.
For every $\ell\in \mathcal{I}$, \eqref{eq:Huber-common-cluster} gives
\begin{align*}
    \langle H-P^{[\ell]}, V \rangle
    &=
    \langle H-P^{[\ell]}, (H-P^{[h_0]})/s \rangle
    =
    \langle (H- P^{[h_0]}) + (P^{[h_0]} - P^{[\ell]}) , (H-P^{[h_0]})/s \rangle \\
    &=
    s^{-1} \|H-P^{[h_0]} \|_{\mathrm{F}}^2
    + s^{-1} \langle P^{[h_0]} - P^{[\ell]} , H-P^{[h_0]} \rangle \\
    &\geq
    s - s^{-1} \|P^{[h_0]} - P^{[\ell]} \|_{\mathrm{F}} \|H-P^{[h_0]} \|_{\mathrm{F}}
    \geq
    s - R/8 ,
\end{align*}
and similarly, $\|H-P^{[\ell]} \|_{\mathrm{F}} \leq s+ R/8$.
Consequently, with \eqref{eq:Huber-gradient}, we have
\begin{align*}
\big\langle
\nabla_H \rho_R\big( \|H-P^{[\ell]} \|_{\mathrm{F}} \big), V
\big\rangle 
\geq
(s-R/8) \min\Big\{1, \frac{R}{s+R/8} \Big\} 
\geq 3R/8 ,
\end{align*}
where the last inequality is direct to verify by considering $s\leq 7R/8$ and $s> 7R/8$, respectively.
For every $\ell\notin \mathcal{I}$, \eqref{eq:Huber-gradient} and the Cauchy--Schwarz inequality give
\[
\big\langle
\nabla_H \rho_R\big( \|H-P^{[\ell]} \|_{\mathrm{F}} \big), V
\big\rangle 
\geq
- \| \nabla_H \rho_R\big( \|H-P^{[\ell]} \|_{\mathrm{F}} \big) \|_{\mathrm{F}}
\|V\|_{\mathrm{F}}
\geq 
-R ,
\]
which, combined with \eqref{eq:Huber-common-cluster-size}, we obtain
\begin{align}
    \langle \nabla F(H),V \rangle
    &=
    \sum_{\ell\in \mathcal{I}}
    \big\langle
    \nabla_H \rho_R\big(\|H-P^{[\ell]} \|_{\mathrm{F}} \big), V
    \big\rangle 
    +
    \sum_{\ell\notin \mathcal{I}}
    \big\langle
    \nabla_H \rho_R\big(\|H-P^{[\ell]} \|_{\mathrm{F}} \big), V
    \big\rangle  \notag \\
    &\geq
    |\mathcal{I}| \cdot \frac{3R}{8}
    + (B -|\mathcal{I}|) \cdot (-R)
    \geq
    \frac{3B}{4} \frac{3R}{8} - \frac{B}{4} R
    =
    \frac{BR}{32} > 0 .
    \label{eq:Huber-positive-derivative}
\end{align}
For every fixed symmetric matrix $V$ with $\|V\|_{\mathrm{F}} = 1$, we have
\[
\frac{d}{d s} F\big( P^{[h_0]} + s V \big)
=
\left\langle
\nabla F\big( P^{[h_0]} + s V \big), V
\right\rangle.
\]
Hence \eqref{eq:Huber-positive-derivative} implies that $s \mapsto F(P^{[h_0]} + s V)$ is strictly increasing for $s \geq R/2$, and no minimizer of $F$ can lie outside or on the boundary (recall the Huber loss is a $C^1$ function and hence the function $s\mapsto F(P^{[h_0]} + sV)$ is differentiable at $R/2$) of the $R/2$-ball centered at $P^{[h_0]}$. 
The same argument applies to $F'$ as its terms indexed by $\mathcal{I}$ coincide with those of $F$, while every remaining Huber term still has gradient norm at most $R$. Altogether, we can conclude
\begin{equation}
\label{eq:Huber-minimizer-localization}
\|\widehat{H}- P^{[h_0]}\|_{\mathrm{F}} < R/2 ,
\quad
\|\widehat{H}'- P^{[h_0]}\|_{\mathrm{F}} < R/2 .
\end{equation}

It remains to compare the two minimizers. For $t\in[0,1]$, let $H_t:= (1-t) \widehat{H} + t \widehat{H}'$. 
By \eqref{eq:Huber-common-cluster} and \eqref{eq:Huber-minimizer-localization}, for every $\ell\in \mathcal{I}$,
\[
\| H_t -P^{[\ell]} \|_{\mathrm{F}}
\leq
\| H_t -P^{[h_0]} \|_{\mathrm{F}} + \|P^{[h_0]} -P^{[\ell]}\|_{\mathrm{F}}
< R/2+R/8
< R .
\]
Then by the Huber loss function \eqref{eq:Huber-loss}, the terms indexed by $\mathcal{I}$ (i.e., $\rho_R \big(\|\widehat{H}'-P^{[\ell]} \|_{\mathrm{F}} \big)$ and $\rho_R \big(\|\widehat{H}-P^{[\ell]} \|_{\mathrm{F}} \big)$) admit the quadratic loss function form throughout the segment joining $\widehat{H}$ and $\widehat{H}'$. Each such term hence satisfies
\[
\nabla_H \rho_R \big(\|\widehat{H}'-P^{[\ell]} \|_{\mathrm{F}}\big)
-\nabla_H \rho_R \big(\|\widehat{H}-P^{[\ell]}\|_{\mathrm{F}}\big)
=
\widehat{H}' - \widehat{H} .
\]
For every $\ell\notin \mathcal{I}$, convexity of the Huber loss function gives
\[
\Big\langle
\nabla_H \rho_R \big( \|\widehat{H}' -P^{[\ell]} \|_{\mathrm{F}} \big)
-\nabla_H \rho_R \big( \|\widehat{H} -P^{[\ell]} \|_{\mathrm{F}} \big),
\widehat{H}' - \widehat{H}
\Big\rangle
\geq 0 .
\]
Summing over all terms and using \eqref{eq:Huber-common-cluster-size}, we obtain
\begin{align}
    &\;\quad
    \frac{3B}{4} \|\widehat{H}'-\widehat{H}\|_{\mathrm{F}}^2
    \leq
    |\mathcal{I}| \big\langle 
    \widehat{H}' - \widehat{H}, \widehat{H}' - \widehat{H}
    \big\rangle \notag \\
    &\leq
    \sum_{\ell\in \mathcal{I}}
    \big\langle
    \nabla_H \rho_R \big(\|\widehat{H}'-P^{[\ell]} \|_{\mathrm{F}}\big)
    -\nabla_H \rho_R \big(\|\widehat{H}-P^{[\ell]}\|_{\mathrm{F}}\big) ,
    \widehat{H}' - \widehat{H}
    \big\rangle \\
    \leq &
    \big\langle
    \nabla F(\widehat{H}') - \nabla F(\widehat{H}) ,
    \widehat{H}' - \widehat{H}
    \big\rangle.
    \label{eq:Huber-local-strong-convexity}
\end{align}

Since two objectives $F$ and $F'$ are differentiable and are minimized over the linear space of symmetric matrices, their first-order optimality conditions are
\begin{equation}
\label{eq:Huber-first-order-conditions}
\nabla F(\widehat{H})=0,
\quad
\nabla F'(\widehat{H}')=0.
\end{equation}
By \eqref{eq:two-Huber-objectives}, $F$ and $F'$ differ only in the term indexed by $i_0$. Thus, \eqref{eq:Huber-first-order-conditions} imply
\begin{align}
    \|\nabla F(\widehat{H}')\|_{\mathrm{F}}
    &=
    \|\nabla F(\widehat{H}')-\nabla F'(\widehat{H}')\|_{\mathrm{F}} \nonumber\\
    &=
    \Big\|
    \nabla_H \rho_R\big(\|\widehat{H}'-P^{[i_0]} \|_{\mathrm{F}}\big)
    - \nabla_H \rho_R\big(\|\widehat{H}'-P^{[i_0]\prime} \|_{\mathrm{F}} \big)
    \Big\|_{\mathrm{F}}
    \leq 
    2R ,
    \label{eq:Huber-gradient-perturbation}
\end{align}
where the last inequality used the gradient upper bounded by $R$ as mentioned below \eqref{eq:Huber-gradient}.
Combining \eqref{eq:Huber-local-strong-convexity}, 
\eqref{eq:Huber-first-order-conditions}, 
\eqref{eq:Huber-gradient-perturbation},
and the Cauchy--Schwarz inequality, we obtain
\[
\frac{3B}{4} \|\widehat{H}'-\widehat{H}\|_{\mathrm{F}}^2
\leq
\| \nabla F(\widehat{H}') - 0 \|_{\mathrm{F}}
\cdot \|\widehat{H}'-\widehat{H}\|_{\mathrm{F}}
\leq
2R \|\widehat{H}'-\widehat{H}\|_{\mathrm{F}} .
\]
If $\widehat{H}'= \widehat{H}$, \eqref{eq:Huber-aggregate-sensitivity} is immediate. Otherwise, dividing both sides by $\|\widehat{H}'-\widehat{H} \|_{\mathrm{F}}$ yields
\[
\|\widehat{H}'-\widehat{H}\|_{\mathrm{F}}
\leq
\frac{8R}{3B}
\leq
\frac{4R}{B} ,
\]
which proves \eqref{eq:Huber-aggregate-sensitivity}.
This completes the proof of the lemma.
\end{proof}

\begin{proof}[Proof of Proposition~\ref{prop:standard-FDP}]
For notational convenience, we suppress the client and mode indices (i.e., $j$ and $k$ in Algorithm~\ref{alg:fed-est}) below.
Fix a client and condition on a realization of its data-independent random partition. Fix any pair of neighboring datasets and denote their block-projector lists, consistency scores, and Huber aggregates by
\[
\{P^{[\ell]}\}_{\ell=1}^B ,
\quad q ,
\quad \widehat{H} ,
\quad\text{and}\quad
\{P^{[\ell]\prime}\}_{\ell=1}^B ,
\quad q' ,
\quad \widehat{H}' ,
\]
respectively. 
A replacement of one tensor observation changes at most one block projector, so Lemma~\ref{lem:Huber-sensitivity} applies. 

For a dataset with consistency score $q$ (i.e., as $q_{j,k}$ computed in the procedure in Algorithm~\ref{alg:fed-est}), define the outcome
\begin{equation}
\label{eq:private-test-outcome}
T :=
\mathbbm{1}\Bigg\{
q+W
\geq
\frac{7B}{8} + 2 + \frac{4}{\varepsilon} \log\left(\frac{2}{\delta}
\right) \Bigg\},
\quad
W \sim \mathrm{Lap}(4/\varepsilon) .
\end{equation}
The corresponding outcome for the neighboring dataset is defined analogously using $q'$. 
The Gaussian and Haar randomness used after the test are independent of $W$. 
By \eqref{eq:Huber-score-sensitivity}, $q$ has sensitivity at most two. The Laplace mechanism \citep[][Theorem~3.6]{DworkRoth2014} therefore shows that $q+W$, and hence its post-processing $T$, is $(\varepsilon/2, 0)$-CDP.
To explicitly account for the branch selected by the stability test $T$ in \eqref{eq:private-test-outcome}, we prove privacy for the augmented local output $(T,\widetilde{U})$. The actual local estimator $\widetilde{U}$ is a post-processing of this augmented output.

We first consider the calibration of the data under the Gaussian mechanism, i.e., when $T=1$. 
For each pair $(a,b) \in [p_k] \times [p_k]$, let $E_{ab}$ denote the $p_k \times p_k$ matrix with a one in position $(a,b)$ and zeros elsewhere.
With respect to the Frobenius inner product, an orthonormal basis of the symmetric matrices is 
\[
E_{aa} ,
\quad
\frac{E_{ab} + E_{ba}}{\sqrt{2}} ,
\quad 
1 \leq a < b \leq p_k .
\]
It is direct to check that under this basis, the coordinates of $Z= \tau (G+G^\top) /\sqrt{2}$ are mutually independent $\mathcal{N}(0, 2\tau^2)$ random variables.
This follows from that the diagonal entry $Z_{aa} = \sqrt{2} \tau G_{aa}$ has variance $2\tau^2$, while each off-diagonal entry $Z_{ab}$ with $a \neq b$ has variance $\tau^2$; the basis element $(E_{ab}+E_{ba}) /\sqrt{2}$ combines the symmetric pair $(Z_{ab}, Z_{ba})$ with the factor $1/\sqrt{2}$, giving the same variance $2\tau^2$ in the orthonormal basis.
By the definition of $\tau$ in Algorithm~\ref{alg:fed-est},
\begin{equation}
\label{eq:Gaussian-branch-calibration}
2 \tau^2
=
\frac{2(4R/B)^2}{(\varepsilon/2)^2}
\log\left( \frac{2.5}{\delta} \right) .
\end{equation}
Consequently, whenever two Huber aggregates are at Frobenius distance at most $4R/B$, the corresponding sensitivity is $4R/B$, the Gaussian mechanism \citep[][Theorem~A.1]{DworkRoth2014} applied to them is $(\varepsilon/2, \delta/2)$-CDP, and applying singular value decomposition afterward preserves this guarantee.

We now consider two exhaustive cases for the fixed neighboring pair.
First, suppose $q\geq 7B/8$ and $q'\geq 7B/8$. Then by \eqref{eq:Huber-aggregate-sensitivity} in Lemma~\ref{lem:Huber-sensitivity},
\[
\|\widehat{H} -\widehat{H}'\|_{\mathrm{F}}
\leq 4R/B .
\]
Conditional on $T=0$, the output of Algorithm~\ref{alg:fed-est} is drawn from the same data-independent Haar distribution under both datasets and is therefore $(0,0)$-CDP. Conditional on $T=1$, the output of the local estimation in Algorithm~\ref{alg:fed-est} is $\mathrm{SVD}_{r_k}(\widehat{H} +Z)$, which, by \eqref{eq:Gaussian-branch-calibration}, is $(\varepsilon/2, \delta/2)$-CDP for this neighboring pair. 
Hence for either possible value of $T$, Step~2(b) is $(\varepsilon/2, \delta/2)$-CDP. Together with the $(\varepsilon/2, 0)$-CDP score $q$, sequential composition \citep[][Theorem~3.16]{DworkRoth2014} implies that the augmented output $(T, \widetilde{U})$ is $(\varepsilon, \delta/2)$-CDP. 

Second, suppose that at least one of the scores $q$ and $q'$ is smaller than $7B/8$. By \eqref{eq:Huber-score-sensitivity}, both $q$ and $q'$ are smaller than $7B/8+2$. It follows from \eqref{eq:private-test-outcome} that, under either dataset, passing the test with $T=1$ at least requires
\[
W \geq
\frac{4}{\varepsilon} \log\left( \frac{2}{\delta} \right) .
\]
Since $W \sim \mathrm{Lap}(4/\varepsilon)$, we have
\begin{equation}
\label{eq:unstable-test-pass-probability}
\P(T=1)
\leq
\P\left( W
\geq
\frac{4}{\varepsilon} \log\left( \frac{2}{\delta} \right) \right)
=
\frac{1}{2} \exp\Bigg\{
-\frac{(4/\varepsilon) \log(2/\delta)}{4/\varepsilon}
\Bigg\}
=
\frac{\delta}{4}.
\end{equation}
Let $a$ and $a'$ be the probabilities of passing the test with $T=1$ under the two neighboring datasets. Then \eqref{eq:unstable-test-pass-probability} gives $a,a' \leq \delta/4$. 
Let $Q$ be the common distribution of the augmented output $(T,\widetilde{U})$ when the test fails,
so that by construction, under $Q$, it holds that $T=0$ and $\widetilde{U}$ has a Haar distribution which is data-independent.
Let $Q_1$ and $Q_1'$ be the conditional distributions of the augmented output when the test passes under the two datasets, and these two conditional distributions may be arbitrary. 
Then note that the two unconditional output distributions, denoted by $\mu$ and $\mu'$, satisfy
\begin{equation}
\label{eq:unstable-output-mixtures}
\mu = (1-a)Q + a Q_1 ,
\quad
\mu' = (1-a')Q + a'Q_1' .
\end{equation}
Using \eqref{eq:unstable-output-mixtures} and the fact that the total variation distance between two probability measures is at most one, we have
\begin{align*}
d_{\mathrm{TV}}(\mu,\mu')
&=
\sup_A |\mu(A) - \mu'(A)|
=
\sup_A \left| a[Q_1(A) - Q(A)] + a'[Q(A) - Q_1'(A)] \right| \\
&\leq
a \cdot d_{\mathrm{TV}}(Q_1,Q)
+a' \cdot d_{\mathrm{TV}}(Q_1',Q) 
\leq
a + a'
\leq
\delta/2 .
\end{align*}
Hence, for every measurable set $A$ in the augmented output space,
\[
\mu(A)
\leq
\mu'(A) + \delta/2
=
e^0\mu'(A) + \delta/2 ,
\]
which implies that the augmented output $(T,\widetilde{U})$ is $(0,\delta/2)$-CDP for this neighboring pair. 
This case covers all configurations in which at least one of $q$ and $q'$ is below $7B/8$, including the case where both are below and the mixed case where exactly one is below. In all such configurations, the sensitivity bound ensures that both scores are below $7B/8 + 2$, so the pass-probability bound $a, a' \leq \delta/4$ remains valid.

Finally, the above discussion covers every fixed neighboring pair. In the first case, where both $q \geq 7B/8$ and $q' \geq 7B/8$, the augmented local output $(T, \widetilde{U})$ is $(\varepsilon, \delta/2)$-CDP. 
In the second case, where at least one of $q$ and $q'$ is below $7B/8$, it is $(0, \delta/2)$-CDP, which is stronger than $(\varepsilon, \delta/2)$-CDP since $e^{0} = 1 \leq e^{\varepsilon}$.  Hence in either case the augmented local output satisfies $(\varepsilon, \delta/2)$-CDP.
This conclusion holds for every fixed realization of the data-independent partition with the same privacy parameters, and averaging over the random partition preserves DP. Reintroducing the index $j$, the local release of client $j$ is $(\varepsilon_j, \delta_j)$-CDP with respect to $\mathcal{D}_j$. 
For fixed data from all other clients, their messages and every server-side operation are post-processing relative to $\mathcal{D}_j$. Thus the final federated estimator satisfies $(\boldsymbol\varepsilon, \boldsymbol\delta)$-FDP.
\end{proof}

\subsection{Auxiliary lemmas with proofs}

We first establish three auxiliary results: 
(i) the first shows that random partitioning preserves the factor strength in every block, with high probability;
(ii) the second controls the blockwise empirical projectors;
and (iii) the third gives the local stochastic bound used in the federated aggregation.

\begin{lemma}[Factor balance under random partitioning]
\label{lem:block-balance}
Let Assumption~\ref{ass:latent} holds. 
Fix $j \in [m]$ and $k \in [K]$.
Let $\Lambda_{j,(k)}$ be the full-sample mode-$k$ factor covariance defined in \eqref{def-lambda-1}, let $b_j$ and $B_j$ be set as in Theorem~\ref{thm:fed-upper}, and let $\mathcal{I}_{j,1}, \ldots, \mathcal{I}_{j,B_j}$ be the random blocks constructed in Algorithm~\ref{alg:fed-est}. 
Then with probability at least 
\[
1 - 2 r_k B_j \exp\Big\{-\frac{c b_j}{C_\lambda^2} \Big\} ,
\]
it holds that
\begin{align}\label{eq:block-balance} \frac{1}{2} \Lambda_{j,(k)}
\preceq 
\frac{1}{|\mathcal{I}_{j,\ell}|} \sum_{t\in\mathcal{I}_{j,\ell}}
\mathcal{M}_k(\mathbf{F}_t^{(j)}) \mathcal{M}_k(\mathbf{F}_t^{(j)})^\top
\preceq 
\frac{3}{2} \Lambda_{j,(k)}, \quad \forall \ell \in[B_j] ,
\end{align}
where $A \preceq B$ for matrices $A,B$ means $B-A$ is positive semi-definite;
the probability is taken over the random partition; 
$C_\lambda>0$ is the absolute constant in  Assumption~\ref{ass:latent}; and $c>0$ is an absolute constant.
\end{lemma}

\begin{proof}[Proof of Lemma~\ref{lem:block-balance}]
By Assumption~\ref{ass:latent}, we have $\lambda_{j,r_k}^{(k)} >0$ and hence $\Lambda_{j,(k)}$ is positive definite. For each $t\in[n_j]$, define
the normalized mode-$k$ factor matrix as
\begin{equation}
\label{eq:normalized-factor-matrix}
Y_{j,k,t} :=
\Lambda_{j,(k)}^{-1/2}
\mathcal{M}_k(\mathbf{F}_t^{(j)})
\mathcal{M}_k(\mathbf{F}_t^{(j)})^\top
\Lambda_{j,(k)}^{-1/2} .
\end{equation}
By the definition of $\Lambda_{j,(k)}$ in \eqref{def-lambda-1},
\begin{equation}
\label{eq:normalized-factor-mean}
\frac{1}{n_j} \sum_{t=1}^{n_j} Y_{j,k,t}
= I_{r_k} .
\end{equation}
Moreover, the definition of $\lambda_{j,\max}^{(k)}$ in
\eqref{def-lambda-2} gives
\[
0 \preceq
\mathcal{M}_k(\mathbf{F}_t^{(j)}) \mathcal{M}_k(\mathbf{F}_t^{(j)})^\top
\preceq
\lambda_{j,\max}^{(k)} I_{r_k} .
\]
Hence by Assumption~\ref{ass:latent}, we have
\begin{align}
    0 \preceq  Y_{j,k,t} 
    \preceq
    \lambda_{j,\max}^{(k)} \Lambda_{j,(k)}^{-1} 
    \preceq
    \frac{\lambda_{j,\max}^{(k)}}{\lambda_{j,r_k}^{(k)}} I_{r_k}
    \preceq
    C_\lambda^2 I_{r_k} .
    \label{eq:normalized-factor-bound}
\end{align}

Fix $\ell\in[B_j]$. Then conditional on the deterministic factors, the
matrices $Y_{j,k,1}, \dots, Y_{j,k,n_j}$ form a fixed finite population, while the block $\mathcal{I}_{j,\ell}$ is random. 
For every $t\in[n_j]$, define the
inclusion indicator $Z_{j,\ell,t}^{\mathrm{part}} := \mathbbm{1}\{t\in \mathcal{I}_{j,\ell}\}$.
By construction via random permutation, $\mathcal{I}_{j,\ell}$ is a uniformly sampled subset of $[n_j]$ of size $|\mathcal{I}_{j,\ell}|$. Hence we have
\begin{equation}
\label{eq:block-inclusion-probability}
\P_{\pi_j} \{ Z_{j,\ell,t}^{\mathrm{part}} =1 \}
=
\frac{\binom{n_j -1}{|\mathcal{I}_{j,\ell}| -1}}
{\binom{n_j}{|\mathcal{I}_{j,\ell}|}}
=
\frac{|\mathcal{I}_{j,\ell}|}{n_j} ,
\end{equation}
where $\P_{\pi_j}$ denotes probability over the random permutation.
Since the matrices $Y_{j,k,t}$ are fixed conditional on the factors,
\eqref{eq:normalized-factor-mean} and \eqref{eq:block-inclusion-probability} give
\begin{align}
\E_{\pi_j}\Bigg\{
\frac{1}{|\mathcal{I}_{j,\ell}|}
\sum_{t\in\mathcal{I}_{j,\ell}} Y_{j,k,t}
\Bigg\} 
= \frac{1}{|\mathcal{I}_{j,\ell}|}
\sum_{t=1}^{n_j}
\E_{\pi_j}(Z_{j,\ell,t}^{\mathrm{part}})
Y_{j,k,t} 
=
\frac{1}{n_j} \sum_{t=1}^{n_j} Y_{j,k,t} 
= I_{r_k} ,
\label{eq:block-factor-conditional-mean}
\end{align}
where $\E_{\pi_j}$ denotes expectation over the random permutation.

Note that the indicators $Z_{j,\ell,t}^{\mathrm{part}}$ are dependent with each other, since the sampling is without replacement. 
We address this below.
Conditional on the factors, the matrices indexed by $\mathcal{I}_{j,\ell}$ constitute a uniform sample without replacement from a finite population whose mean is $I_{r_k}$ and whose elements satisfy \eqref{eq:normalized-factor-bound}. 
To justify the use of the
matrix Chernoff inequality (see e.g.\ Chapter~5.1 of \citeauthor{Tropp2015}, \citeyear{Tropp2015}), set $s=|\mathcal{I}_{j,\ell}|$ and let
$\xi_1, \dots, \xi_s \in [n_j]$ be independent and uniformly distributed on $[n_j]$. 
Then $Y_{j,k,\xi_1}, \dots, Y_{j,k,\xi_s}$ form the corresponding samples with replacement. 
Recall \eqref{eq:normalized-factor-mean} and \eqref{eq:normalized-factor-bound} which remain to hold, such that
\[
\E(Y_{j,k,\xi_a}) = I_{r_k} ,
\quad
0 \preceq
Y_{j,k,\xi_a}
\preceq
C_\lambda^2 I_{r_k} .
\]
Applying Lemma~5.4.1 in \cite{Tropp2015} gives,
for every $\theta\in\R$,
\begin{equation}
\label{eq:with-replacement-trace-mgf}
\E\Bigg[ \mathrm{tr}\Big\{ \exp\Big(
\theta\sum_{a=1}^s Y_{j,k,\xi_a}
\Big) \Big\} \Bigg]
\leq
r_k \exp\Bigg\{
\frac{s}{C_\lambda^2}
\big(e^{\theta C_\lambda^2} -1 \big)
\Bigg\} .
\end{equation}
By the arguments at the end of \cite{gross2010note}, we have
\begin{equation}
\label{eq:without-replacement-trace-mgf}
\E_{\pi_j}\Bigg[ \mathrm{tr}\Big\{ \exp\Big(
\theta\sum_{t\in\mathcal I_{j,\ell}}Y_{j,k,t}
\Big) \Big\} \Bigg]
\leq
\E\Bigg[ \mathrm{tr}\Big\{ \exp\Big(
\theta\sum_{a=1}^s Y_{j,k,\xi_a}
\Big) \Big\} \Bigg] ,
\quad
\theta\in\R ,
\end{equation}
so that it suffices to consider the bounds under the samples with replacement, which we constructed earlier.
Consequently, although the matrices selected without replacement are dependent, \eqref{eq:with-replacement-trace-mgf} and \eqref{eq:without-replacement-trace-mgf} give their sum the same trace moment-generating-function upper bound as the independent sum used by Lemma~5.4.1 in \cite{Tropp2015}.

Then for every $u\geq 0$, 
Proposition~3.2.1 of \cite{Tropp2015}
gives
\begin{align}
    &\;\quad
    \P_{\pi_j}\bigg\{
    \lambda_{\max}\Big( \sum_{t\in \mathcal{I}_{j,\ell}} Y_{j,k,t} \Big)
    \geq (1+u)s
    \bigg\} \nonumber \\
    &\leq
    \inf_{\theta>0}
    e^{-\theta(1+u)s}
    \E_{\pi_j}\Bigg[ \mathrm{tr}\Big\{ \exp\Big(
    \theta\sum_{t\in\mathcal I_{j,\ell}}Y_{j,k,t}
    \Big) \Big\} \Bigg]
    \leq 
    r_k \Bigg\{
    \frac{e^u}{(1+u)^{1+u}}
    \Bigg\}^{s/ C_\lambda^2} ,
    \label{eq:block-Chernoff-upper}
\end{align}
where the last inequality used \eqref{eq:with-replacement-trace-mgf} and \eqref{eq:without-replacement-trace-mgf} by taking $\theta = \log(1+u)/ C_\lambda^2$.
Similarly, for every $u\in[0,1)$, applying the same argument with
$-\theta$ and taking $\theta = -\log(1-u)/ C_\lambda^2$ gives
\begin{align}
& \P_{\pi_j}\bigg\{
\lambda_{\min}\Big( \sum_{t\in\mathcal{I}_{j,\ell}} Y_{j,k,t} \Big)
\leq (1-u)s
\bigg\}
\leq 
r_k \Bigg\{
\frac{e^{-u}}{(1-u)^{1-u}}
\Bigg\}^{s/ C_\lambda^2} .
\label{eq:block-Chernoff-lower}
\end{align}
Thus Lemma~5.4.1 in \cite{Tropp2015} applies under the above events.
Taking $u=1/2$ in \eqref{eq:block-Chernoff-upper} and \eqref{eq:block-Chernoff-lower} and combining the two tails gives
\begin{align}
\P_{\pi_j}\Bigg\{
\Bigg\| \frac{1}{|\mathcal{I}_{j,\ell}|} \sum_{t\in\mathcal{I}_{j,\ell}} Y_{j,k,t} - I_{r_k} \Bigg\|
> \frac{1}{2}
\Bigg\} 
\leq
2 r_k \exp\Bigg\{
-\frac{c s}{C_\lambda^2}
\Bigg\}
\leq
2 r_k \exp\Bigg\{
-\frac{c b_j}{C_\lambda^2}
\Bigg\} ,
\label{eq:one-block-factor-balance}
\end{align}
where $c >0$ is an absolute constant and the last inequality follows
from setting $s= |\mathcal{I}_{j,\ell}|\geq b_j$.
This completes the proof of this lemma.
\end{proof}

\begin{lemma}[Uniform block-projector concentration]
\label{lem:block-projector}
Suppose that Assumptions~\ref{ass:loading}--\ref{ass:noise} 
and the conditions of Theorem~\ref{thm:fed-upper} hold.
Condition on any realization of the random partitions for which \eqref{eq:block-balance} holds for every client and block.
Let $R_{j,k}$ be defined in \eqref{eq:block-radius}, 
and suppose $R_{j,k} \leq c_{\mathrm{rad}}$ for every $j \in [m]$,
where $c_{\mathrm{rad}} > 0$ is a sufficiently small absolute constant.
Then for a sufficiently large choice of $C_R$, there exist absolute constants $C, c > 0$ such that, with conditional probability at least $1 - C (n \vee p)^{-c}$,
\begin{equation}
\label{eq:all-blocks-good}
\max_{j\in[m]} \max_{\ell\in[B_j]}
\frac{\big\| \widehat\Pi_{j,k}^{[\ell]} - A_kA_k^\top \big\|_{\mathrm{F}} }{R_{j,k}}
\leq
\frac{1}{16} .
\end{equation}
\end{lemma}

\begin{proof}[Proof of Lemma~\ref{lem:block-projector}]
Fix $j\in[m]$ and $\ell\in[B_j]$, and condition on the random partition and the factors. 
For notational simplicity within this proof, write
$\mathcal{I} = \mathcal{I}_{j,\ell}$, $s= |\mathcal{I}|$, and $A_{-k} = \bigotimes_{h\in[K] \setminus\{k\}} A_h$ (where the Kronecker product is always taken in descending order unless otherwise specified).
Also, for each $t\in\mathcal{I}$, define $D_t:= \mathcal{M}_k(\mathbf{F}_t^{(j)})$, $E_t:= \mathcal{M}_k(\mathbf{E}_t^{(j)})$, and $Y_t:= \mathcal{M}_k(\mathbf{X}_t^{(j)})$.
Let block sample covariance used to construct $\widehat\Pi_{j,k}^{[\ell]}$ be
\begin{equation}
\label{eq:block-sample-covariance-proof}
\widehat\Sigma_{j,k}^{[\ell]}
:= \frac{1}{s} \sum_{t\in\mathcal{I}} Y_tY_t^\top .
\end{equation}
Also define $S_D := s^{-1} \sum_{t\in\mathcal I} D_t D_t^\top$.
Then expanding each summand in \eqref{eq:block-sample-covariance-proof},
as the mode-$k$ matricization of $\mathbf{X}_t^{(j)}$ is $Y_t = A_k D_t A_{-k}^\top + E_t$,
and by $A_{-k}^\top A_{-k}  = \bigotimes_{h\in[K] \setminus\{k\}} (A_h^\top A_h) = I_{r_{-k}}$ from Assumption~\ref{ass:loading},
we obtain
\begin{align}
\widehat\Sigma_{j,k}^{[\ell]} &=
A_k S_D A_k^\top 
+ \frac{1}{s} \sum_{t\in\mathcal{I}}
\big( A_k D_t A_{-k}^\top E_t^\top + E_t A_{-k} D_t^\top A_k^\top \big)
+ \frac{1}{s} \sum_{t\in\mathcal{I}} E_t E_t^\top .
\label{eq:block-covariance-expansion}
\end{align}
As the entries of $E_t$ are independent
$\mathcal{N}(0,\sigma^2)$ random variables, we have $\mathbb{E}(E_t E_t^\top) = \sigma^2p_{-k} I_{p_k}$ and
$\mathbb{E}\big( A_k D_t A_{-k}^\top E_t^\top \mid D_t \big) = 0$.
Hence conditional on the factors and the partition, we have
\begin{equation}
\label{eq:block-conditional-mean}
\Sigma^\circ :=
\mathbb{E} \left[ \widehat\Sigma_{j,k}^{[\ell]}
\mid
\{D_t\}_{t\in\mathcal{I}}, \mathcal{I} \right]
=
A_k S_D A_k^\top
+ \sigma^2 p_{-k} I_{p_k} .
\end{equation}

Define the centered perturbation $\Delta := \widehat\Sigma_{j,k}^{[\ell]} -\Sigma^\circ$,
which, by \eqref{eq:block-covariance-expansion} and
\eqref{eq:block-conditional-mean}, reduces to
\begin{align*}
    \Delta &=
    \frac{1}{s} \sum_{t\in\mathcal{I}}
    \big( A_k D_t A_{-k}^\top E_t^\top
    + E_t A_{-k} D_t^\top A_k^\top \big) +
    \frac{1}{s} \sum_{t\in\mathcal{I}}
    \big( E_t E_t^\top - \sigma^2 p_{-k} I_{p_k} \big) 
    =:
    \mathcal{C} + \mathcal{C}^\top
    + \mathcal{W} ,
\end{align*}
where
\[
\mathcal{C} :=
\frac{1}{s} \sum_{t\in\mathcal{I}} A_k D_t A_{-k}^\top E_t^\top ,
\quad
\mathcal{W} := 
\frac{1}{s} \sum_{t\in\mathcal{I}}
\big( E_t E_t^\top - \sigma^2 p_{-k} I_{p_k} \big) .
\]
Note that the triangle inequality gives
\begin{equation}
\label{eq:block-perturbation-components}
\|\Delta\| \leq
2 \|\mathcal{C}\| + \|\mathcal{W}\| .
\end{equation}

Consider $\mathcal{C}$ first.
Define $\widetilde{E}_t := E_t A_{-k}$,
so that $\mathcal{C} = s^{-1} A_k \sum_{t\in\mathcal{I}} D_t \widetilde{E}_t^\top$.
Since $A_{-k}$ has orthonormal columns, the entries of $\widetilde{E}_t \in \mathbb{R}^{p_k\times r_{-k}}$ are independent $\mathcal{N}(0, \sigma^2)$ random variables. 
Hence conditional on the factors and the partition, the columns of $s^{-1} \sum_{t\in\mathcal{I}} D_t \widetilde{E}_t^\top$ are independent Gaussian vectors with mean zero and covariance matrix
\[
\frac{\sigma^2}{s^2}
\sum_{t\in\mathcal{I}} D_t D_t^\top
= \frac{\sigma^2}{s} S_D .
\]
Altogether, we have
\begin{equation}
\label{eq:block-cross-gaussian-representation}
\mathcal{C}
\overset{d}{=}
\frac{\sigma}{\sqrt{s}} A_k S_D^{1/2} G_1 ,
\end{equation}
where $G_1\in \mathbb{R}^{r_k\times p_k}$ has independent standard Gaussian entries. For every $u>0$, the standard Gaussian spectral-norm bound \citep[e.g.][Theorem~4.4.5]{vershynin2018high} gives
\[
\|G_1\| \leq
\sqrt{p_k} + \sqrt{r_k} + \sqrt{2u} ,
\]
with probability at least $1- \exp(-u)$. In addition, \eqref{eq:block-balance} in Lemma~\ref{lem:block-balance} implies $\|S_D\| \leq (3/2) \lambda_{j,1}^{(k)}$.
Therefore, it holds with probability at least $1- \exp(-u)$ that 
\begin{equation}
\label{eq:block-cross-concentration}
\|\mathcal{C}\| \leq
C \sigma \sqrt{\lambda_{j,1}^{(k)}}
\sqrt{\frac{p_k + r_k + u}{s}} .
\end{equation}

Next, consider $\mathcal{W}$. Stack the noise matrices column-wise as $\mathcal{E} := [E_t: t\in \mathcal{I}] \in \mathbb{R}^{p_k\times s p_{-k}}$.
Note that we can write $\mathcal{E}= \sigma G_2$, where $G_2 \in \mathbb{R}^{p_k\times s p_{-k}}$ has independent standard Gaussian entries. Hence,
\[
\mathcal{W} =
\frac{\sigma^2}{s}
\big( G_2 G_2^\top - s p_{-k} I_{p_k} \big) .
\]
The Gaussian sample-covariance bound in
\citet{vershynin2018high} (see Equation~(4.22) below Theorem~4.6.1 within) gives, with probability at least
$1- 2\exp(-u)$,
\begin{align}
    \|\mathcal{W}\|
    &\leq
    \frac{C \sigma^2}{s} (sp_{-k})
    \left\{
    \sqrt{\frac{p_k}{sp_{-k}}}
    + \frac{\sqrt{u}}{\sqrt{sp_{-k}}}
    + \frac{p_k}{sp_{-k}}
    + \frac{u}{sp_{-k}}
    \right\} \notag \\
    &=
    C \sigma^2 \left\{
    \sqrt{\frac{p_{-k} (p_k+u)}{s}}
    + \frac{p_k + u}{s}
    \right\} .
    \label{eq:block-Wishart-concentration}
\end{align}
Combining \eqref{eq:block-perturbation-components},
\eqref{eq:block-cross-concentration}, 
and \eqref{eq:block-Wishart-concentration}, we obtain that, with probability at least $1- 3\exp(-u)$,
\begin{align}
    \|\Delta\| &\leq
    C\Bigg[
    \sigma \sqrt{\lambda_{j,1}^{(k)}}
    \sqrt{\frac{p_k + r_k + u}{s}}
    + \sigma^2 \left\{
    \sqrt{\frac{p_{-k} (p_k + u)}{s}}
    + \frac{p_k + u}{s}
    \right\}
    \Bigg] .
    \label{eq:block-total-perturbation-bound}
\end{align}

It remains to convert the above bound into a projector perturbation bound. Recall Lemma~\ref{lem:block-balance}, it holds with high probability that
\[
\frac{1}{2} \Lambda_{j,(k)}
\preceq
S_D
\preceq
\frac{3}{2} \Lambda_{j,(k)} .
\]
It follows from \eqref{eq:block-conditional-mean} that the leading $r_k$-dimensional eigenspace of $\Sigma^\circ$ is $\mathrm{span}(A_k)$, with spectral projector $P_k = A_k A_k^\top$.
Consequently,
\begin{align}
    \lambda_{r_k}(\Sigma^\circ)
    - \lambda_{r_k +1}(\Sigma^\circ)
    &=
    \lambda_{r_k}(S_D) 
    \geq
    \frac{1}{2} \lambda_{j,r_k}^{(k)} .
    \label{eq:block-population-eigengap}
\end{align}

As $\widehat\Pi_{j,k}^{[\ell]}$ is the spectral projector onto the
leading $r_k$-dimensional eigenspace of $\widehat\Sigma_{j,k}^{[\ell]} = \Sigma^\circ + \Delta$,
the Davis--Kahan theorem \citep[Theorem~2]{yu2015useful} and \eqref{eq:block-population-eigengap} give
\begin{align}
    \big\| \widehat\Pi_{j,k}^{[\ell]} - A_k A_k^\top \big\|_{\mathrm{F}}
    &\leq
    C \sqrt{r_k}
    \frac{\|\Delta\|}{\lambda_{r_k}(S_D)} 
    \leq
    \frac{C \sqrt{r_k}}{\lambda_{j,r_k}^{(k)}}
    \big( 2 \|\mathcal{C}\| + \|\mathcal{W}\| \big) .
    \label{eq:block-Davis-Kahan-components}
\end{align}
Then substituting \eqref{eq:block-cross-concentration} and \eqref{eq:block-Wishart-concentration} into \eqref{eq:block-Davis-Kahan-components} yields
\begin{align}
    \big\| \widehat\Pi_{j,k}^{[\ell]} - A_k A_k^\top \big\|_{\mathrm{F}}
    \leq
    C \sqrt{r_k} \Bigg[
    \frac{\sigma \sqrt{\lambda_{j,1}^{(k)}}}{\lambda_{j,r_k}^{(k)}}
    \sqrt{\frac{p_k + r_k + u}{s}}
    + \frac{\sigma^2}{\lambda_{j,r_k}^{(k)}}
    \left\{
    \sqrt{\frac{p_{-k} (p_k+u)}{s}}
    + \frac{p_k + u}{s}
    \right\}
    \Bigg] .
    \label{eq:block-projector-fixed-bound}
\end{align}

Fix an absolute constant $c > 0$ and, for each client $j \in [m]$,
take $u_j = \log\{ 16 m n_j (n \vee p)^c \}$.
Since $m \leq n$, $n_j \leq n$, and $n \vee p \geq 2$, it holds that $c \log(n \vee p) \leq u_j \leq (c + 6) \log(n \vee p)$,
and hence $u_j \asymp \log(n \vee p)$, with comparison constants depending only on the fixed constant $c$.
Also, as $s\geq b_j$ (just by construction in Algorithm~\ref{alg:fed-est}),
with a sufficiently large choice of $C_R$, the definition of $R_{j,k}$ in \eqref{eq:block-radius} implies
\[
\big\| \widehat\Pi_{j,k}^{[\ell]} - A_k A_k^\top \big\|_{\mathrm{F}}
\leq
\frac{R_{j,k}}{16} 
\]
for this fixed block, except on an event of probability at most $3\exp(-u_j)$. 
Recall the condition $R_{j,k} \leq c_{\mathrm{rad}}$, for sufficiently small $c_{\mathrm{rad}}$,
Finally, as $B_j\leq n_j$, a union bound over all clients and blocks gives
\begin{align*}
    & \mathbb{P}\left\{
    \max_{j\in[m]} \max_{\ell\in[B_j]}
    \frac{\big\| \widehat\Pi_{j,k}^{[\ell]} - A_k A_k^\top \big\|_{\mathrm{F}} }{R_{j,k} }
    > \frac{1}{16} \right\}
    \leq
    \sum_{j=1}^m 3 B_j \exp(-u_j)
    \leq
    \sum_{j=1}^m \frac{3 n_j}{16 m n_j(n\vee p)^c}
    \leq
    C(n\vee p)^{-c} ,
\end{align*}
as desired.
This shows \eqref{eq:all-blocks-good} and thus concludes the proof of this lemma.
\end{proof}

\begin{lemma}[Local ideal release]
\label{lem:local-ideal}
Let all conditions in Theorem~\ref{thm:fed-upper} hold. 
Then for every client~$j$, there exists a random projector $\breve\Pi_{j,k}$ with the following properties:
\begin{enumerate}[label=(\roman*)]
\item the projectors $\breve\Pi_{1,k}, \ldots, \breve\Pi_{m,k}$ are independent;
\item It holds with probability at least $1-C(n\vee p)^{-c}$ that $\breve\Pi_{j,k} = \widetilde{U}_{j,k} \widetilde{U}_{j,k}^\top$ for every $j\in[m]$;
\item $(I_{p_k} - A_k A_k^\top) \E(\breve\Pi_{j,k}) A_k A_k^\top = 0$;
\item $\|\E(\breve\Pi_{j,k}) - A_k A_k^\top\| \leq c$ for a sufficiently small absolute constant $c$;
\item it holds almost surely that $\big\| \breve\Pi_{j,k} - \E(\breve\Pi_{j,k}) \big\|_{\mathrm{F}} \leq C \sqrt{r_k} \Psi_{j,k}^{\mathrm{loc}}$.
\end{enumerate}
Here the expectations in parts~(iii)--(v) are taken over the observation noises, the data-independent random partitions, and the auxiliary Gaussian randomness used in the construction.
\end{lemma}

\begin{proof}[Proof of Lemma~\ref{lem:local-ideal}]
Fix any $j\in[m]$ and write $P:= A_k A_k^\top$ for the target projector. Recall from our procedure in Algorithm~\ref{alg:fed-est} that $\mathcal{I}_{j,1}, \ldots, \mathcal{I}_{j,B_j}$ are the random blocks, 
$\widehat\Pi_{j,k}^{[\ell]}$ is the empirical projector computed from block $\mathcal{I}_{j,\ell}$, 
and $Z_{j,k}$ is the Gaussian orthogonal ensemble perturbation with scale $\tau_{j,k}$. 
Without loss of generality, we may generate $Z_{j,k}$ before the private test and ignore it when the test fails,
and note that this coupling does not change the distribution of the released estimator. 
Let $\mathcal{P}_j :=(\mathcal{I}_{j,1}, \ldots, \mathcal{I}_{j,B_j})$ denote the random partition and let $\mathcal{B}_j$ denote the event that \eqref{eq:block-balance} holds for all $\ell \in[B_j]$. Define the truncated block projector
\begin{equation}
\label{eq:truncated-block-projector}
\overline\Pi_{j,k}^{[\ell]}
:=
\begin{cases}
\widehat\Pi_{j,k}^{[\ell]}, & \mathcal{B}_j \text{ holds and } 
\|\widehat\Pi_{j,k}^{[\ell]} - P\|_{\mathrm{F}} \leq R_{j,k}/16 , \\
P,  & \text{otherwise} .
\end{cases}
\end{equation}
It is noteworthy that this construction will be used only in the proof as it depends on the unknown $P$.
By \eqref{eq:truncated-block-projector}, it holds almost surely that
\begin{equation}
\label{eq:truncated-block-bound}
\|\overline\Pi_{j,k}^{[\ell]} - P\|_{\mathrm{F}}
\leq R_{j,k} /16 .
\end{equation}
Moreover, if $\mathcal{B}_j$ and \eqref{eq:all-blocks-good} hold, then $\overline\Pi_{j,k}^{[\ell]} = \widehat\Pi_{j,k}^{[\ell]}$ for every $\ell\in[B_j]$.
Conditional on $\mathcal{P}_j$, the blocks are disjoint. Assumption~\ref{ass:noise} therefore implies that $\overline\Pi_{j,k}^{[1]}, \ldots, \overline\Pi_{j,k}^{[B_j]}$ are conditionally independent. 
Define their conditional average and centered summands by
\begin{equation}
\label{eq:conditional-block-mean}
\overline M_{j,k} :=
\E\Bigg\{
\frac{1}{B_j} \sum_{\ell=1}^{B_j} \overline\Pi_{j,k}^{[\ell]}
\mathrel{\Big|} \mathcal{P}_j
\Bigg\} ,
\quad
X_{j,k}^{[\ell]} :=
\overline\Pi_{j,k}^{[\ell]} - \E\big\{\overline\Pi_{j,k}^{[\ell]} 
\mid \mathcal{P}_j \big\} ,
\end{equation}
respectively.
The matrices $X_{j,k}^{[\ell]}$ are conditionally independent and centered and also, by \eqref{eq:truncated-block-bound} and Jensen's inequality:
\begin{equation}
\label{eq:centered-block-bound}
\|X_{j,k}^{[\ell]}\|_{\mathrm{F}}
\leq
\|\overline\Pi_{j,k}^{[\ell]} - P\|_{\mathrm{F}}
+ \E\big\{ \|\overline\Pi_{j,k}^{[\ell]} - P\|_{\mathrm{F}} \mid \mathcal{P}_j \big\}
\leq R_{j,k}/8 .
\end{equation}
The symmetric $p_k\times p_k$ matrices equipped with the Frobenius inner product form a Hilbert space. 
Hence applying Theorem~3.5 of \cite{pinelis1994optimum}, 
conditionally on $\mathcal{P}_j$ to \eqref{eq:conditional-block-mean} and \eqref{eq:centered-block-bound}, gives for every $u>0$ that
\begin{equation}
\label{eq:Hilbert-Hoeffding-local}
\P\Bigg\{
\Bigg\| \frac{1}{B_j} \sum_{\ell=1}^{B_j} X_{j,k}^{[\ell]}
\Bigg\|_{\mathrm{F}}
> C R_{j,k} \sqrt{\frac{u}{B_j}}
\mathrel{\Big|} \mathcal{P}_j
\Bigg\}
\leq
2 e^{-u}.
\end{equation}
Since $m\leq n\leq n\vee p$, taking $u=(c+1) \log(n\vee p)$ in \eqref{eq:Hilbert-Hoeffding-local} and then averaging over $\mathcal{P}_j$, we have
\begin{align}
    \Bigg\|
    \frac{1}{B_j} \sum_{\ell=1}^{B_j}  \overline\Pi_{j,k}^{[\ell]} - \overline M_{j,k}
    \Bigg\|_{\mathrm{F}}
    &\leq
    C R_{j,k} \sqrt{\frac{\log(n\vee p)}{B_j}} ,
    \label{eq:centered-block-average}
\end{align}
with probability at least $1 - C(n\vee p)^{-c}/m$. Also, the standard Gaussian spectral-norm bound \citep[][Equation~(4.22)]{vershynin2018high} gives
\begin{equation}
\label{eq:GOE-tail}
\|Z_{j,k}\|
\leq
C \tau_{j,k} \sqrt{p_k + \log(n\vee p)} ,
\end{equation}
with probability at least $1 - C(n\vee p)^{-c}/m$. Let $\mathcal{E}_{j,k}$ denote the intersection of the events in \eqref{eq:centered-block-average} and \eqref{eq:GOE-tail}, so that
\begin{equation}
\label{eq:ideal-good-event-probability}
\P(\mathcal{E}_{j,k}^c)
\leq C (n\vee p)^{-c} /m .
\end{equation}
We next identify the leading eigenspace of the conditional mean $\overline{M}_{j,k}$. Let $O:= 2P - I_{p_k}$. As $P$ is an orthogonal projector, $O^\top O=I_{p_k}$, $OP=PO=P$, $(I_{p_k} -P) O = -(I_{p_k} -P)$, and $OA_k=A_k$.
By \eqref{eq:fed-model}, left multiplication of the mode-$k$ matricization by $O$ leaves its signal component unchanged, and hence with Assumption~\ref{ass:noise}, we have
\[
O \mathcal{M}_k(\mathbf{E}_t^{(j)}) \stackrel{d}{=} 
\mathcal{M}_k(\mathbf{E}_t^{(j)}) ,
\]
so that conditional on $\mathcal{P}_j$, $O \widehat\Pi_{j,k}^{[\ell]} O \stackrel{d}{=} \widehat\Pi_{j,k}^{[\ell]}$. 
By $OPO=P$ and the fact that the Frobenius norm is invariant to orthogonal transformation, the truncation in \eqref{eq:truncated-block-projector} is equivariant under conjugation by $O$. We then arrive at
\begin{equation}
\label{eq:truncated-reflection-invariance}
O \overline\Pi_{j,k}^{[\ell]} O
\stackrel{d}{=}
\overline\Pi_{j,k}^{[\ell]} ,
\quad
O\overline{M}_{j,k} O = \overline{M}_{j,k}  .
\end{equation}
By $OP=P$ and $(I_{p_k} -P) O = -(I_{p_k} -P)$ in \eqref{eq:truncated-reflection-invariance}, we have
\begin{equation}
\label{eq:conditional-mean-cross-block}
(I_{p_k} -P) \overline{M}_{j,k} P
=
-(I_{p_k} -P) \overline{M}_{j,k} P
= 0 .
\end{equation}
Hence $\overline{M}_{j,k}$ is block diagonal with respect to $\mathrm{span}(A_k) \oplus \mathrm{span}(A_k)^\perp$.
Furthermore, \eqref{eq:truncated-block-bound} and Jensen's inequality imply
\begin{equation}
\label{eq:conditional-mean-close}
\|\overline{M}_{j,k} - P\|
\leq
\|\overline{M}_{j,k} -P\|_{\mathrm{F}}
\leq R_{j,k}/16 .
\end{equation}
Combining \eqref{eq:conditional-mean-cross-block} and \eqref{eq:conditional-mean-close}, the leading $r_k$-dimensional eigenspace of $\overline{M}_{j,k}$ is $\mathrm{span}(A_k)$ and its eigengap is at least $1-R_{j,k}/8$.
Note that the lower bound in \eqref{eq:upper-conditions-1} and our choice of block size choice ensure $R_{j,k}$ is sufficiently small, so this eigengap is bounded below by an absolute positive constant.
Next, define the ideal loading matrix and its projector by
\begin{equation}
\label{eq:ideal-private-projector}
U_{j,k}^{\circ} :=
\mathrm{SVD}_{r_k}\Bigg\{
\frac{1}{B_j} \sum_{\ell=1}^{B_j} \overline\Pi_{j,k}^{[\ell]} + Z_{j,k}
\Bigg\} ,
\quad
\Pi_{j,k}^{\circ}:= 
U_{j,k}^{\circ} (U_{j,k}^{\circ})^\top .
\end{equation}
On $\mathcal{E}_{j,k}$, write the matrix inside the SVD above as $\overline{M}_{j,k} + \Delta_{j,k}$, where
\[
\Delta_{j,k} :=
\frac{1}{B_j} \sum_{\ell=1}^{B_j} \overline\Pi_{j,k}^{[\ell]} - \overline{M}_{j,k} + Z_{j,k} .
\]
Combining \eqref{eq:centered-block-average}, \eqref{eq:GOE-tail}, and \eqref{eq:Psi-local}, we have $\|\Delta_{j,k}\| \leq C \sqrt{r_k} \Psi_{j,k}^{\mathrm{loc}}$. 
Then by Davis--Kahan theorem \citep[Theorem~2]{yu2015useful}, we have
\begin{align}
    \|\Pi_{j,k}^{\circ} -P\|_{\mathrm{F}}
    &\leq
    C \| \Delta_{j,k} P\|_{\mathrm{F}}
    \leq
    C \Bigg\{ \Bigg\|
    \frac{1}{B_j} \sum_{\ell=1}^{B_j} \overline\Pi_{j,k}^{[\ell]} - \overline{M}_{j,k}
    \Bigg\|_{\mathrm{F}}
    + \sqrt{r_k} \|Z_{j,k}\|
    \Bigg\} \nonumber\\
    &\leq
    C \Bigg\{
    R_{j,k} \sqrt{\frac{\log(n\vee p)}{B_j}}
    + \sqrt{r_k} \tau_{j,k} \sqrt{p_k + \log(n\vee p)}
    \Bigg\}
    =
    C \sqrt{r_k} \Psi_{j,k}^{\mathrm{loc}} .
    \label{eq:ideal-projector-bound}
\end{align}
Now define projector
\begin{equation}
\label{eq:proof-only-projector}
\breve\Pi_{j,k} :=
\begin{cases}
\Pi_{j,k}^{\circ} , & \mathcal{E}_{j,k} 
\text{ holds} ,\\
P , & \mathcal{E}_{j,k}^c 
\text{ holds}.
\end{cases}
\end{equation}
We next verify parts~(iii)--(v) in Lemma~\ref{lem:local-ideal}. 
The GOE distribution satisfies $O Z_{j,k} O \stackrel{d}{=} Z_{j,k}$. Together with \eqref{eq:truncated-reflection-invariance}, this shows that $O \Pi_{j,k}^{\circ} O \stackrel{d}{=} \Pi_{j,k}^{\circ}$ conditional on $\mathcal{P}_j$. 
The event $\mathcal{E}_{j,k}$ is also invariant under this conjugation because it is defined by the Frobenius and spectral norms in \eqref{eq:centered-block-average} and \eqref{eq:GOE-tail}. Since $O P O=P$, \eqref{eq:proof-only-projector} gives $O \breve\Pi_{j,k} O \stackrel{d}{=} \breve\Pi_{j,k}$. Taking expectations and using $OP=P$ and $(I_{p_k} -P)O = -(I_{p_k} -P)$, we obtain
\begin{equation}
\label{eq:ideal-expectation-cross-block}
(I_{p_k} -P) \E(\breve\Pi_{j,k}) P =0 ,
\end{equation}
which is part~(iii).
By \eqref{eq:ideal-projector-bound} and \eqref{eq:proof-only-projector}, it holds almost surely that
\begin{equation}
\label{eq:proof-only-radius}
\|\breve\Pi_{j,k} -P\|_{\mathrm{F}}
\leq 
C \sqrt{r_k} \Psi_{j,k}^{\mathrm{loc}} .
\end{equation}
By Jensen's inequality, \eqref{eq:proof-only-radius}, and the last condition in \eqref{eq:upper-conditions}, we have
\[
\|\E(\breve\Pi_{j,k}) -P\|
\leq
\E\|\breve\Pi_{j,k} -P\|_{\mathrm{F}}
\leq 
C \sqrt{r_k} \Psi_{j,k}^{\mathrm{loc}} 
\leq c ,
\]
after choosing the constant $C_0$ in \eqref{eq:upper-conditions-1} and \eqref{eq:upper-conditions-2} sufficiently large.
This proves part~(iv). 
Also by the triangle and Jensen inequalities, we have
\begin{align*}
    \|\breve\Pi_{j,k} -\E(\breve\Pi_{j,k}) \|_{\mathrm{F}}
    &\leq
    \|\breve\Pi_{j,k} -P\|_{\mathrm{F}}
    + \|\E(\breve\Pi_{j,k}) -P\|_{\mathrm{F}}
    \leq 
    C \sqrt{r_k} \Psi_{j,k}^{\mathrm{loc}} ,
\end{align*}
which proves part~(v).

We next prove part~(i) in Lemma~\ref{lem:local-ideal}. 
For different clients, the data noises are independent by Assumption~\ref{ass:noise}, and the random partitions and GOE matrices are generated independently by construction. Note that all constructions in \eqref{eq:truncated-block-projector}--\eqref{eq:proof-only-projector} use only the corresponding client's random variables and the fixed matrix $P$. Hence $\breve\Pi_{1,k}, \ldots, \breve\Pi_{m,k}$ are independent.

It remains to establish the simultaneous coupling in part~(ii) of Lemma~\ref{lem:local-ideal}. 
On the intersection of $\mathcal{B}_j$ and the event in \eqref{eq:all-blocks-good}, the definition \eqref{eq:truncated-block-projector} gives $\overline\Pi_{j,k}^{[\ell]} = \widehat\Pi_{j,k}^{[\ell]}$
for every $\ell \in [B_j]$.
Moreover, for any $\ell,h \in[B_j]$,
\[
\|\widehat\Pi_{j,k}^{[\ell]} -\widehat\Pi_{j,k}^{[h]}\|_{\mathrm{F}}
\leq
\|\widehat\Pi_{j,k}^{[\ell]} -P\|_{\mathrm{F}}
+ \|\widehat\Pi_{j,k}^{[h]} -P\|_{\mathrm{F}}
\leq
R_{j,k}/8 .
\]
Hence, the score defined in Algorithm~\ref{alg:fed-est} satisfies $q_{j,k}=B_j$. 
Let $H_{j,k}^{0} := B_j^{-1} \sum_{\ell=1}^{B_j} \widehat\Pi_{j,k}^{[\ell]}$. 
For every $\ell\in[B_j]$, we have
\[
\|H_{j,k}^{0} -\widehat\Pi_{j,k}^{[\ell]} \|_{\mathrm{F}}
\leq
\frac{1}{B_j} \sum_{h=1}^{B_j} \|\widehat\Pi_{j,k}^{[h]} -\widehat\Pi_{j,k}^{[\ell]} \|_{\mathrm{F}}
\leq 
R_{j,k} /8
<R_{j,k} .
\]
Thus, every Huber term is quadratic at $H_{j,k}^{0}$, and the gradient of the Huber objective at this point is $\sum_{\ell=1}^{B_j} (H_{j,k}^{0} -\widehat\Pi_{j,k}^{[\ell]}) =0$.
The convexity of the objective in the minimization problem of Algorithm~\ref{alg:fed-est} and its strict convexity in a neighborhood of $H_{j,k}^{0}$ imply that its selected minimizer is
\begin{equation}
\label{eq:Huber-equals-block-average}
\widehat{H}_{j,k} = H_{j,k}^{0}
= \frac{1}{B_j} \sum_{\ell=1}^{B_j} \widehat\Pi_{j,k}^{[\ell]} .
\end{equation}
When $q_{j,k}=B_j$, the test in Algorithm~\ref{alg:fed-est} fails only if
\[
W_{j,k} <
2 + \frac{4}{\varepsilon_j} \log\left(\frac{2}{\delta_j}\right) - \frac{B_j}{8} .
\]
Note that the first condition in \eqref{eq:upper-conditions} makes the right-hand side at most $-4 \log\{4m (n\vee p)^c\} /\varepsilon_j$. Since $W_{j,k} \sim \mathrm{Lap}(4 /\varepsilon_j)$,
\begin{equation}
\label{eq:ideal-test-failure}
\P\{\text{the test of client~$j$ fails} \mid q_{j,k}=B_j\}
\leq
\frac{1}{8 m (n\vee p)^c}
\leq
\frac{(n\vee p)^{-c}}{4 m} .
\end{equation}
Under $\mathcal{B}_j$, \eqref{eq:all-blocks-good}, $\mathcal{E}_{j,k}$, and the event that the test passes, then 
\eqref{eq:ideal-private-projector}, 
\eqref{eq:proof-only-projector}, 
and \eqref{eq:Huber-equals-block-average} altogether give $\breve\Pi_{j,k} = \Pi_{j,k}^{\circ} = \widetilde{U}_{j,k} \widetilde{U}_{j,k}^\top$.
By Lemma~\ref{lem:block-balance} and our choice of $b_j$ as in Theorem~\ref{thm:fed-upper}, a sufficiently large $C_b$ makes the probability that $\mathcal{B}_j$ fails for at least one client at most $C (n\vee p)^{-c}$.
Finally, note that Lemma~\ref{lem:block-projector} controls \eqref{eq:all-blocks-good}, 
and also, \eqref{eq:ideal-good-event-probability} and a union bound control $\cap_{j=1}^m \mathcal{E}_{j,k}$, 
and note that \eqref{eq:ideal-test-failure} controls the private tests.
Altogether, we conclude
\[
\P\big\{\breve\Pi_{j,k} = \widetilde{U}_{j,k} \widetilde{U}_{j,k}^\top, \ \forall j\in[m] \big\}
\geq
1 - C (n\vee p)^{-c} ,
\]
which proves part~(ii) of this lemma and thus completes this proof.
\end{proof}

\begin{lemma}[Record-level private information contraction for Gaussian models]
\label{lem:score-contraction}
Let $\Theta \subseteq \R^d$ be open and let $\mathcal{D}=(Y_1,\ldots,Y_N)$ consist of independent records satisfying
\[
Y_t \sim \mathcal{N}\{\mu_t(\theta), \sigma^2 I_q\}, 
\quad t\in [N],
\]
where $\theta\in \Theta$, $\sigma^2>0$, and each $\mu_t: \Theta \mapsto \R^q$ is continuously differentiable. 
Write $B_t(\theta) :=D_\theta \mu_t(\theta) \in\R^{q\times d}$ for its Jacobian.
Suppose that for some $\lambda>0$ and absolute constants $C\geq c>0$,
uniformly over $\theta\in\Theta$, that
\begin{equation}
\label{eq:score-contraction-jacobian-bounds}
cN\lambda I_d \preceq 
\sum_{t=1}^N B_t(\theta)^\top B_t(\theta) 
\preceq
CN\lambda I_d 
\quad \mbox{and} \quad 
\max_{t\in[N]} \|B_t(\theta)\|^2 \leq C\lambda .
\end{equation}

Let $Z= \mathscr{M}(\mathcal{D})$, where $\mathscr{M}$ is a parameter-independent $(\varepsilon, \delta)$-centrally differentially private mechanism under replacement of one record, taking values in a standard Borel space. 
Suppose that there exists an absolute constant $\varepsilon_0, c_0 >0$ such that $0<\varepsilon \leq \varepsilon_0$, $0<\delta \leq \exp(-1)$, and
\begin{equation}
\label{eq:score-contraction-delta-condition}
d\delta \log(1/\delta) \leq 
c_0 \sqrt{N}\varepsilon^2.
\end{equation}
Let $I_Z(\theta)$ denote the Fisher information matrix for $\theta$ in the induced distribution of $Z$. 
Then uniformly over $\theta\in\Theta$, there exists an absolute constant $C'>0$ such that  
\begin{equation}
\label{eq:score-contraction-bound}
\mathrm{tr}\big\{ I_Z(\theta) \big\}
\leq 
C' \frac{\lambda}{\sigma^2} \big\{ (N d) \wedge (N^2 \varepsilon^2)
\big\}.
\end{equation}
\end{lemma}

\begin{proof}[Proof of Lemma~\ref{lem:score-contraction}]
We divide the proof into four steps.

\medskip\noindent\textbf{Step 1. Complete-data information.}
Let $p_{t,\theta}$ be the density of $Y_t$. Its score is
\[
S_t(\theta) := \nabla_\theta \log(p_{t,\theta}(Y_t)) = \frac{1}{\sigma^2} B_t(\theta)^\top \{Y_t -\mu_t(\theta)\} .
\]
Thus, $S_t(\theta)$ is centered Gaussian with covariance $\sigma^{-2} B_t(\theta)^\top B_t(\theta)$. 
Since the records are independent, the complete-data score is $\sum_{t=1}^N S_t(\theta)$ and its Fisher information is
\begin{equation}\label{eq-corollary-1.1}
J(\theta) := \mathbb{E}_\theta\bigg\{ 
\bigg( \sum_{t=1}^N S_t(\theta) \bigg)
\bigg( \sum_{t=1}^N S_t(\theta) \bigg)^\top 
\bigg\} 
= 
\frac{1}{\sigma^2} \sum_{t=1}^N B_t(\theta)^\top B_t(\theta),
\end{equation}
where the last equality follows from the fact that the record scores are
independent and centered. 
Combining \eqref{eq:score-contraction-jacobian-bounds} and \eqref{eq-corollary-1.1}, we have
\begin{equation}\label{eq-corollary-1.2}
c N \frac{\lambda}{\sigma^2} I_d 
\preceq 
J(\theta) 
\preceq 
C N \frac{\lambda}{\sigma^2} I_d.
\end{equation}

\medskip\noindent\textbf{Step 2. Score of the private output.}
Let $p_\theta(\mathbf{y}) =\prod_{t=1}^N p_{t,\theta}(y_t)$. 
Write $k(z \mid \mathbf{y})$ for the conditional density or mass function of the output given $\mathcal{D} =\mathbf{y}$, and let $q_\theta$ be the density or mass function of $Z$. 
Since the mechanism does not depend on $\theta$, neither does $k(z \mid \mathbf{y})$, and
\[
q_\theta(z) = \int k(z \mid \mathbf{y}) p_\theta(\mathbf{y}) d\mathbf{y} .
\]
Differentiating under the integral gives that
\begin{align}\label{eq-corollary-2.1}
    \nabla_\theta \log(q_\theta(z)) &= 
    \frac{1}{q_\theta(z)}
    \int k(z \mid \mathbf{y}) \nabla_\theta p_\theta(\mathbf{y}) d\mathbf{y}
    \nonumber\\
    &=
    \int \sum_{t=1}^N \nabla_\theta \log(p_{t,\theta}(y_t))  \frac{k(z \mid \mathbf{y}) p_\theta(\mathbf{y})}{q_\theta(z)} d\mathbf{y}
    =
    \mathbb{E}_\theta\bigg\{ \sum_{t=1}^N S_t(\theta) \bigg| Z=z \bigg\},
\end{align}
where the second equality follows from 
\[
\nabla_\theta p_\theta(\mathbf{y}) = p_\theta(\mathbf{y}) \sum_{t=1}^N \nabla_\theta \log(p_{t,\theta}(y_t)) ,
\]
and the last follows since $k(z \mid \mathbf{y}) p_\theta(\mathbf{y}) /q_\theta(z)$ is the conditional density of $\mathcal{D}$ given $Z=z$. 
Consequently, \eqref{eq-corollary-2.1} gives
\begin{equation}\label{eq-corollary-2.2}
I_Z(\theta) = \mathbb{E}_\theta\bigg[ 
\mathbb{E}_\theta\bigg\{ \sum_{t=1}^N S_t(\theta) \bigg| Z\bigg\} 
\mathbb{E}_\theta\bigg\{ \sum_{t=1}^N S_t(\theta) \bigg| Z\bigg\}^{\top} 
\bigg].
\end{equation}

\medskip\noindent\textbf{Step 3. Private-score contraction.} 
We apply Corollary~B.4 in  \citet{cai2026van} to bound the normalized Fisher information $\operatorname{tr}\{J(\theta)^{-1} I_Z(\theta)\}$ retained in the privatized output $Z$.

Note that \eqref{eq-corollary-1.2} implies $J(\theta)\succ 0$. Since the complete-data score is centered Gaussian with covariance $J(\theta)$,
\[
J(\theta)^{-1/2} \sum_{t=1}^N S_t(\theta) \sim \mathcal{N}(0,I_d) .
\]
For every $\|v\|_2=1$, we have
\[
\mathrm{Var}_\theta \{v^\top J(\theta)^{-1/2} S_t(\theta)\} = 
\frac{1}{\sigma^2} v^\top J(\theta)^{-1/2} B_t(\theta)^\top B_t(\theta) J(\theta)^{-1/2} v 
\leq \frac{C}{cN},
\]
where the inequality follows from \eqref{eq:score-contraction-jacobian-bounds} and \eqref{eq-corollary-1.2}. 
Since this random variable is Gaussian, its sub-Gaussian norm is bounded by a constant multiple of its standard deviation. 
Hence there exists an absolute constant $C_1>0$ such that 
\begin{equation}\label{eq-corollary-b4-1}
\sup_{\|v\|_2=1} \bigg\|v^\top J(\theta)^{-1/2} \sum_{t=1}^N  S_t(\theta) \bigg\|_{\psi_2}
\leq C_1, 
\quad  
\sup_{\|v\|_2=1} \big\|v^\top J(\theta)^{-1/2} S_t(\theta) \big\|_{\psi_2} 
\leq  \frac{C_1}{\sqrt{N}}, 
\quad 
\forall t \in [N].
\end{equation}

By Corollary~B.4 in \citet{cai2026van}, together with \eqref{eq-corollary-b4-1} and $\mathrm{rank}\{J(\theta)\}=d$, there exists an absolute constant $C_2 >0$ such that
\begin{align}\label{eq-corollary-3.1}
\mathrm{tr}\{J(\theta)^{-1} I_Z(\theta)\}
\leq
C_2\big\{ N(e^\varepsilon -1)^2 + C_1\delta d \sqrt{N} \log(1/\delta) \big\}.
\end{align}
Since $0< \varepsilon \leq \varepsilon_0$, we have 
\begin{equation}\label{eq-corollary-3.2}
(e^\varepsilon -1)^2 = \bigg(\int_0^\varepsilon e^t dt \bigg)^2
\leq
e^{2\varepsilon_0}\varepsilon^2.
\end{equation}
Combining \eqref{eq-corollary-3.1} and \eqref{eq-corollary-3.2} implies
\begin{equation}\label{eq-corollary-3.3}
\mathrm{tr}\{J(\theta)^{-1} I_Z(\theta)\}
\leq 
C_2\big\{ N e^{2\varepsilon_0} \varepsilon^2 + C_1\delta d\sqrt{N} \log(1/\delta) \big\} \leq 
C_3 N \varepsilon^2 ,
\end{equation}
where $C_3 = C_2 e^{2\varepsilon_0} + C_2  C_1 c_0$ is an absolute constant and the last inequality follows from \eqref{eq:score-contraction-delta-condition}.

\medskip\noindent\textbf{Step 4. Completion of the proof.}
Since $I_Z(\theta) \succeq 0$, $J(\theta)^{-1/2} I_Z(\theta)J(\theta)^{-1/2} \succeq 0$. 
Hence combining \eqref{eq-corollary-1.2} and \eqref{eq-corollary-3.3}, we have
\begin{align}
\label{eq-corollary-4.1}
    \mathrm{tr}\{I_Z(\theta) \} &=
    \mathrm{tr}\{ J(\theta)J(\theta)^{-1/2} I_Z(\theta)J(\theta)^{-1/2} \} \nonumber\\
    &\leq 
    \lambda_1\{J(\theta)\} \mathrm{tr}\{J(\theta)^{-1} I_Z(\theta)\}
    \leq
    C C_3 \frac{N^2 \varepsilon^2 \lambda}{\sigma^2}.
\end{align}

For the non-private bound, combining  \eqref{eq-corollary-1.1} and \eqref{eq-corollary-2.2}, we have
\begin{align*}
    \mathbb{E}_\theta\bigg[ \mathrm{Var}_\theta \bigg\{ \sum_{t=1}^N S_t(\theta) \bigg| Z \bigg\}\bigg]
    &=
    \mathbb{E}_\theta\bigg\{ \bigg(\sum_{t=1}^N S_t(\theta)\bigg) \bigg(\sum_{t=1}^N S_t(\theta)\bigg)^\top \bigg\} \\
    &\;\quad 
    - \mathbb{E}_\theta\bigg[ \mathbb{E}_\theta\bigg\{ \sum_{t=1}^N S_t(\theta) \bigg| Z \bigg\} 
    \mathbb{E}_\theta\bigg\{ \sum_{t=1}^N S_t(\theta) \bigg| Z \bigg\}^{\top} \bigg] \\
    &=
    J(\theta) - I_Z(\theta) 
    \succeq 0.
\end{align*}
Thus, it follows that 
\begin{equation}\label{eq-corollary-4.2}
\mathrm{tr}\{I_Z(\theta)\} 
\leq 
\mathrm{tr}\{J(\theta)\} 
\leq 
C \frac{N d \lambda}{\sigma^2} ,
\end{equation}
where the last inequality follows from \eqref{eq-corollary-1.2}. 
Combining \eqref{eq-corollary-4.1} and \eqref{eq-corollary-4.2}, we conclude \eqref{eq:score-contraction-bound},
which completes this proof.
\end{proof}

\subsection{Proofs of theorems and corollaries}

\begin{proof}[Proof of Theorem~\ref{thm:fed-upper}]
Let $P:= A_k A_k^\top$ and let $\breve\Pi_{j,k}$ be the projector constructed in Lemma~\ref{lem:local-ideal}. 
Recall from \eqref{eq:Psi-local} that $w_{j,k} \geq 0$ and $\sum_{j=1}^m w_{j,k} =1$. 
Define the population and random aggregation matrices as
\begin{equation}
\label{eq:ideal-server-matrices}
M_k :=
\sum_{j=1}^m w_{j,k} \E(\breve\Pi_{j,k}),
\quad
\breve{M}_k :=
\sum_{j=1}^m w_{j,k} \breve\Pi_{j,k}.
\end{equation}
We first identify the leading eigenspace of $M_k$. By Lemma~\ref{lem:local-ideal}~(iii), we have
\begin{align}
(I_{p_k} -P) M_k P &=
\sum_{j=1}^m w_{j,k} (I_{p_k} -P)\E(\breve\Pi_{j,k}) P
= 0 .
\label{eq:server-population-cross-block}
\end{align}
Since $M_k$ is symmetric, \eqref{eq:server-population-cross-block} also gives $P M_k (I_{p_k} -P) =0$. Hence $M_k$ is block diagonal with respect to $\mathrm{span}(A_k) \oplus \mathrm{span}(A_k)^\perp$. Also, by Lemma~\ref{lem:local-ideal}~(iv), we have
\begin{equation}
\label{eq:server-population-close}
\|M_k-P\| \leq
\sum_{j=1}^m w_{j,k} \|\E(\breve\Pi_{j,k}) -P\|
\leq c .
\end{equation}
For every unit vector $u\in \mathrm{span}(A_k)$ and every unit vector $v\in \mathrm{span}(A_k)^\perp$, \eqref{eq:server-population-close} gives $u^\top M_ku \geq 1-c$ and $v^\top M_k v \leq c$.
Consequently, if the absolute constant $c$ in Lemma~\ref{lem:local-ideal}~(iv) is sufficiently small, the leading $r_k$-dimensional eigenspace of $M_k$ is exactly $\mathrm{span}(A_k)$ and
\begin{equation}
\label{eq:server-population-gap}
\lambda_{r_k}(M_k) -\lambda_{r_k+1}(M_k)
\geq 1 -2c
\geq \frac{1}{2} .
\end{equation}

Consider $\breve{M}_k -M_k$ next. 
For $j\in[m]$, set $Y_{j,k} := w_{j,k} \big\{\breve\Pi_{j,k}-\E(\breve\Pi_{j,k}) \big\}$.
Then by Lemma~\ref{lem:local-ideal}~(i), the matrices $Y_{1,k},\ldots, Y_{m,k}$ are independent, and $\E(Y_{j,k})=0$. 
Lemma~\ref{lem:local-ideal}~(v) gives the almost-sure bound such that
\begin{equation}
\label{eq:weighted-ideal-bound}
\|Y_{j,k}\|_{\mathrm{F}}
\leq
C w_{j,k} \sqrt{r_k} \Psi_{j,k}^{\mathrm{loc}} .
\end{equation}
By similar arguments for \eqref{eq:Hilbert-Hoeffding-local} and \eqref{eq:centered-block-average}, applying Theorem~3.5 of \cite{pinelis1994optimum} in the space of symmetric matrices equipped with the Frobenius norm to \eqref{eq:weighted-ideal-bound}, for every $x>0$, we obtain
\begin{equation}
\label{eq:server-Hilbert-Hoeffding}
\P\Bigg\{
\|\breve{M}_k -M_k\|_{\mathrm{F}}^2
>
C r_k x \sum_{j=1}^m w_{j,k}^2 (\Psi_{j,k}^{\mathrm{loc}})^2
\Bigg\}
\leq
2 e^{-x} .
\end{equation}
By the definition of the weights in \eqref{eq:Psi-local}, we have
\begin{align}
    \sum_{j=1}^m w_{j,k}^2 (\Psi_{j,k}^{\mathrm{loc}})^2
    &=
    \frac{ \sum_{j=1}^m (\Psi_{j,k}^{\mathrm{loc}})^{-4}(\Psi_{j,k}^{\mathrm{loc}})^2}{\big\{ \sum_{j=1}^m (\Psi_{j,k}^{\mathrm{loc}})^{-2} \big\}^2}
    =
    \frac{1}{\sum_{j=1}^m (\Psi_{j,k}^{\mathrm{loc}})^{-2}} .
    \label{eq:inverse-square-weight-identity}
\end{align}
Setting $x=\log(n\vee p)$ in \eqref{eq:server-Hilbert-Hoeffding}, together with \eqref{eq:inverse-square-weight-identity}, it then holds with probability at least $1-2(n\vee p)^{-1}$ that
\begin{equation}
\label{eq:ideal-server-concentration}
\|\breve{M}_k -M_k\|_{\mathrm{F}}^2
\leq
\frac{C r_k \log(n\vee p)}
{\sum_{j=1}^m (\Psi_{j,k}^{\mathrm{loc}})^{-2}} .
\end{equation}
Let $\breve{A}_k:= \mathrm{SVD}_{r_k}(\breve{M}_k)$ and $\breve{P}_k := \breve{A}_k \breve{A}_k^\top$. 
The population eigengap in \eqref{eq:server-population-gap} and the Davis--Kahan theorem \citep[Theorem~2]{yu2015useful} imply
\begin{equation}
\label{eq:ideal-server-Davis-Kahan}
\|\breve{P}_k -P\|_{\mathrm{F}}^2
\leq
C \|\breve{M}_k -M_k\|_{\mathrm{F}}^2.
\end{equation}
Combining \eqref{eq:ideal-server-concentration} and \eqref{eq:ideal-server-Davis-Kahan}, 
it holds with probability at least $1 -2(n\vee p)^{-1}$ that
\begin{equation}
\label{eq:ideal-server-projector-rate}
\|\breve{P}_k -P\|_{\mathrm{F}}^2
\leq
\frac{C r_k \log(n\vee p)}{\sum_{j=1}^m (\Psi_{j,k}^{\mathrm{loc}})^{-2}} .
\end{equation}
In addition, note that for any two rank-$r_k$ orthogonal projectors $P_1$ and $P_2$, we have $\|P_1 -P_2\|_{\mathrm{F}}^2 = 2 r_k -2\mathrm{tr}(P_1 P_2) \leq 2 r_k$.
Hence the right-hand side of \eqref{eq:ideal-server-projector-rate} may be replaced by its minimum with $2 r_k$. It remains to connect the above result to our procedure in Algorithm~\ref{alg:fed-est}. 
To this end, define the matrix aggregated by the server as $\widetilde{M}_k := \sum_{j=1}^m w_{j,k} \widetilde{U}_{j,k} \widetilde{U}_{j,k}^\top$.
By Lemma~\ref{lem:local-ideal}~(ii), it holds with probability at least $1-C(n\vee p)^{-c}$ that 
\[
\breve\Pi_{j,k} =
\widetilde{U}_{j,k} \widetilde{U}_{j,k}^\top ,
\quad
\forall j\in[m] . 
\]
Note that on this event, $\breve{M}_k=\widetilde{M}_k$. 
Hence by the fixed tie-breaking rule stated in Section~\ref{subsec:notation},
$\breve{P}_k = \widehat{A}_k^{\mathrm{fed}} (\widehat{A}_k^{\mathrm{fed}})^\top$.
Combining this event with \eqref{eq:ideal-server-projector-rate} by a union bound, we obtain that, with probability at least $1-2(n\vee p)^{-1}-C(n\vee p)^{-c}$,
\[
\big\|\widehat{A}_k^{\mathrm{fed}} (\widehat{A}_k^{\mathrm{fed}})^\top -A_k A_k^\top \big\|_{\mathrm{F}}^2
\leq
\frac{C r_k \log(n\vee p)}{\sum_{j=1}^m (\Psi_{j,k}^{\mathrm{loc}})^{-2}}
\wedge (2r_k) ,
\]
which is \eqref{eq:fed-upper}.
This completes the proof of this theorem.
\end{proof}

\begin{proof}[Proof of Theorem~\ref{thm:fed-lower-target}]
Our proof proceeds in three steps. 
First, we construct an $r_k (p_k -r_k)$-dimensional Gaussian sub-model in which only the mode-$k$ loading space varies. 
Second, we bound the Fisher information retained in each private client message and aggregate these bounds over all clients. 
Finally, we apply the multivariate van Trees inequality and derive the desired minimax lower bound.

\medskip\noindent\textbf{Step 1.} 
In this step, we construct a $r_k (p_k -r_k)$-dimensional Gaussian sub-model that remains inside $\mathfrak{P}_k^{\mathrm{fed}}(\lambda, \sigma^2)$.

For each $h\in [K] \setminus\{k\}$, fix $A_h\in \mathbb{O}_{p_h\times r_h}$ satisfying Assumption~\ref{ass:loading}, 
and recall $A_{-k} =\bigotimes_{h\in[K] \setminus\{k\}} A_h$. 
Since $A_h^\top A_h=I_{r_h}$ for every $h\in [K]$, we have
\[
A_{-k}^\top A_{-k}  = \bigotimes_{h\in[K] \setminus\{k\}} A_h^\top A_h = I_{r_{-k}}.
\]

We next construct a local family of mode-$k$ loading matrices. For any $H \in \R^{(p_k-r_k) \times r_k}$, define
\begin{equation}\label{eq:lower-1.1.1}
A_k(H) := 
\begin{pmatrix}
I_{r_k}\\ H 
\end{pmatrix} 
\big( I_{r_k} + H^\top H \big)^{-1/2},
\end{equation}
Since $A_k(H)^\top A_k(H) =I_{r_k}$, we have that $A_k(H)\in \mathbb{O}_{p_k\times r_k}$ for every $H\in \R^{(p_k -r_k)\times r_k}$. Let $P_k(H)= A_k(H)A_k(H)^\top$ and $d_k= r_k (p_k -r_k)$, where $d_k$ is the number of free coordinates in $H$.

\medskip\noindent\textbf{Step 1.1.} 
In this sub-step, we show that there exist absolute constants  $C_1' \geq c_1' >0$ such that, 
for all $H,H' \in\R^{(p_k -r_k)\times r_k}$ satisfying $\|H\|_{\mathrm{F}}, \|H'\|_{\mathrm{F}} \leq 1$,
\begin{equation}\label{eq:lower-1.1.2}
c_1' \|H-H'\|_{\mathrm{F}} 
\leq
\|P_k(H) - P_k(H')\|_{\mathrm{F}} 
\leq 
C_1' \|H-H'\|_{\mathrm{F}} .
\end{equation}

Fix $H,H'\in\R^{(p_k-r_k)\times r_k}$ with $\|H\|_{\mathrm{F}}, \|H'\|_{\mathrm{F}} \leq 1$. Let
\[
X(H) := 
\begin{pmatrix}
I_{r_k}\\
H
\end{pmatrix}
\quad \mbox{and} \quad 
R(H) := \big( I_{r_k} + H^\top H \big)^{-1/2}.
\]
Then $A_k(H)= X(H) R(H)$ and
\begin{equation}\label{eq:lower-def-pk-h}
P_k(H) =
\begin{pmatrix}
R(H)^2 & R(H)^2 H^\top\\
H R(H)^2 & H R(H)^2 H^\top
\end{pmatrix}.
\end{equation}
Since $R(H)^2= (I_{r_k} + H^\top H)^{-1}$ and $R(H')^2 =(I_{r_k} +H'^\top H')^{-1}$,
\begin{equation}\label{eq:lower-1.1.3}
\|R(H)^{-2}\| \leq 2, \quad
\|R(H')^{-2}\|\leq 2, \quad 
\|R(H)^2\|\leq 1, \quad
\|H'R(H')^2\| \leq \|H'\| \|R(H')^2\| \leq 1 .
\end{equation}
Moreover, 
\begin{align}\label{eq:lower-1.1.4}
    \|R(H)^{-2} -R(H')^{-2}\|_{\mathrm{F}}  
    &= 
    \|R(H)^{-2} \{R(H')^2 -R(H)^2\} R(H')^{-2}\|_{\mathrm{F}}  \nonumber\\
    &\leq
    \|R(H)^{-2}\| \| R(H')^2 -R(H)^2\|_{\mathrm{F}} \|R(H')^{-2}\| \nonumber\\
    &\leq 
    4 \| R(H')^2 -R(H)^2\|_{\mathrm{F}} ,
\end{align} 
where the last inequality follows from \eqref{eq:lower-1.1.3}.
Therefore, we have
\begin{align}\label{eq:lower-1.1.5}
    \|H-H'\|_{\mathrm{F}} 
    &=
    \big\| \{HR(H)^2 -H'R(H')^2\} R(H)^{-2} 
    + H'R(H')^2 \{R(H)^{-2} -R(H')^{-2}\} \big\|_{\mathrm{F}} \nonumber\\
    &\leq
    \|R(H)^{-2} \| \|HR(H)^2-H'R(H')^2\|_{\mathrm{F}}
    + \| H'R(H')^2 \| \|R(H)^{-2} -R(H')^{-2}\|_{\mathrm{F}} \nonumber\\
    &\leq 
    2 \|HR(H)^2 -H' R(H')^2\|_{\mathrm{F}} 
    + 4 \|R(H)^{2} -R(H')^{2}\|_{\mathrm{F}}
    \leq 
    6\|P_k(H) -P_k(H')\|_{\mathrm{F}} ,
\end{align}
where the second inequality follows from \eqref{eq:lower-1.1.3} and \eqref{eq:lower-1.1.4}, and the last inequality follows from~\eqref{eq:lower-def-pk-h}.

For any positive definite matrices $U,V \succeq I_{r_k}$, the integral representation of the inverse square root gives
\begin{align}\label{eq:lower-1.1.6}
    \|U^{-1/2}-V^{-1/2}\|_{\mathrm{F}}
    &=
    \bigg\| \frac{1}{\pi} \int_0^\infty t^{-1/2} (U + tI_{r_k})^{-1} dt 
    - \frac{1}{\pi} \int_0^\infty t^{-1/2}(V + tI_{r_k})^{-1} dt \bigg\|_{\mathrm{F}} \nonumber\\
    &=
    \bigg\| \frac{1}{\pi} \int_0^\infty t^{-1/2} (U + tI_{r_k})^{-1} (V-U)(V + tI_{r_k})^{-1}dt  \bigg\|_{\mathrm{F}} \nonumber\\
    &\leq
    \frac{1}{\pi} \int_0^\infty t^{-1/2} \| (U + tI_{r_k})^{-1} \| \| V-U\|_{\mathrm{F}} \| (V + tI_{r_k})^{-1}\| dt  \nonumber\\
    &\leq
    \frac{1}{\pi} \int_0^\infty \frac{t^{-1/2}}{(1+t)^2}dt \|U-V\|_{\mathrm{F}} 
    =  
    \frac{1}{2} \|U-V\|_{\mathrm{F}} .
\end{align}
where the second inequality follows from $U,V \succeq I_{r_k}$, which implies that
\[
\|(U + tI_{r_k})^{-1}\| \leq 
\frac{1}{1+t} 
\quad \mbox{and} \quad
\|(V + tI_{r_k})^{-1}\| 
\leq \frac{1}{1+t}.
\]

Applying \eqref{eq:lower-1.1.6} with $U = I_{r_k} + H^\top H$ and $V = I_{r_k} + H'^\top H'$ gives
\begin{align}\label{eq:lower-1.1.7}
    \|R(H)-R(H')\|_{\mathrm{F}}
    &\leq 
    \frac{1}{2} \|H^\top H-H'^\top H'\|_{\mathrm{F}} 
    = 
    \frac{1}{2} \|H^\top(H-H') + (H-H')^\top H'\|_{\mathrm{F}} \nonumber\\
    &\leq 
    \frac{1}{2} \{\|H\|+\|H'\|\} \|H-H'\|_{\mathrm{F}} 
    \leq
    \|H-H'\|_{\mathrm{F}} ,
\end{align}
where the last inequality follows from $\|H\|, \|H'\| \leq 1$.
Moreover, since
\[
X(H)-X(H') =\begin{pmatrix}
0\\
H-H'
\end{pmatrix},
\]
we have
\begin{align}
    & \|X(H)-X(H')\|_{\mathrm{F}} = \|H-H'\|_{\mathrm{F}} ,
    \label{eq:lower-1.1.8} \\
    & \|X(H')\|^2 = \|X(H')^\top X(H')\| = \|I_{r_k} + H'^\top H'\| 
    \leq 2 .
    \label{eq:lower-1.1.9}
\end{align}
Hence altogether, we arrive at
\begin{align}\label{eq:lower-1.1.10}
    \| P_k(H)-P_k(H') \|_{\mathrm{F}}
    &=
    \big\| \{A_k(H)-A_k(H')\} A_k(H)^\top + A_k(H') \{A_k(H)-A_k(H')\}^\top \big\|_{\mathrm{F}} \nonumber\\
    &\leq
    \big(\| A_k(H)\| + \| A_k(H')\|\big)  \| A_k(H)-A_k(H')\|_{\mathrm{F}} 
    = 
    2\| A_k(H)-A_k(H')\|_{\mathrm{F}}  \nonumber\\ 
    &=
    2 \big\| \{X(H)-X(H')\}R(H) + X(H')\{R(H)-R(H')\} \big\|_{\mathrm{F}} \nonumber\\
    &\leq 
    2 \|X(H)-X(H')\|_{\mathrm{F}} \|R(H)\| + 2\|X(H')\|\|R(H)-R(H')\|_{\mathrm{F}}
    \nonumber\\
    &\leq
    2 (1+\sqrt{2}) \|H-H'\|_{\mathrm{F}} ,
\end{align}
where the second equality follows from $A_k(H), A_k(H')\in \mathbb{O}_{p_k\times r_k}$, and the last inequality follows from \eqref{eq:lower-1.1.7}, \eqref{eq:lower-1.1.8}, \eqref{eq:lower-1.1.9} and $\| R(H) \|= \|(I_{r_k} + H^\top H)^{-1/2}\| \leq 1$.
Combining \eqref{eq:lower-1.1.5} and \eqref{eq:lower-1.1.10}, we complete the proof of \eqref{eq:lower-1.1.2} with constants $c_1' = 1/6$ and $C_1' = 2 (1 +\sqrt{2})$.

\medskip\noindent\textbf{Step 1.2.} For $H,G\in \R^{(p_k-r_k)\times r_k}$, define the directional derivatives
\[
\dot{A}_k(H)[G] := 
\frac{d}{ds} A_k(H + sG) \bigg|_{s=0} 
\quad \mbox{and} \quad 
\dot{P}_k(H)[G] := \frac{d}{ds} P_k(H + sG) \bigg|_{s=0}.
\]
In this sub-step, we show that there exist absolute constants $C_2'\geq c_2' >0$ such that, 
for every $\|H\|_{\mathrm{F}} \leq 1/2$ and every $G\in \R^{(p_k-r_k)\times r_k}$,
\begin{equation}\label{eq:lower-1.2.1}
c_2' \|G\|_{\mathrm{F}} 
\leq 
\|\dot{A}_k(H)[G]\|_{\mathrm{F}} 
\leq C_2' \|G\|_{\mathrm{F}} .
\end{equation}

Fix $\|H\|_{\mathrm{F}} \leq 1/2$ and $G\in \R^{(p_k-r_k)\times r_k}$. For all sufficiently small $|s|$, both $H$ and $H + sG$ belong to the unit ball. Therefore,
\eqref{eq:lower-1.1.2} gives
\[
c_1' |s| \|G\|_{\mathrm{F}} 
\leq
\|P_k(H + sG) -P_k(H)\|_{\mathrm{F}} .
\]
Moreover, the calculation leading to \eqref{eq:lower-1.1.10} gives
\[
\|A_k(H + sG) -A_k(H)\|_{\mathrm{F}} 
\leq
(1+\sqrt{2}) |s| \|G\|_{\mathrm{F}} .
\]
Dividing these inequalities by $|s|$ and letting $s\to 0$ gives 
\begin{equation}\label{eq:lower-1.2.2}
c_1' \|G\|_{\mathrm{F}} 
\leq
\|\dot{P}_k(H)[G]\|_{\mathrm{F}} ,
\quad \mbox{and} \quad 
\|\dot{A}_k(H)[G]\|_{\mathrm{F}} 
\leq 
(1+\sqrt{2}) \|G\|_{\mathrm{F}} .
\end{equation}

Differentiating $P_k(H) = A_k(H) A_k(H)^\top$ in the direction~$G$ gives
\[
\dot{P}_k(H)[G] = \dot{A}_k(H)[G] A_k(H)^\top + A_k(H) \dot{A}_k(H)[G]^\top .
\]
Since $\|A_k(H)\|=1$, we have
\begin{equation}\label{eq:lower-1.2.3}
\|\dot{P}_k(H)[G]\|_{\mathrm{F}} 
\leq 
2 \|\dot{A}_k(H)[G]\|_{\mathrm{F}} .
\end{equation}
Combining \eqref{eq:lower-1.2.2} and \eqref{eq:lower-1.2.3}, we have 
\begin{equation}\label{eq:lower-1.2.4}
\frac{c_1'}{2} \|G\|_{\mathrm{F}} 
\leq
\|\dot{A}_k(H)[G]\|_{\mathrm{F}} 
\leq 
(1+\sqrt{2}) \|G\|_{\mathrm{F}} .
\end{equation}
Thus, \eqref{eq:lower-1.2.1} holds with $c_2'= c_1'/2$ and $C_2' = 1+\sqrt{2}$.

\medskip\noindent\textbf{Step 1.3.}
In this sub-step, we construct core factors.
Fix $j\in[m]$. We construct its core factors as follows. For each $\ell\in [n_j r_{-k}]$, define
\[
i_\ell := 1+ \big\{(\ell-1) \bmod r_k\big\}
,\quad 
a_i := \big|\{\ell\in [n_j r_{-k}]: i_\ell=i\} \big|
,\quad i\in[r_k].
\]
The indices $i_\ell$ cycle through $[r_k]$, and hence
\[
a_i \in \bigg\{ \bigg\lfloor \frac{n_j r_{-k}}{r_k} \bigg\rfloor, 
\bigg\lceil \frac{n_j r_{-k}}{r_k} \bigg\rceil \bigg\}.
\]
Since $n_j r_{-k} \geq r_k$, we have
\begin{equation}\label{eq:lower-1.3.1}
a_i \geq \bigg\lfloor \frac{n_j r_{-k}}{r_k} \bigg\rfloor 
\geq \frac{n_j r_{-k}}{2r_k} .
\end{equation}
Let $e_i$ denote the $i$th canonical basis vector of $\R^{r_k}$ and define
\[
u_\ell := \bigg( \frac{n_j\lambda}{a_{i_\ell}} \bigg)^{1/2} e_{i_\ell}
, \quad \forall \ell\in [n_j r_{-k}].
\]
For each $j\in[m]$ and $t\in [n_j]$, define
\[
D_{j,t} := \big[ u_{(t-1) r_{-k} +1} , \ldots, u_{t r_{-k}} \big] 
\in \R^{r_k\times r_{-k}} .
\]

Since each coordinate $i\in [r_k]$ occurs exactly $a_i$ times, we have
\begin{equation}\label{eq:lower-1.3.2}
\frac{1}{n_j} \sum_{t=1}^{n_j} D_{j,t} D_{j,t}^\top
=
\frac{1}{n_j} \sum_{\ell=1}^{n_j r_{-k}} u_\ell u_\ell^\top 
=
\frac{1}{n_j} \sum_{i=1}^{r_k} a_i \frac{n_j\lambda}{a_i} e_i e_i^\top 
=
\lambda I_{r_k} .
\end{equation}
Moreover, for every $t\in[n_j]$, we obtain
\begin{equation}\label{eq:lower-1.3.3}
\|D_{j,t}\|^2 \leq
\|D_{j,t}\|_{\mathrm{F}}^2 
=
\sum_{\ell= (t-1) r_{-k} +1}^{t r_{-k}}  \|u_\ell\|_2^2 
=
\sum_{\ell =(t-1) r_{-k} +1}^{t r_{-k}}
\frac{n_j\lambda}{a_{i_\ell}} 
\leq 
\frac{n_j r_{-k} \lambda}{\min_{i\in[r_k]} a_i} 
\leq 
2 r_k \lambda = C_r\lambda ,
\end{equation}
where the last inequality follows from \eqref{eq:lower-1.3.1}, and $C_r = 2 r_k >0$ is a constant since the factor ranks are fixed.

For each $j\in [m]$ and $t\in [n_j]$, let $\mathbf{F}_t^{(j)}$ be the core tensor satisfying $\mathcal{M}_k(\mathbf{F}_t^{(j)}) = D_{j,t}$.  
Since mode-$k$ matricization is a bijection, this core tensor is uniquely determined. For each $H$ satisfying $\|H\|_{\mathrm{F}} < 1/2$, define
\[
\mathbf{S}_t^{(j)} := \mathbf{F}_t^{(j)} \times_1 A_1 \times_2 \cdots \times_k A_k(H) \times_{k + 1} \cdots \times_K A_K, 
\quad 
\mathbf{X}_t^{(j)} := \mathbf{S}_t^{(j)} + \mathbf{E}_t^{(j)},
\]
where the entries of $\mathbf{E}_t^{(j)}$ are independent $\mathcal{N}(0, \sigma^2)$ random variables.
For notational convenience,
write $Y_{j,t} := \mathcal{M}_k(\mathbf{X}_t^{(j)})$ and $E_{j,t} := \mathcal{M}_k(\mathbf{E}_t^{(j)})$. 
Then we have 
\begin{equation}\label{eq:lower-1.3.4}
Y_{j,t} = A_k(H) D_{j,t} A_{-k}^\top + E_{j,t}.
\end{equation}
Moreover, \eqref{eq:lower-1.3.2} and \eqref{eq:lower-1.3.3} imply
$\Lambda_{j,(k)}=\lambda I_{r_k}$ and $\max_t \|D_{j,t}\|^2 \leq C_r\lambda$. 
Hence the constructed model is a subclass of $\mathfrak{P}_k^{\mathrm{fed}}(\lambda, \sigma^2)$.

\medskip\noindent\textbf{Step 2.}
In this step, we first fix a client $j \in[m]$, verify the conditions of Lemma~\ref{lem:score-contraction}, and apply the lemma to bound the Fisher information in its private message. We then aggregate these information bounds over all clients.

\medskip\noindent\textbf{Step 2.1.}
In this sub-step, we derive the score and Fisher information of one
record.
Fix $j\in[m]$. Since Lemma~\ref{lem:score-contraction} is formulated for a vector-valued parameter, for every $H \in\R^{(p_k-r_k)\times r_k}$, write $\theta := \mathrm{vec}(H) \in\R^{d_k}$. 
Note that $\|\theta\|_2 = \|H\|_{\mathrm{F}}$. For every $t\in[n_j]$, define the mean vector
\[
\mu_{j,t}(H) := \mathrm{vec} \big\{ A_k(H) D_{j,t} A_{-k}^\top \big\} 
\in\R^{p_k p_{-k}} .
\]
It follows from \eqref{eq:lower-1.3.4} that
\begin{equation}\label{eq:lower-2.1.1}
\mathrm{vec}(Y_{j,t}) \sim 
\mathcal{N}\big\{\mu_{j,t}(H), \sigma^2 I_{p_k p_{-k}}\big\}, 
\quad 
\forall t \in [n_j] ,
\end{equation}
and these Gaussian records are independent across $t$.

Define the Jacobian of $\mu_{j,t}(H)$ with respect to $\theta = \mathrm{vec}(H)$ by
\[
B_{j,t}(H) = D_\theta \mu_{j,t}(H) =
\bigg( 
\frac{\partial \mu_{j,t,a}(H)}{\partial \theta_b}
\bigg)_{a \in [p_k p_{-k}], b\in[d_k]}
\in \R^{p_kp_{-k}\times d_k} .
\]
For any $G\in \R^{(p_k-r_k)\times r_k}$, let
$g:= \mathrm{vec}(G)$. By the definition of the directional derivative, we have
\begin{align}\label{eq:lower-2.1.2}
B_{j,t}(H)g = \frac{d}{ds} \mu_{j,t}(H + sG) \bigg|_{s=0} 
= \mathrm{vec}\big\{ \dot{A}_k(H)[G] D_{j,t} A_{-k}^\top \big\}.
\end{align}

Let $f_{j,t,H}$ denote the density of $\mathrm{vec}(Y_{j,t})$ under the parameter $H$. We have
\[
\log\big( f_{j,t,H}\{\mathrm{vec}(Y_{j,t})\} \big)
= 
-\frac{1}{2\sigma^2} \|\mathrm{vec}(Y_{j,t}) -\mu_{j,t}(H)\|_2^2 .
\]
Correspondingly, the score of the $t$th record with respect to $\theta=\mathrm{vec}(H)$ is
\begin{align}\label{eq:lower-2.1.3}
    S_{j,t}(H) = \nabla_\theta \log\big( f_{j,t,H}\{\mathrm{vec}(Y_{j,t})\} \big)
    &=
    \frac{1}{\sigma^2} B_{j,t}(H)^\top \{\mathrm{vec}(Y_{j,t}) -\mu_{j,t}(H)\}
    \nonumber\\
    &=
    \frac{1}{\sigma^2} B_{j,t}(H)^\top \mathrm{vec}(E_{j,t}) 
    \in\R^{d_k}.
\end{align}
Since the entries of $E_{j,t}$ are independent $\mathcal{N}(0, \sigma^2)$ random variables, $S_{j,t}(H)$ is a centered Gaussian vector. Moreover,
\begin{align}\label{eq:lower-2.1.4}
    \mathbb{E}_H\{ S_{j,t}(H) S_{j,t}(H)^\top \}
    &= 
    \frac{1}{\sigma^4} B_{j,t}(H)^\top 
    \mathbb{E}_H\{ \mathrm{vec}(E_{j,t}) \mathrm{vec}(E_{j,t})^\top \} B_{j,t}(H) \nonumber\\
    &= 
    \frac{1}{\sigma^2} B_{j,t}(H)^\top B_{j,t}(H) .
\end{align}

\medskip\noindent\textbf{Step 2.2.}
In this sub-step, we verify the conditions in Lemma~\ref{lem:score-contraction}.
For any $g = \mathrm{vec}(G) \in\R^{d_k}$, we have
\begin{align}\label{eq:lower-2.2.1}
    &\;\quad
    g^\top \bigg\{ \sum_{t=1}^{n_j} B_{j,t}(H)^\top B_{j,t}(H) \bigg\} g
    =
    \sum_{t=1}^{n_j} \|B_{j,t}(H)g\|_2^2 
    =
    \sum_{t=1}^{n_j} \|\dot{A}_k(H)[G] D_{j,t} A_{-k}^\top \|_{\mathrm{F}}^2 \nonumber\\
    &=
    \sum_{t=1}^{n_j} \|\dot{A}_k(H)[G]D_{j,t} \|_{\mathrm{F}}^2
    =
    \mathrm{tr}\bigg[ \dot{A}_k(H)[G]^\top \dot{A}_k(H)[G] \bigg\{ \sum_{t=1}^{n_j} D_{j,t} D_{j,t}^\top \bigg\} \bigg] 
    =
    n_j\lambda \|\dot{A}_k(H)[G]\|_{\mathrm{F}}^2 ,
\end{align}
where the first second equality follows from \eqref{eq:lower-2.1.2}, and the last equality follows from \eqref{eq:lower-1.3.2}.

Since $\|g\|_2 = \|G\|_{\mathrm{F}}$, combining \eqref{eq:lower-1.2.1} and \eqref{eq:lower-2.2.1} gives 
\[
(c_2')^2 n_j \lambda \|g\|_2^2 
\leq 
g^\top \bigg\{ \sum_{t=1}^{n_j} B_{j,t}(H)^\top B_{j,t}(H) \bigg\} g 
\leq
(C_2')^2 n_j \lambda\|g\|_2^2 .
\]
Since this holds for every $g\in\R^{d_k}$, we have that 
\begin{equation}\label{eq:lower-2.2.2}
(c_2')^2 n_j\lambda I_{d_k} 
\preceq \sum_{t=1}^{n_j} B_{j,t}(H)^\top B_{j,t}(H) 
\preceq (C_2')^2 n_j\lambda I_{d_k}.
\end{equation}

We also note that 
\begin{align*}
    \|B_{j,t}(H)\|^2 &=
    \sup_{\|g\|_2=1} \| B_{j,t}(H)g \|_2^2 = \sup_{\|G\|_{\mathrm{F}}=1} \| \dot{A}_k(H)[G] D_{j,t} A_{-k}^\top \|_{\mathrm{F}}^2
    = 
    \sup_{\|G\|_{\mathrm{F}}=1} \| \dot{A}_k(H)[G] D_{j,t} \|_{\mathrm{F}}^2 \\
    &\leq 
    \sup_{\|G\|_{\mathrm{F}}=1} \| \dot{A}_k(H)[G] \|_{\mathrm{F}}^2 \|D_{j,t}\|^2 
    \leq
    (C_2')^2 C_r\lambda  ,
\end{align*}
where the second equality follows from  \eqref{eq:lower-2.1.2}, and the last inequality follows from \eqref{eq:lower-1.2.1} and \eqref{eq:lower-1.3.3}. 
Consequently, we obtain
\begin{equation}\label{eq:lower-2.2.3}
\max_{t\in[n_j]} \|B_{j,t}(H)\|^2 
\leq (C_2')^2 C_r\lambda .
\end{equation}

\medskip\noindent\textbf{Step 2.3.}
In this sub-step, we apply Lemma~\ref{lem:score-contraction} to the private message released by client~$j$.

Let $\mathcal{D}_j:= \big(\mathbf{X}_1^{(j)}, \ldots, \mathbf{X}_{n_j}^{(j)} \big)$ be the complete dataset of client~$j$, and let $Z_j:= \mathscr{M}_j(\mathcal{D}_j)$ denote its private message.
We now verify the assumptions of Lemma~\ref{lem:score-contraction}. First,
\[
\bigg\{ \mathrm{vec}(H): \|H\|_{\mathrm{F}}< \frac{1}{2} \bigg\}
\]
is an open subset of $\R^{d_k}$. 
Moreover, $H \mapsto A_k(H)$ is continuously differentiable. Hence $H\mapsto\mu_{j,t}(H)$ is continuously differentiable.
Second, \eqref{eq:lower-2.1.1} shows that the observations are independent Gaussian vectors with covariance $\sigma^2 I_{p_k p_{-k}}$.
Equations \eqref{eq:lower-2.2.2} and \eqref{eq:lower-2.2.3} verify the two Jacobian conditions of Lemma~\ref{lem:score-contraction}, uniformly over $\|H\|_{\mathrm{F}} < 1/2$.
Third, by the definition of FDP, $\mathscr{M}_j$ is $(\varepsilon_j, \delta_j)$-centrally differentially private. 
Finally, the condition on $\delta$ in Lemma~\ref{lem:score-contraction} becomes
\[
d_k \delta_j \log(1/\delta_j) \leq c_0 \sqrt{n_j} \varepsilon_j^2 ,
\]
which follows from the condition in the statement of Theorem~\ref{thm:fed-lower-target}.
Therefore, Lemma~\ref{lem:score-contraction} gives that there exists an absolute constant $C' >0$ such that 
\begin{equation}\label{eq:private-information-client}
\mathrm{tr}\{ I_{Z_j}(H) \} 
\leq
C' \frac{\lambda}{\sigma^2} \big\{ (n_j d_k) \wedge (n_j^2 \varepsilon_j^2) \big\} ,
\end{equation}
uniformly over $\|H\|_{\mathrm{F}} < 1/2$, where $I_{Z_j}(H)$ denotes the Fisher information matrix for $\theta=\mathrm{vec}(H)$ in the induced distribution of $Z_j$.

\medskip\noindent\textbf{Step 2.4.}
In this sub-step, we aggregate the Fisher information bounds over all clients.
Let $\mathbf{Z}:= (Z_1, \ldots, Z_m)$ denote the collection of private messages received by the server. 
Conditional on $H$, the client datasets are independent, 
and $Z_j$ depends only on $\mathcal{D}_j$ and the local randomization of client~$j$. Therefore, $Z_1, \ldots, Z_m$ are independent conditional on $H$.

Let $q_{j,H}$ denote the density of $Z_j$ under the parameter $H$. The joint density of $\mathbf{Z}$ is
\[
q_H(z_1, \ldots, z_m) = \prod_{j=1}^m q_{j,H}(z_j) .
\]
Consequently, the score of $\mathbf{Z}$ with respect to $\theta=\mathrm{vec}(H)$ is
\[
S_{\mathbf{Z}}(H) = \nabla_\theta \log\big( q_H(Z_1, \ldots, Z_m) \big)
= \sum_{j=1}^m \nabla_\theta \log\big( q_{j,H}(Z_j) \big) .
\]
Since the individual message scores are independent and centered under $H$, the Fisher information matrix for $\theta = \mathrm{vec}(H)$ in the induced distribution of $\mathbf{Z}$, denoted by $I_{\mathbf{Z}}(H)$, satisfies
\begin{equation}\label{eq:lower-2.4.1}
I_{\mathbf{Z}}(H) = \sum_{j=1}^m I_{Z_j}(H).
\end{equation}
Taking the trace in \eqref{eq:lower-2.4.1} and applying \eqref{eq:private-information-client}, we obtain
\begin{equation}\label{eq:lower-2.4.2}
\mathrm{tr}\{ I_{\mathbf{Z}}(H) \} 
\leq
C' \frac{\lambda}{\sigma^2} \sum_{j=1}^m \big\{ (n_j d_k) \wedge (n_j^2 \varepsilon_j^2) \big\} .
\end{equation}

\medskip\noindent\textbf{Step 3.} In this final step, we first reduce the projector loss to the squared Euclidean loss for the local coordinate. We then construct a product cosine prior, apply the multivariate van Trees inequality, and express the resulting lower bound in terms of $\underline\Psi_{j,k}^{\mathrm{loc}}$.

\medskip\noindent\textbf{Step 3.1.}
In this sub-step, we reduce estimation of the loading projector to estimation of the local coordinate $H$.

Fix $0< \rho <1/2$ and let $\widehat{A}_k$ be any admissible estimator based on $\mathbf{Z}$. Define $\widehat{H}$ as a measurable selection from
\[
\underset{\|H'\|_{\mathrm{F}} \leq \rho}{\mathrm{argmin}} \| \widehat{A}_k\widehat{A}_k^\top -P_k(H') \|_{\mathrm{F}} .
\]
Such a minimizer exists as the set $\{H': \|H'\|_{\mathrm{F}} \leq \rho\}$ is compact and $P_k(H')$ is continuous in $H'$.

For every $H$ satisfying $\|H\|_{\mathrm{F}} \leq \rho$, both $H$ and $\widehat{H}$ belong to the unit ball. 
By \eqref{eq:lower-1.1.2}, the triangle inequality, and the fact that $\widehat{H}$ is the minimizer, we have
\begin{align*}
    &\;\quad
    c_1' \|\widehat{H} -H\|_{\mathrm{F}}
    \leq 
    \|P_k(\widehat{H}) -P_k(H)\|_{\mathrm{F}} \\
    &\leq 
    \| P_k(\widehat{H}) -\widehat{A}_k\widehat{A}_k^\top \|_{\mathrm{F}} + \| \widehat{A}_k\widehat{A}_k^\top -P_k(H) \|_{\mathrm{F}} 
    \leq 
    2 \| \widehat{A}_k\widehat{A}_k^\top -P_k(H) \|_{\mathrm{F}} .
\end{align*}
Consequently, we have that 
\begin{equation}\label{eq:lower-3.1.1}
\frac{(c_1')^2}{4} \|\widehat{H} -H\|_{\mathrm{F}}^2 
\leq
\| \widehat{A}_k\widehat{A}_k^\top -P_k(H) \|_{\mathrm{F}}^2 .
\end{equation}

\medskip\noindent\textbf{Step 3.2.} In this sub-step, we construct a prior for $\theta=\mathrm{vec}(H)$ and calculate its Fisher information.
Let $a:= \rho /\sqrt{d_k}$ and define the density
\[
\pi_a(u) := \frac{1}{a} \cos^2\bigg( \frac{\pi u}{2a} \bigg) \mathbbm{1}\{|u|\leq a\} .
\]
Since
\[
\frac{1}{a} \int_{-a}^a \cos^2\bigg( \frac{\pi u}{2a} \bigg) du = 1,
\]
it is a probability density. 
Place independent priors with density $\pi_a$ on the $d_k$ coordinates of $\theta=\mathrm{vec}(H)$, and denote the resulting product prior by $\pi_\rho$. 
If $\theta$ belongs to the support of $\pi_\rho$, then
\[
\|H\|_{\mathrm{F}}^2 = \|\theta\|_2^2 
\leq d_k a^2 = \rho^2 .
\]
Thus, $\pi_\rho$ is supported on $\{H: \|H\|_{\mathrm{F}} \leq \rho\}$. For $|u|<a$, the score of the one-dimensional prior is
\[
\frac{d}{du} \log(\pi_a(u)) =
-\frac{\pi}{a} \tan\bigg( \frac{\pi u}{2a} \bigg) .
\]
Hence its Fisher information is
\begin{equation}\label{eq:lower-3.2.1}
\int_{-a}^a \bigg\{ \frac{d}{du} \log(\pi_a(u)) \bigg\}^2 \pi_a(u) du 
=
\frac{\pi^2}{a^3} \int_{-a}^a \sin^2\bigg( \frac{\pi u}{2a} \bigg) du 
= \frac{\pi^2}{a^2} .
\end{equation}
Let $I_{\pi_\rho}$ denote the Fisher information matrix of the product prior. Since its coordinates are independent, by \eqref{eq:lower-3.2.1}, we have 
\begin{equation}\label{eq:lower-3.2.3}
\mathrm{tr}\{I_{\pi_\rho}\} = 
d_k \frac{\pi^2}{a^2} 
= \frac{\pi^2 d_k^2}{\rho^2} .
\end{equation}

\medskip\noindent\textbf{Step 3.3.}
In this sub-step, we apply the classical multivariate van Trees inequality \citep{gill1995applications}, in the form stated in Lemma~D.2 in \citet{cai2026minimax}.

We first verify the required conditions. The joint sub-model obtained by stacking the Gaussian records in \eqref{eq:lower-2.1.1} has a smooth mean and a fixed positive definite covariance matrix. 
Hence, it satisfies the differentiability conditions in Assumption~D.1 of \citet{cai2026minimax}. 
Since $\mathbf{Z}$ is obtained through a parameter-independent Markov kernel, the induced experiment is also differentiable in quadratic mean. 
Moreover, the product cosine prior $\pi_\rho$ is continuously differentiable, vanishes on the boundary of its support, and has finite Fisher information. Thus, the conditions of \citet[Lemma~D.2]{cai2026minimax} are satisfied.

Applying Lemma~D.2 in \citet{cai2026minimax} and using
$\|\mathrm{vec}(\widehat{H}) -\mathrm{vec}(H)\|_2 = \|\widehat{H}-H\|_{\mathrm{F}}$, we obtain
\begin{align}\label{eq:lower-3.3.1}
    \int \mathbb{E}_H \|\widehat{H}-H\|_{\mathrm{F}}^2 d\pi_\rho(H)
    &\geq 
    \frac{d_k^2}{\int \mathrm{tr}\{I_{\mathbf{Z}}(H)\} d\pi_\rho(H) + \mathrm{tr}\{I_{\pi_\rho}\}}  \nonumber\\
    &\geq 
    \frac{d_k^2}{C'\frac{\lambda}{\sigma^2} \sum_{j=1}^m \big\{ (n_j d_k) \wedge (n_j^2 \varepsilon_j^2) \big\}
    + (\pi^2 d_k^2)/\rho^2 } ,
\end{align}
where the last inequality follows from \eqref{eq:lower-2.4.2} and \eqref{eq:lower-3.2.3}. 
Choose
\begin{equation}\label{eq:lower-3.3.2}
\rho^2 = 
c_\rho \bigg( 
\frac{\sigma^2 d_k^2}{\lambda \sum_{j=1}^m \big\{ (n_j d_k) \wedge (n_j^2 \varepsilon_j^2) \big\} }
\wedge1 \bigg),
\end{equation}
where $c_\rho = 1/16$. Since $u^{-1} \leq (u\wedge 1)^{-1}$ for every $u>0$, combining \eqref{eq:lower-3.3.1} and \eqref{eq:lower-3.3.2} gives
\begin{align}\label{eq:lower-3.3.3}
    \int \mathbb{E}_H \|\widehat{H}-H\|_{\mathrm{F}}^2 d\pi_\rho(H)
    &\geq
    C_2 \bigg( \frac{\sigma^2 d_k^2}{\lambda \sum_{j=1}^m \big\{ (n_j d_k) \wedge (n_j^2 \varepsilon_j^2) \big\} }
    \wedge 1 \bigg) \nonumber\\
    &\geq
    \frac{C_2}{2} \Bigg\{ \frac{1}{\sum_{j=1}^m \Big( \frac{\sigma^2 d_k}{\lambda n_j} + \frac{\sigma^2 d_k^2}{\lambda n_j^2 \varepsilon_j^2} \Big)^{-1}}
    \wedge 2 \Bigg\} ,
\end{align}
where $C_2= 1/(C' + 16\pi^2) >0$ is an absolute constant, and the last inequality follows from the fact that for any $x, y > 0$,
\[
x\wedge y \leq 2(x^{-1} + y^{-1})^{-1} .
\]
Moreover, since $p_k \geq 2r_k$,
\[
d_k = r_k(p_k-r_k)  \geq \frac{r_k p_k}{2}
\]
It follows from the definition of $\underline\Psi_{j,k}^{\mathrm{loc}}$ that
\begin{equation}\label{eq:lower-3.3.4}
\frac{\sigma^2 d_k}{\lambda n_j} + \frac{\sigma^2 d_k^2}{\lambda n_j^2 \varepsilon_j^2} 
\geq  
\frac{r_k (\underline\Psi_{j,k}^{\mathrm{loc}})^2}{4} .
\end{equation}
Combining \eqref{eq:lower-3.1.1}, \eqref{eq:lower-3.3.3}, and \eqref{eq:lower-3.3.4}, we obtain
\begin{align}\label{eq:lower-3.3.5}
    \int \mathbb{E}_H \| \widehat{A}_k\widehat{A}_k^\top -P_k(H) \|_{\mathrm{F}}^2
    d\pi_\rho(H)  
    &\geq
    \frac{C_2(c_1')^2}{32} \Bigg\{ 
    \frac{r_k}{\sum_{j=1}^m (\underline\Psi_{j,k}^{\mathrm{loc}})^{-2}} 
    \wedge 8 \Bigg\} \nonumber\\
    &\geq
    c_2 \Bigg\{
    \frac{r_k}{\sum_{j=1}^m (\underline\Psi_{j,k}^{\mathrm{loc}})^{-2}} 
    \wedge (2 r_k) \Bigg\} ,
\end{align}
where the last inequality follows since under the fixed-rank assumption, $r_k\leq r_{\max}$ for some fixed absolute constant
$r_{\max}$ and, for every $x\geq 0$,
\[
x\wedge 8 \geq \min\bigg\{1, \frac{4}{r_{\max}} \bigg\} (x\wedge (2 r_k)).
\]
with 
\[
c_2 = \frac{C_2(c_1')^2 \min\{1, 4/r_{\max}\}}{32},
\]
The support of $\pi_\rho$ is contained in the Gaussian subclass constructed in \textbf{Step 1}. Hence for every admissible estimator,
\begin{align*}
    \sup_{\mathscr{P}\in \mathfrak{P}_k^{\mathrm{fed}}(\lambda, \sigma^2)}
    \mathbb{E}_{\mathscr{P}} \| \widehat{A}_k\widehat{A}_k^\top -A_kA_k^\top \|_{\mathrm{F}}^2 
    &\geq 
    \int \mathbb{E}_H \| \widehat{A}_k\widehat{A}_k^\top -P_k(H) \|_{\mathrm{F}}^2 d\pi_\rho(H) \\
    &\geq
    c_2 \Bigg\{ 
    \frac{r_k}{\sum_{j=1}^m (\underline\Psi_{j,k}^{\mathrm{loc}})^{-2}} 
    \wedge (2 r_k) \Bigg\},
\end{align*}
which completes the proof. 
\end{proof}

\section{Proofs for adaptive tuning}
\label{sec:adaptive-proofs}

\subsection{Auxiliary lemmas with proofs}

For the adaptive calibration arguments, define
\begin{equation}
\label{eq:adaptive-intrinsic-radius}
\varrho_{j,k} :=
\sqrt{r_k}
\left[
\frac{\sigma}{\sqrt{\lambda_{j,1}^{(k)}}}
\sqrt{\frac{p_k}{b_j}}
+
\frac{\sigma^2}{\lambda_{j,1}^{(k)}}
\left\{
\sqrt{\frac{p}{b_j}} + \frac{p_k}{b_j}
\right\}
\right].
\end{equation}
Let $C_p > 0$ be the fixed dimension constant used in Theorem~\ref{thm:adaptive-upper}, and define $C_{\mathrm{cmp}} := 1 + C_p^{-1}$.
Under $p_k \geq C_p \log(n \vee p)$, the radius in \eqref{eq:block-radius} satisfies
\begin{equation}
\label{eq:adaptive-radius-comparison}
C_R \varrho_{j,k} \leq  R_{j,k}
\leq 
C_{\mathrm{cmp}} C_R \varrho_{j,k}.
\end{equation}

\begin{lemma}
\label{lem:adaptive-model-calibration}
Fix $j\in[m]$, and let $\mathcal{P}_j^{(0)}$ be the sigma-algebra generated by its pilot partition and deterministic factors. Define
\begin{equation}
\label{eq:adaptive-pair-cdf}
G_{j,k}(d) :=
\frac{1}{(B_j^{(0)})^2}
\sum_{\ell=1}^{B_j^{(0)}}
\sum_{h=1}^{B_j^{(0)}}
\P\left\{
\big\| \widehat\Pi_{j,k}^{(0),[\ell]} -\widehat\Pi_{j,k}^{(0),[h]} \big\|_{\mathrm{F}}
\leq d
\mathrel{\Big|}
\mathcal{P}_j^{(0)}
\right\}.
\end{equation}
Suppose that Assumptions~\ref{ass:loading}--\ref{ass:noise} hold.
Let
\[
b_j =
\left\lceil C_b \{ r_k + \log(n \vee p) \} \right\rceil,
\qquad
p_k \geq 2 r_k,
\qquad
p_k \geq C_p \log(n \vee p),
\qquad
\frac{\lambda_{j,1}^{(k)}}{\sigma^2}
\geq C_0 ( p_k + \sqrt{p} ).
\]
Also assume $B_j^{(0)} \geq 8$ and \eqref{eq:block-balance} holds for every pilot block.
Then for sufficiently large absolute constants $C_b$, $C_p$, and $C_0$,
there exist absolute constants $0 < c_- < c_+ < \infty$
such that
\begin{equation}
\label{eq:adaptive-model-calibration}
G_{j,k}(c_- \varrho_{j,k}) \leq \frac{3}{4},
\quad
G_{j,k}(c_+ \varrho_{j,k}) \geq \frac{15}{16} .
\end{equation}
\end{lemma}

\begin{proof}[Proof of Lemma~\ref{lem:adaptive-model-calibration}]
Condition on $\mathcal{P}_j^{(0)}$, and we suppress the indices $j$ and $k$ throughout this proof. 
Write $P:=A_kA_k^\top$ and $q_k:=p_k-r_k$
and choose $U_k\in \mathbb{O}_{p_k\times q_k}$ such that $U_kU_k^\top = I_{p_k} -P$. 
Within this proof, we write $\mathcal{I}_\ell := \mathcal{I}_{j,\ell}^{(0)}$ for the $\ell$th pilot block, suppressing the client index and the pilot superscript.
Then for a pilot block $\mathcal{I}_\ell$, set
\[
s_\ell = |\mathcal{I}_\ell| ,
\quad
S_\ell = \frac{1}{s_\ell}
\sum_{t\in\mathcal{I}_\ell} D_tD_t^\top ,
\]
where $D_t:= \mathcal{M}_k(\mathbf{F}_t^{(j)})$. 
By \eqref{eq:block-balance}, we have
\begin{equation}
\label{eq:adaptive-factor-scale}
\frac{1}{2} \Lambda_{j,(k)}
\preceq 
S_\ell
\preceq
\frac{3}{2} \Lambda_{j,(k)},
\quad
b_j\leq s_\ell\leq 2 b_j.
\end{equation}

\medskip\noindent\textbf{Step 1.}
In this step, we first establish an upper bound for the distance between two pilot block projectors.
Apply \eqref{eq:block-projector-fixed-bound} with a fixed sufficiently large tail parameter. Assumption~\ref{ass:latent}, \eqref{eq:adaptive-factor-scale}, and $p_k\gtrsim \log(n\vee p)$ give
\begin{equation}
\label{eq:adaptive-one-block-upper}
\P\left\{
\|\widehat\Pi_{j,k}^{(0),[\ell]}-P\|_{\mathrm{F}}
> C \varrho_{j,k}
\mathrel{\Big|}
\mathcal{P}_j^{(0)}
\right\}
\leq 
\frac{1}{64} .
\end{equation}
For two distinct blocks, the triangle inequality and \eqref{eq:adaptive-one-block-upper} therefore give, after choosing
$c_+$ sufficiently large,
\begin{equation}
\label{eq:adaptive-pair-upper}
\P\left\{
\big\| \widehat\Pi_{j,k}^{(0),[\ell]} - \widehat\Pi_{j,k}^{(0),[h]}
\big\|_{\mathrm{F}}
\leq 
c_+ \varrho_{j,k}
\mathrel{\Big|}
\mathcal{P}_j^{(0)}
\right\}
\geq
\frac{31}{32}.
\end{equation}

\medskip\noindent\textbf{Step 2.}
Next, we derive a first-order expansion of each pilot-block spectral
projector around the true loading-space projector.
For a pilot block $\mathcal{I}_\ell$, let
\[
\widehat\Sigma_\ell :=
\frac{1}{s_\ell} \sum_{t\in \mathcal{I}_\ell} \mathcal{M}_k(\mathbf{X}_t^{(j)})
\mathcal{M}_k(\mathbf{X}_t^{(j)})^\top .
\]
Conditional on $\mathcal{P}_j^{(0)}$, we have
\[
\Sigma_\ell^\star := \mathbb{E}[\widehat\Sigma_\ell]
= A_k S_\ell A_k^\top + \sigma^2 p_{-k} I_{p_k} .
\]
Next, set $\Delta_\ell = \widehat\Sigma_\ell-\Sigma_\ell^\star$
and $P= A_kA_k^\top$.
Since $A_k^\top A_k=I_{r_k}$ and
$U_kU_k^\top = I_{p_k} -P$, we have
\[
\begin{pmatrix}
A_k & U_k
\end{pmatrix}^\top
\Sigma_\ell^\star
\begin{pmatrix}
A_k & U_k
\end{pmatrix}
=
\begin{pmatrix}
S_\ell + \sigma^2 p_{-k} I_{r_k} & 0 \\
0 & \sigma^2 p_{-k} I_{q_k}
\end{pmatrix} .
\]
Thus $P$ is the leading rank-$r_k$ spectral projector of $\Sigma_\ell^\star$, with eigengap $\lambda_{r_k}(S_\ell)$.

For $u\in [0,1]$, let $\Pi_\ell(u)$ be the leading rank-$r_k$ spectral projector of $\Sigma_\ell^\star + u\Delta_\ell$. 
For the remainder of the proof, we work on the event
\begin{equation}\label{eq-lemma-adap-eigen-event}
\mathcal{E}_{\mathrm{gap}} := 
\bigg\{ \|\Delta_\ell\| \leq \frac{1}{4} \lambda_{r_k}(S_\ell) \bigg\} .
\end{equation}
Then Weyl's inequality gives
\begin{align}\label{eq-lemma-adap-eigengap}
    \lambda_{r_k}(\Sigma_\ell^\star + u\Delta_\ell)
    - \lambda_{r_k+1}(\Sigma_\ell^\star + u\Delta_\ell)
    &\geq
    \lambda_{r_k}(S_\ell) - 2u \|\Delta_\ell\|
    \geq
    \frac{1}{2} \lambda_{r_k}(S_\ell) ,
    \quad u\in[0,1].
\end{align}
The uniform eigengap ensures that $\Pi_\ell(u)$ is twice differentiable along the path, and
\[
\Pi_\ell(0) = P,
\quad
\Pi_\ell(1) = \widehat\Pi_{j,k}^{(0),[\ell]}.
\]

We first calculate $\dot\Pi_\ell(0)$. Differentiating
\begin{equation}
\label{eq-lemma-adap-0}
(\Sigma_\ell^\star + u\Delta_\ell) \Pi_\ell(u)
=
\Pi_\ell(u)(\Sigma_\ell^\star + u\Delta_\ell),
\quad
\Pi_\ell(u)^2 = \Pi_\ell(u),
\end{equation}
at $u=0$ gives
\begin{equation}\label{eq-lemma-adap-1}
\Delta_\ell P + \Sigma_\ell^\star \dot\Pi_\ell(0)
=
\dot\Pi_\ell(0) \Sigma_\ell^\star + P\Delta_\ell ,
\quad
\dot\Pi_\ell(0)P + P \dot\Pi_\ell(0)
=
\dot\Pi_\ell(0) .
\end{equation}
The second identity implies
\begin{equation}\label{eq-lemma-adap-2}
P\dot\Pi_\ell(0)P =0 ,
\quad
(I_{p_k} -P) \dot\Pi_\ell(0) (I_{p_k}-P) = 0.
\end{equation}
To determine the off-diagonal blocks, multiply the first differentiated identity in \eqref{eq-lemma-adap-1} on the left by $I_{p_k} -P$ and on the right by $P$. The term containing $P\Delta_\ell$
vanishes because $(I_{p_k} -P) P = 0$.
Hence we have
\begin{align*}
    & (I_{p_k} -P) \Delta_\ell P
    + (I_{p_k} -P) \Sigma_\ell^\star \dot\Pi_\ell(0)P =
    (I_{p_k} -P) \dot\Pi_\ell(0) \Sigma_\ell^\star P .
\end{align*}
By $\Sigma_\ell^\star = A_k S_\ell A_k^\top + \sigma^2 p_{-k} I_{p_k}$,
$PA_k=A_k$, 
and $(I_{p_k} -P) A_k =0$, we have
\[
(I_{p_k} -P) \Sigma_\ell^\star
=
\sigma^2 p_{-k} (I_{p_k} -P) ,
\quad
\Sigma_\ell^\star P
=
A_k S_\ell A_k^\top + \sigma^2 p_{-k} P .
\]
Substituting these expressions into the preceding identity gives
\begin{align*}
    & (I_{p_k} -P) \Delta_\ell P
    + \sigma^2 p_{-k} (I_{p_k} -P) \dot\Pi_\ell(0) P =
    (I_{p_k} -P) \dot\Pi_\ell(0) A_k S_\ell A_k^\top
    + \sigma^2 p_{-k} (I_{p_k} -P) \dot\Pi_\ell(0) P .
\end{align*}
Hence
\[
(I_{p_k} -P) \Delta_\ell P
=
(I_{p_k} -P) \dot\Pi_\ell(0) A_k S_\ell A_k^\top .
\]
Right-multiplying both sides by $A_k S_\ell^{-1} A_k^\top$ and using $PA_k=A_k$ and $A_k^\top A_k=I_{r_k}$ gives
\[
(I_{p_k} -P) \Delta_\ell A_k S_\ell^{-1} A_k^\top
=
(I_{p_k} -P) \dot\Pi_\ell(0) P .
\]
Each $\Pi_\ell(u)$ is a symmetric projector, so $\dot\Pi_\ell(0)$ is symmetric. 
As $\Delta_\ell$ is also symmetric, transposing the last equality gives
\[
P \dot\Pi_\ell(0) (I_{p_k} -P)
=
A_k S_\ell^{-1} A_k^\top \Delta_\ell (I_{p_k} -P) .
\]
Finally, the two diagonal blocks of $\dot\Pi_\ell(0)$ were shown above to vanish. Therefore
\begin{align*}
    \dot\Pi_\ell(0) &=
    (I_{p_k} -P) \dot\Pi_\ell(0) P
    + P \dot\Pi_\ell(0)(I_{p_k} -P) \\
    &=
    (I_{p_k} -P) \Delta_\ell A_k S_\ell^{-1} A_k^\top
    + A_k S_\ell^{-1} A_k^\top \Delta_\ell (I_{p_k} -P) .
\end{align*}
Using $I_{p_k}-P = U_kU_k^\top$, define
\begin{align}\label{eq-lemma-adap-def-L}
    & K_\ell :=
    U_k^\top \Delta_\ell A_k S_\ell^{-1},
    \quad
    \mathcal{L}_\ell :=
    \dot\Pi_\ell(0)
    = U_k K_\ell A_k^\top
    + A_k K_\ell^\top U_k^\top , 
    \quad
    \mathcal{R}_\ell :=
    \widehat\Pi_{j,k}^{(0),[\ell]} -P -\mathcal{L}_\ell .
\end{align}
Then we can write
\begin{equation}
\label{eq:adaptive-projector-expansion}
\widehat\Pi_{j,k}^{(0),[\ell]} -P =
\mathcal{L}_\ell + \mathcal{R}_\ell.
\end{equation}

\medskip\noindent\textbf{Step 3.}
We establish uniform operator- and Frobenius-norm bounds for $\dot\Pi_\ell(u)$ along the perturbation path. These bounds will be used to control the remainder $\mathcal{R}_\ell$.

Fix $u\in[0,1]$.
Differentiating \eqref{eq-lemma-adap-0} gives
\begin{align}
    & \Pi_\ell(u) \dot\Pi_\ell(u)
    + \dot\Pi_\ell(u) \Pi_\ell(u)
    =
    \dot\Pi_\ell(u) , 
    \label{eq-lemma-adap-3} \\
    & \Delta_\ell\Pi_\ell(u)
    + (\Sigma_\ell^\star + u\Delta_\ell) \dot\Pi_\ell(u)
    =
    \dot\Pi_\ell(u)(\Sigma_\ell^\star + u\Delta_\ell)
    + \Pi_\ell(u) \Delta_\ell .
    \label{eq-lemma-adap-3.1}
\end{align}
Multiplying \eqref{eq-lemma-adap-3} on both sides by $\Pi_\ell(u)$ gives
\begin{equation}\label{eq-lemma-adap-4}
\Pi_\ell(u) \dot\Pi_\ell(u) \Pi_\ell(u) =0.
\end{equation}
Similarly, multiplying \eqref{eq-lemma-adap-3} on both sides by $I_{p_k}-\Pi_\ell(u)$ gives
\begin{equation}\label{eq-lemma-adap-5}
(I_{p_k}-\Pi_\ell(u)) \dot\Pi_\ell(u)
(I_{p_k}-\Pi_\ell(u)) = 0.
\end{equation}
Thus the two diagonal blocks of $\dot\Pi_\ell(u)$, relative to the range of $\Pi_\ell(u)$ and its orthogonal complement, are zero.
Since $\dot\Pi_\ell(u)$ is symmetric, its two remaining off-diagonal blocks are transposes of each other.

Fix $u\in[0,1]$. After taking the derivatives above, choose an orthonormal eigen-basis of $\Sigma_\ell^\star + u\Delta_\ell$, with its leading $r_k$ eigenvectors first. 
All matrix entries below are read in this basis, which is held fixed.
In this basis, $\Sigma_\ell^\star + u\Delta_\ell$ is diagonal. 
Also, $\Pi_\ell(u)$ has diagonal entries one in the first $r_k$ positions and zero in the remaining positions. 
Equations \eqref{eq-lemma-adap-4} and \eqref{eq-lemma-adap-5} show that the two diagonal blocks of $\dot\Pi_\ell(u)$ are zero.
As $\Pi_\ell(u)$ is symmetric for every $u$, $\dot\Pi_\ell(u)$ is symmetric as well.

Now take row $a >r_k$ and column $i\leq r_k$. Note that
\begin{align*}
    &
    [\Delta_\ell\Pi_\ell(u)]_{ai}
    =
    [\Delta_\ell]_{ai} 
    , \quad
    [(\Sigma_\ell^\star + u\Delta_\ell) \dot\Pi_\ell(u)]_{ai}
    =
    \lambda_a(\Sigma_\ell^\star + u\Delta_\ell) [\dot\Pi_\ell(u)]_{ai}, \\
    &
    [\dot\Pi_\ell(u) (\Sigma_\ell^\star + u\Delta_\ell)]_{ai}
    =
    [\dot\Pi_\ell(u)]_{ai} \lambda_i(\Sigma_\ell^\star + u\Delta_\ell)   
    , \quad
    [\Pi_\ell(u)\Delta_\ell]_{ai}
    = 0 ,
\end{align*}
where the first equality holds because, in the chosen eigen-basis, right multiplication by $\Pi_\ell(u)$ leaves the $i$th column unchanged, the second and third equalities hold because left and right multiplication by the diagonal matrix $\Sigma_\ell^\star + u\Delta_\ell$ scale the $a$th row and the $i$th column, respectively, by their corresponding eigenvalues, and the last equality holds because the $a$th row of $\Pi_\ell(u)$ is zero.
Substituting these four entries into \eqref{eq-lemma-adap-3.1} gives
\[
[\Delta_\ell]_{ai}
+ \lambda_a(\Sigma_\ell^\star + u\Delta_\ell) 
[\dot\Pi_\ell(u)]_{ai}
=
[\dot\Pi_\ell(u)]_{ai}
\lambda_i(\Sigma_\ell^\star + u\Delta_\ell).
\]
Hence we have
\[
[\dot\Pi_\ell(u)]_{ai}
=
\frac{[\Delta_\ell]_{ai}}
{\lambda_i(\Sigma_\ell^\star + u\Delta_\ell)
- \lambda_a(\Sigma_\ell^\star + u\Delta_\ell)}.
\]
As $i\leq r_k <a$, \eqref{eq-lemma-adap-eigengap} gives
\begin{equation}\label{eq-lemma-adap-6}
|[\dot\Pi_\ell(u)]_{ai}|
\leq
\frac{2 |[\Delta_\ell]_{ai}|}
{\lambda_{r_k}(S_\ell)} .
\end{equation}

Hence we have
\begin{align}
    \|\dot\Pi_\ell(u)\|_{\mathrm{F}}^2
    &=
    2 \sum_{i\leq r_k <a}
    |[\dot\Pi_\ell(u)]_{ai}|^2 
    \leq
    \frac{8}{\lambda_{r_k}(S_\ell)^2}
    \sum_{i\leq r_k<a}|[\Delta_\ell]_{ai}|^2 \nonumber \\
    &=
    \frac{8}{\lambda_{r_k}(S_\ell)^2}
    \| (I_{p_k}-\Pi_\ell(u))
    \Delta_\ell \Pi_\ell(u) \|_{\mathrm{F}}^2
    \leq
    \frac{8}{\lambda_{r_k}(S_\ell)^2}
    \|\Delta_\ell\|^2 \|\Pi_\ell(u)\|_{\mathrm{F}}^2 
    =
    \frac{8 r_k \|\Delta_\ell\|^2}{\lambda_{r_k}(S_\ell)^2} ,
    \label{eq-lemma-adap-7}
\end{align}
where the first inequality used \eqref{eq-lemma-adap-6},
the second equality follows from the fact that $(I_{p_k}-\Pi_\ell(u)) \Delta_\ell \Pi_\ell(u)$ selects exactly the entries with $a >r_k$ and $i\leq r_k$, 
and the last equality follows from  $\|\Pi_\ell(u)\|_{\mathrm{F}}^2 = \mathrm{rank}\{\Pi_\ell(u)\} =r_k$.

For the operator norm, a similar argument yields
\begin{equation}\label{eq-lemma-adap-8}
\|
\dot\Pi_\ell(u) \| \leq
\frac{2 \|\Delta_\ell\|}
{\lambda_{r_k}(S_\ell)} .
\end{equation}
Combining \eqref{eq-lemma-adap-7} and \eqref{eq-lemma-adap-8}, and using the invariance of both norms under orthogonal changes of basis, on the event $\mathcal{E}_{\mathrm{gap}}$, we obtain the following bounds in the original basis, uniformly over $u\in[0,1]$:
\begin{equation}\label{eq-lemma-adap-9}
\|\dot\Pi_\ell(u)\|
\leq
\frac{2 \|\Delta_\ell\|}
{\lambda_{r_k}(S_\ell)},
\quad
\|\dot\Pi_\ell(u)\|_{\mathrm{F}}
\leq
\frac{2 \sqrt{2 r_k} \|\Delta_\ell\|}
{\lambda_{r_k}(S_\ell)} .
\end{equation}

\medskip\noindent\textbf{Step 4.}
In this step, we show that, with high conditional probability, the second-order
projector remainder $\mathcal R_\ell$ is a small fraction of the intrinsic radius $\varrho_{j,k}$.

Differentiating \eqref{eq-lemma-adap-3} and \eqref{eq-lemma-adap-3.1} once more gives
\begin{align*}
    & \Pi_\ell(u) \ddot\Pi_\ell(u)
    + \ddot\Pi_\ell(u) \Pi_\ell(u)
    + 2 \dot\Pi_\ell(u)^2
    = 
    \ddot\Pi_\ell(u) , 
    \\
    & (\Sigma_\ell^\star + u\Delta_\ell) \ddot\Pi_\ell(u) 
    - \ddot\Pi_\ell(u)(\Sigma_\ell^\star + u\Delta_\ell)
    =
    2 \{ \dot\Pi_\ell(u) \Delta_\ell
    - \Delta_\ell \dot\Pi_\ell(u) \} .
\end{align*}
Applying the same blockwise argument as \textbf{Step 3}, with the first identity controlling the diagonal blocks and the second controlling the off-diagonal blocks through the eigengap, together with \eqref{eq-lemma-adap-9} yields that
\[
\|\ddot\Pi_\ell(u)\|_{\mathrm{F}}
\leq
C_2 \sqrt{r_k}
\frac{\|\Delta_\ell\|^2}
{\lambda_{r_k}(S_\ell)^2} ,
\quad u\in[0,1] ,
\]
where $C_2 >0$ is an absolute constant.
Since $\Pi_\ell(0)=P$,
$\Pi_\ell(1) = \widehat\Pi_{j,k}^{(0),[\ell]}$,
and $\dot\Pi_\ell(0) = \mathcal{L}_\ell$, Taylor's formula yields
\[
\mathcal{R}_\ell =
\int_0^1 (1-u) \ddot\Pi_\ell(u) du .
\]
Hence, on the event $\mathcal{E}_{\mathrm{gap}}$, it holds that
\begin{equation}
\label{eq:adaptive-projector-remainder}
\|\mathcal{R}_\ell\|_{\mathrm{F}}
\leq
C \sqrt{r_k}
\frac{\|\Delta_\ell\|^2}
{\lambda_{r_k}(S_\ell)^2} .
\end{equation}

It remains to verify the event  $\mathcal{E}_{\mathrm{gap}}$ defined in \eqref{eq-lemma-adap-eigen-event}. Define
\[
\zeta_{j,k} :=
\frac{\sigma \sqrt{\lambda_{j,1}^{(k)}}}{\lambda_{j,r_k}^{(k)}}
\sqrt{\frac{p_k}{b_j}}
+
\frac{\sigma^2}{\lambda_{j,r_k}^{(k)}}
\left\{ \sqrt{\frac{p}{b_j}} + \frac{p_k}{b_j} \right\}.
\]
By Assumption~\ref{ass:latent}
and the condition stated in Lemma~\ref{lem:adaptive-model-calibration}, 
we have
\[
\zeta_{j,k} \leq C_3 \left( C_0^{-1/2} + C_0^{-1} \right) ,
\]
where $C_3 >0$ is an absolute constant.
By \eqref{eq:block-total-perturbation-bound}, there exists an absolute constant $C_4 >0$ such that 
\begin{equation}\label{eq-lemma-adap-10}
\P\left\{
\frac{\|\Delta_\ell\|}{\lambda_{r_k}(S_\ell)}
> C_4 \zeta_{j,k}
\mathrel{\Big|} \mathcal{P}_j^{(0)}
\right\}
\leq
\frac{1}{96} .
\end{equation}
We choose $C_0$ sufficiently large that $C_4 \zeta_{j,k}\leq 1/4$. 
On the complement of the exceptional event in \eqref{eq-lemma-adap-10},
\[
\frac{\|\Delta_\ell\|}{\lambda_{r_k}(S_\ell)}
\leq 
C_4 \zeta_{j,k}
\leq \frac{1}{4} .
\]
Hence, conditional on $\mathcal{P}_j^{(0)}$, the event $\mathcal{E}_{\mathrm{gap}}$ defined in \eqref{eq-lemma-adap-eigen-event} holds with probability at least $95/96$. 
On this event, \eqref{eq:adaptive-projector-remainder} gives
\begin{equation}\label{eq-lemma-adap-11}
\|\mathcal{R}_\ell\|_{\mathrm{F}}
\leq
C_2 \sqrt{r_k}
\left\{
\frac{\|\Delta_\ell\|}{\lambda_{r_k}(S_\ell)}
\right\}^2
\leq 
C_5 \sqrt{r_k} \zeta_{j,k}^2 ,
\end{equation}
where $C_5 >0$ is an absolute constant.
Note that 
\begin{equation}\label{eq-lemma-adap-12}
\sqrt{r_k} \zeta_{j,k}
=
\frac{\lambda_{j,1}^{(k)}}{\lambda_{j,r_k}^{(k)}} \varrho_{j,k} .
\end{equation}
Combining Assumption~\ref{ass:latent}, \eqref{eq-lemma-adap-11}, and \eqref{eq-lemma-adap-12}, on the event $\mathcal{E}_{\mathrm{gap}}$, we have
\[
\|\mathcal{R}_\ell\|_{\mathrm{F}} \leq C_6 \zeta_{j,k}\varrho_{j,k} ,
\]
where $C_6 >0$ is an absolute constant.

Fix a sufficiently small absolute constant $c_1>0$.
Since $\zeta_{j,k}$ can be made uniformly small by increasing $C_0$, choose $C_0$ further
so that $C_6 \zeta_{j,k}\leq c_1$.  
Hence we have
\begin{equation}\label{eq:adaptive-remainder-control}
\P\left\{
\|\mathcal{R}_\ell\|_{\mathrm{F}}
> c_1\varrho_{j,k}
\mathrel{\Big|} \mathcal{P}_j^{(0)}
\right\}
\leq
\frac{1}{96} .
\end{equation}

\medskip\noindent\textbf{Step 5.}
In this step, we establish a lower bound of order $\varrho_{j,k}$ for the distance between projectors from any two distinct pilot blocks, with probability at least $2/3$ conditional on $\mathcal{P}_j^{(0)}$.

Fix two distinct pilot blocks $\mathcal{I}_\ell$ and $\mathcal{I}_h$. By the definition of $\mathcal{L}_\ell$ in \eqref{eq-lemma-adap-def-L}, we have
\[
\mathcal{L}_\ell -\mathcal{L}_h
=
U_k (K_\ell-K_h) A_k^\top
+ A_k (K_\ell-K_h)^\top U_k^\top.
\]
Hence
\begin{align*}
    \|\mathcal{L}_\ell -\mathcal{L}_h\|_{\mathrm{F}}^2
    &=
    \left\| U_k (K_\ell-K_h) A_k^\top
    + A_k (K_\ell-K_h)^\top U_k^\top
    \right\|_{\mathrm{F}}^2
    =
    2 \|K_\ell-K_h\|_{\mathrm{F}}^2 .
\end{align*}
Moreover, \eqref{eq:adaptive-projector-expansion} gives
\[
\widehat\Pi_{j,k}^{(0),[\ell]}
- \widehat\Pi_{j,k}^{(0),[h]}
=
(\mathcal{L}_\ell - \mathcal{L}_h)
+ (\mathcal{R}_\ell - \mathcal{R}_h),
\]
and therefore
\begin{equation}\label{eq-lemma-adap-13}
\left\| \widehat\Pi_{j,k}^{(0),[\ell]}
- \widehat\Pi_{j,k}^{(0),[h]} \right\|_{\mathrm{F}}
\geq
\sqrt{2} \|K_\ell-K_h\|_{\mathrm{F}}
- \|\mathcal{R}_\ell\|_{\mathrm{F}} 
- \|\mathcal{R}_h\|_{\mathrm{F}} .
\end{equation}

The remainders $\|\mathcal{R}_\ell\|_{\mathrm{F}}$ and $ \|\mathcal{R}_h\|_{\mathrm{F}}$ are controlled by \eqref{eq:adaptive-remainder-control}.
Hence it remains to argue that $\|K_\ell-K_h\|_{\mathrm{F}}$ is unlikely to be much smaller than $\varrho_{j,k}$.

To this end,
we derive a representation of $K_\ell$, and the same argument applies to $K_h$.
Recall that
\[
\mathcal{M}_k(\mathbf{X}_t^{(j)})
= A_k D_t A_{-k}^\top + E_t ,
\quad
A_{-k}
= \bigotimes_{u\in [K] \setminus\{k\}} A_u.
\]
Using $s_\ell^{-1} \sum_{t\in \mathcal{I}_\ell} D_tD_t^\top = S_\ell$,
we can write the covariance perturbation as
\[
\Delta_\ell =
\frac{1}{s_\ell} \sum_{t\in \mathcal{I}_\ell}
\left\{
A_k D_t A_{-k}^\top E_t^\top
+ E_t A_{-k} D_t^\top A_k^\top
+ E_t E_t^\top -\sigma^2 p_{-k} I_{p_k}
\right\}.
\]
Multiplying on the left by $U_k^\top$ and on the right by $A_k$, the first term vanishes as $U_k^\top A_k=0$.
The constant part of the last term vanishes for the same reason. The remaining terms give
\[
U_k^\top \Delta_\ell A_k
=
\frac{1}{s_\ell} \sum_{t\in \mathcal{I}_\ell}
U_k^\top E_t
\left( A_{-k} D_t^\top + E_t^\top A_k \right).
\]
Define
\[
Z_t := 
U_k^\top E_t
\in\mathbb{R}^{q_k\times p_{-k}},
\quad
H_t := A_{-k} D_t^\top + E_t^\top A_k
\in\mathbb{R}^{p_{-k}\times r_k} .
\]
As $K_\ell= U_k^\top \Delta_\ell A_k S_\ell^{-1}$, we have
\begin{equation}
\label{eq:adaptive-tangent-representation}
K_\ell =
\frac{1}{s_\ell} \sum_{t\in \mathcal{I}_\ell} Z_t H_t S_\ell^{-1}.
\end{equation}

We next obtain the conditional distribution of $K_\ell-K_h$.
Conditional on $\mathcal{P}_j^{(0)}$, the $D_t$ are fixed.
The Gaussian projections $U_k^\top E_t$ and $A_k^\top E_t$ have cross-covariance
\[
\mathrm{Cov}\left\{
\mathrm{vec}(U_k^\top E_t),
\mathrm{vec}(A_k^\top E_t)
\right\}
=
\sigma^2 \left(I_{p_{-k}} \otimes U_k^\top A_k\right)
= 0 .
\]
They are therefore independent. In particular, $Z_t$ is independent of $H_t$. Conditional on $\mathcal{P}_j^{(0)}$ and the $H_t$ in block $\mathcal{I}_\ell$, each row of \eqref{eq:adaptive-tangent-representation} is a centered Gaussian vector with covariance
\begin{equation}
\label{eq:adaptive-tangent-Gaussian}
V_\ell :=
\frac{\sigma^2}{s_\ell^2} S_\ell^{-1}
\left( \sum_{t\in \mathcal{I}_\ell} H_t^\top H_t \right)
S_\ell^{-1} .
\end{equation}
The rows are independent. The same statement holds for $K_h$, with row covariance $V_h$.

Since the two blocks are disjoint, $K_\ell$ and $K_h$ are conditionally independent given $\mathcal{P}_j^{(0)}$ and the $H_t$ in both blocks. 
Thus $K_\ell-K_h$ has independent centered Gaussian rows, each with covariance $V_\ell+V_h$. Since it has $q_k$ rows,
\begin{equation}
\label{eq:adaptive-tangent-second-moment}
v_{\ell h}^2 :=
\E\left\{
\|K_\ell-K_h\|_{\mathrm{F}}^2
\mathrel{\Big|}
\mathcal{P}_j^{(0)} ,
\{H_t: t\in \mathcal{I}_\ell \cup \mathcal{I}_h\}
\right\}
=
q_k \mathrm{tr}(V_\ell + V_h) .
\end{equation}
We now show that this conditional second moment is 
of order $\varrho_{j,k}^2$.
For one block, define
\[
Q_\ell :=
\sum_{t\in \mathcal{I}_\ell} \| H_t S_\ell^{-1} \|_{\mathrm{F}}^2 .
\]
Taking the trace in \eqref{eq:adaptive-tangent-Gaussian} gives
\begin{equation}
\label{eq:adaptive-V-trace}
\mathrm{tr}(V_\ell) =
\frac{\sigma^2}{s_\ell^2} Q_\ell .
\end{equation}
Thus it suffices to determine the size of $Q_\ell$.

Conditional on $\mathcal{P}_j^{(0)}$, the matrices $D_t$ and $S_\ell$ are fixed, 
while the matrices $E_t^\top A_k$ have independent $\mathcal{N}(0,\sigma^2)$ entries and are independent across $t$.
Since $H_t = A_{-k} D_t^\top + E_t^\top A_k$, expanding $\|H_t S_\ell^{-1}\|_{\mathrm{F}}^2$ gives
\begin{align*}
    \E\big( \|H_t S_\ell^{-1}\|_{\mathrm{F}}^2
    \mid \mathcal{P}_j^{(0)} \big)
    &= 
    \| A_{-k} D_t^\top S_\ell^{-1} \|_{\mathrm{F}}^2
    + p_{-k} \sigma^2 \tr(S_\ell^{-2}) , \\
    \mathrm{Var}\big( \|H_t S_\ell^{-1}\|_{\mathrm{F}}^2
    \mid \mathcal{P}_j^{(0)} \big)
    &=
    4 \sigma^2
    \| A_{-k} D_t^\top S_\ell^{-2} \|_{\mathrm{F}}^2
    + 2 p_{-k} \sigma^4 \tr(S_\ell^{-4}) .
\end{align*}
Moreover, $A_{-k}^\top A_{-k}=I$
and $\sum_{t\in \mathcal{I}_\ell} D_tD_t^\top = s_\ell S_\ell$ imply
\[
\sum_{t\in \mathcal{I}_\ell}
\| A_{-k} D_t^\top S_\ell^{-1} \|_{\mathrm{F}}^2
=
s_\ell \tr(S_\ell^{-1}),
\quad
\sum_{t\in \mathcal{I}_\ell}
\| A_{-k} D_t^\top S_\ell^{-2} \|_{\mathrm{F}}^2
=
s_\ell \tr(S_\ell^{-3}) .
\]
Summing the preceding mean and variance identities over the independent observations $t\in \mathcal{I}_\ell$ yields
\begin{align}
    \E(Q_\ell \mid \mathcal{P}_j^{(0)})
    &=
    s_\ell \mathrm{tr}(S_\ell^{-1})
    + s_\ell p_{-k} \sigma^2 \mathrm{tr}(S_\ell^{-2}) ,
    \label{eq:adaptive-Q-mean}\\
    \mathrm{Var}(Q_\ell \mid \mathcal{P}_j^{(0)})
    &=
    4 s_\ell \sigma^2 \mathrm{tr}(S_\ell^{-3})
    +
    2 s_\ell p_{-k} \sigma^4 \mathrm{tr}(S_\ell^{-4}) .
    \label{eq:adaptive-Q-variance}
\end{align}

By \eqref{eq:adaptive-V-trace} and \eqref{eq:adaptive-Q-mean}, we have
\begin{align*}
    q_k \E\{\mathrm{tr}(V_\ell) \mid \mathcal{P}_j^{(0)}\}
    &=
    \frac{q_k \sigma^2}{s_\ell}
    \left\{ \mathrm{tr}(S_\ell^{-1})
    + p_{-k} \sigma^2 \mathrm{tr}(S_\ell^{-2}) \right\} .
\end{align*}
By \eqref{eq:adaptive-factor-scale} and Assumption~\ref{ass:latent},
\[
\frac{2 r_k}{3\lambda_{j,1}^{(k)}}
\leq 
\mathrm{tr}(S_\ell^{-1})
\leq 
\frac{2 C_\lambda r_k}{\lambda_{j,1}^{(k)}},
\quad
\frac{4 r_k}{9 (\lambda_{j,1}^{(k)})^2}
\leq 
\mathrm{tr}(S_\ell^{-2})
\leq 
\frac{4 C_\lambda^2 r_k}{(\lambda_{j,1}^{(k)})^2} .
\]
Since $p_k/2 \leq q_k \leq p_k$ and
$b_j\leq s_\ell\leq 2b_j$, it follows that
\begin{equation}\label{eq-lemma-adap-14}
\frac{r_k}{9} \left\{ \frac{\sigma^2 p_k}{b_j \lambda_{j,1}^{(k)}} 
+ \frac{\sigma^4 p}{b_j (\lambda_{j,1}^{(k)})^2} \right\} 
\leq
q_k \E\{\mathrm{tr}(V_\ell) \mid  \mathcal{P}_j^{(0)}\}
\leq
4 C_\lambda^2 r_k \left\{\frac{\sigma^2 p_k}{b_j \lambda_{j,1}^{(k)}} 
+ \frac{\sigma^4 p}{b_j (\lambda_{j,1}^{(k)})^2} \right\} .
\end{equation}
To compare these bounds with the intrinsic radius, recall that
\begin{equation}\label{eq-lemma-adap-15}
\frac{\varrho_{j,k}}{\sqrt{r_k}} =
\frac{\sigma}{\sqrt{\lambda_{j,1}^{(k)}}} \sqrt{\frac{p_k}{b_j}} 
+ \frac{\sigma^2}{\lambda_{j,1}^{(k)}} \sqrt{\frac{p}{b_j}}
+ \frac{\sigma^2p_k}{b_j \lambda_{j,1}^{(k)}} .
\end{equation}
By the condition stated in Lemma~\ref{lem:adaptive-model-calibration}, we have
\[
\frac{\sigma^2 p_k}{b_j \lambda_{j,1}^{(k)}}
\leq 
\frac{1}{C_0 C_b} 
\leq 1.
\]
The last summand in $\varrho_{j,k} /\sqrt{r_k}$ is the square of the first and is therefore no larger than the first. 
Squaring \eqref{eq-lemma-adap-15} gives
\begin{equation}\label{eq-lemma-adap-16}
r_k \left\{ \frac{\sigma^2 p_k}{b_j \lambda_{j,1}^{(k)}}
+ \frac{\sigma^4 p}{b_j (\lambda_{j,1}^{(k)})^2}
\right\} 
\leq 
\varrho_{j,k}^2 
\leq
5 r_k\left\{ \frac{\sigma^2 p_k}{b_j \lambda_{j,1}^{(k)}}
+ \frac{\sigma^4 p}{b_j (\lambda_{j,1}^{(k)})^2}
\right\}.
\end{equation}
Combining \eqref{eq-lemma-adap-14} and \eqref{eq-lemma-adap-16}, we obtain 
\begin{equation}\label{eq:adaptive-V-mean-scale}
\frac{1}{45} \varrho_{j,k}^2 
\leq
q_k \E\{\mathrm{tr}(V_\ell) \mid \mathcal{P}_j^{(0)}\}
\leq 
4 C_\lambda^2 \varrho_{j,k}^2.
\end{equation}
Next, we check that $Q_\ell$ is close to its mean.
By \eqref{eq:adaptive-Q-mean}, 
\eqref{eq:adaptive-Q-variance},
\eqref{eq:adaptive-factor-scale},
and Assumption~\ref{ass:latent}, we have 
\[
\frac{\mathrm{Var}(Q_\ell \mid \mathcal{P}_j^{(0)})}{\{ \E(Q_\ell \mid \mathcal{P}_j^{(0)}) \}^2} 
\leq 
\frac{24 C_\lambda \sigma^2}
{s_\ell r_k \lambda_{j,r_k}^{(k)}}
\leq
\frac{C_7}{C_0 p_k r_k}.
\]
where $C_7 >0$ is an absolute constant and the last inequality used the condition stated in Lemma~\ref{lem:adaptive-model-calibration}.
Hence after increasing $C_0$ if necessary, by Chebyshev's inequality we have
\begin{equation}\label{eq-lemma-adap-17}
\P\left\{
\left| Q_\ell - \E(Q_\ell \mid \mathcal{P}_j^{(0)}) \right| > 
\frac{1}{2} \E(Q_\ell \mid \mathcal{P}_j^{(0)})
\mathrel{\Big|} \mathcal{P}_j^{(0)} \right\}
\leq
4\frac{ \mathrm{Var}(Q_\ell \mid \mathcal{P}_j^{(0)})}{\{ \E(Q_\ell \mid \mathcal{P}_j^{(0)}) \}^2} \leq
\frac{1}{96} .
\end{equation}
Since $q_k \mathrm{tr}(V_\ell) = q_k \sigma^2 Q_\ell /s_\ell^2$,
\eqref{eq:adaptive-V-mean-scale} shows that, outside this exceptional event in \eqref{eq-lemma-adap-17},
\[
\frac{1}{90} \varrho_{j,k}^2
\leq q_k \mathrm{tr}(V_\ell)
\leq 6 C_\lambda^2 \varrho_{j,k}^2 .
\]
The same argument applies to block~$h$. A union bound over the two blocks gives
\begin{equation}
\label{eq:adaptive-tangent-scale}
\P\big( \mathcal{E}_{\ell h}
\mid \mathcal{P}_j^{(0)} \big)
\geq \frac{47}{48} ,
\qquad \text{where \ }
\mathcal{E}_{\ell h}
:=
\left\{
c'' \varrho_{j,k}^2
\leq
2 q_k \mathrm{tr}(V_\ell + V_h)
\leq
C'' \varrho_{j,k}^2
\right\} ,
\end{equation}
for some absolute constants $C''\geq c'' >0$.
In view of \eqref{eq:adaptive-tangent-second-moment},
this establishes the typical conditional second-moment scale of $K_\ell-K_h$.

If $G$ is a centered Gaussian vector with
$v^2 = \E\|G\|_2^2 >0$, Theorem~8 in \cite{carbery2001distributional}, applied to the nonnegative quadratic form $\|G\|_2^2$, gives
\begin{equation}\label{eq:adaptive-Gaussian-small-ball}
\P\{\|G\|_2\leq uv\}\leq C_8 u,
\quad u\in(0,1),
\end{equation}
for an absolute constant $C_8 >0$.
On the event $\mathcal{E}_{\ell h}$ in \eqref{eq:adaptive-tangent-scale},
\eqref{eq:adaptive-tangent-second-moment} implies
\[
v_{\ell h} \geq \sqrt{\frac{c''}{2}} \varrho_{j,k}.
\]
Fix $u\in(0,1)$ and set $c_3 = u \sqrt{c'' /2}$.
Conditional on $\mathcal{P}_j^{(0)}$ and the $H_t$ in the two blocks, apply \eqref{eq:adaptive-Gaussian-small-ball} to $K_\ell-K_h$. 
On the event $\mathcal{E}_{\ell h}$ in \eqref{eq:adaptive-tangent-scale}, this gives
\[
\P\left\{
\|K_\ell-K_h\|_{\mathrm{F}}
\leq c_3 \varrho_{j,k}
\mathrel{\Big|}
\mathcal{P}_j^{(0)}, \{H_t: t\in \mathcal{I}_\ell \cup \mathcal{I}_h\}
\right\} 
\leq C_8 u.
\]
Hence we have
\begin{equation}\label{eq:adaptive-K-small-ball}
\P\left\{\|K_\ell-K_h\|_{\mathrm{F}}
\leq c_3 \varrho_{j,k}
\mathrel{\Big|} \mathcal{P}_j^{(0)} \right\} 
\leq
C_8 u \P(\mathcal{E}_{\ell h} \mid \mathcal{P}_j^{(0)})
+ \P(\mathcal{E}_{\ell h}^{c} \mid \mathcal{P}_j^{(0)})
\leq 
C_8 u + \frac{1}{48}.
\end{equation}

Choose $u\in(0,1)$ sufficiently small  such that $C_8 u + 1/24 \leq 1/3$. 
Choose $c_1 < c_3/\sqrt{2}$ and increase $C_0$ as required by \eqref{eq:adaptive-remainder-control}.
Then
\[
c_- := \sqrt{2} c_3 -2 c_1 >0.
\]
If $\|K_\ell-K_h\|_{\mathrm{F}} > c_3 \varrho_{j,k}$ and both remainders have Frobenius norm at most $c_1 \varrho_{j,k}$ in \eqref{eq:adaptive-remainder-control}, then \eqref{eq-lemma-adap-13} gives
\[
\left\| \widehat\Pi_{j,k}^{(0),[\ell]} - \widehat\Pi_{j,k}^{(0),[h]} \right\|_{\mathrm{F}}
> c_- \varrho_{j,k}.
\]
Hence, by \eqref{eq:adaptive-remainder-control} and \eqref{eq:adaptive-K-small-ball}, we have
\begin{align*}
    &\;\quad
    \P\left\{ \left\| \widehat\Pi_{j,k}^{(0),[\ell]} -\widehat\Pi_{j,k}^{(0),[h]} \right\|_{\mathrm{F}}
    \leq c_- \varrho_{j,k} \mathrel{\Big|} \mathcal{P}_j^{(0)}
    \right\} \\
    &\leq 
    \P\left\{ \|K_\ell-K_h\|_{\mathrm{F}} 
    \leq c_3 \varrho_{j,k} \mathrel{\Big|} \mathcal{P}_j^{(0)}
    \right\}
    + \P\left\{ \|\mathcal{R}_\ell\|_{\mathrm{F}}
    > c_1 \varrho_{j,k} \mathrel{\Big|}
    \mathcal{P}_j^{(0)} \right\}
    + \P\left\{ \|\mathcal{R}_h\|_{\mathrm{F}}
    > c_1 \varrho_{j,k} \mathrel{\Big|}
    \mathcal{P}_j^{(0)} \right\}
    \\
    &\leq
    C_8 u + \frac{1}{48} + \frac{2}{96}
    =
    C_8 u + \frac{1}{24}
    \leq
    \frac{1}{3} .
\end{align*}
Consequently, for every pair of distinct blocks,
\begin{equation}\label{eq:adaptive-pair-lower}
\P\left\{\left\| \widehat\Pi_{j,k}^{(0),[\ell]} - \widehat\Pi_{j,k}^{(0),[h]}  \right\|_{\mathrm{F}}
> c_- \varrho_{j,k} \mathrel{\Big|} \mathcal{P}_j^{(0)}
\right\} 
\geq
\frac{2}{3} .
\end{equation}

\medskip\noindent\textbf{Step 6.} Finally, we conclude this proof as follows.
Since every diagonal distance is zero, the definition of $G_{j,k}$ gives, for $d\geq 0$,
\[
G_{j,k}(d) = \frac{1}{B_j^{(0)}}+ \frac{1}{(B_j^{(0)})^2}
\sum_{\substack{\ell,h \in[B_j^{(0)}] \\ \ell\neq h}}
\P\left\{
\left\| \widehat\Pi_{j,k}^{(0),[\ell]} -\widehat\Pi_{j,k}^{(0),[h]}  \right\|_{\mathrm{F}} 
\leq d \mathrel{\Big|} \mathcal{P}_j^{(0)} 
\right\}.
\]

For each off-diagonal pair, \eqref{eq:adaptive-pair-lower} bounds the probability in this sum by $1/3$ when $d = c_- \cdot \varrho_{j,k}$. Thus, using $B_j^{(0)}\geq 8$, we have
\begin{equation}\label{eq-lemma-adap-final-1}
G_{j,k} (c_- \varrho_{j,k}) 
\leq 
\frac{1}{B_j^{(0)}}
+ \left( 1 -\frac{1}{B_j^{(0)}} \right) \frac{1}{3}
\leq \frac{5}{12}
\leq\frac{3}{4} .
\end{equation}

Similarly, \eqref{eq:adaptive-pair-upper} bounds each off-diagonal probability below by $31/32$ when
$d = c_+ \cdot \varrho_{j,k}$. Therefore,
\begin{equation}\label{eq-lemma-adap-final-2}
G_{j,k} (c_+ \varrho_{j,k}) 
\geq
\frac{1}{B_j^{(0)}}
+ \left( 1 -\frac{1}{B_j^{(0)}} \right) \frac{31}{32}
\geq
\frac{31}{32}
\geq 
\frac{15}{16} .
\end{equation}
Combining \eqref{eq-lemma-adap-final-1} and \eqref{eq-lemma-adap-final-2}, we complete the proof of \eqref{eq:adaptive-model-calibration},
which concludes the proof of this lemma.
\end{proof}

\medskip
\medskip

\begin{lemma}
\label{lem:adaptive-private-calibration}
Suppose that Assumptions~\ref{ass:loading}--\ref{ass:noise} hold.
We use the block sizes and public grid specified in Theorem~\ref{thm:adaptive-upper}.
Suppose also $p_k \geq 2 r_k$,
$p_k \geq C_p \log(n \vee p)$, 
and \eqref{eq:adaptive-upper-conditions-1} holds.
Let $c_-$ and $c_+$ be the constants in Lemma~\ref{lem:adaptive-model-calibration}, and let $C_{\mathrm{cmp}}$ be defined in \eqref{eq:adaptive-radius-comparison}.
Choose the fixed calibration constant $C_A$ and the constant $c_0$ in the upper signal bound so that
\[
C_A
\geq
\frac{C_{\mathrm{cmp}} C_R}{c_-},
\qquad
0 < c_0
\leq
\frac{c_-^2}{C_b + 1}.
\]
Then there exist absolute constants $C_1, C', C, c >0$ such that, if $\varepsilon_j \in(0,1]$ and
\begin{equation} 
\label{eq:adaptive-pilot-size}
n_j \varepsilon_j
\geq 
C_1\{ r_k + \log(n\vee p)\} \log\big\{ m(n\vee p)^{c+1} \big\},
\quad j\in[m],
\end{equation}
then it holds that
\begin{equation} 
\label{eq:adaptive-certificate-sandwich}
\P\big\{ R_{j,k} \leq  \widehat{R}_{j,k} 
\leq C' R_{j,k} ,\mbox{ for all } j\in[m] \big\} 
\geq 
1 - C (n\vee p)^{-c}.
\end{equation}
\end{lemma}

\begin{proof}[Proof of Lemma~\ref{lem:adaptive-private-calibration}]
Define
\[
\mathcal{E}_{\mathrm{bal}}^{(0)} :=
\bigcap_{j=1}^m \bigcap_{\ell=1}^{B_j^{(0)}}
\left\{
\frac{1}{2} \Lambda_{j,(k)}
\preceq
\frac{1}{|\mathcal{I}_{j,\ell}^{(0)}|}
\sum_{t\in \mathcal{I}_{j,\ell}^{(0)}}
\mathcal{M}_k(\mathbf{F}_t^{(j)})
\mathcal{M}_k(\mathbf{F}_t^{(j)})^\top
\preceq
\frac{3}{2} \Lambda_{j,(k)}
\right\} .
\]
Each pilot block is a uniformly selected subset of the corresponding client's records and has size at least $b_j$.
Let $c_{\mathrm{bal}} > 0$ denote the absolute constant in the exponential bound of Lemma~\ref{lem:block-balance}.
By that lemma and a union bound, we have
\begin{align}
    \P\left\{ \left( \mathcal{E}_{\mathrm{bal}}^{(0)} \right)^c \right\}
    &\leq
    2 r_k \sum_{j=1}^m B_j^{(0)}
    \exp\left\{
    - \frac{c_{\mathrm{bal}} b_j}{C_\lambda^2}
    \right\} \nonumber\\
    &\leq
    2 r_k
    (n \vee p)^{-c_{\mathrm{bal}} C_b / C_\lambda^2}
    \cdot \sum_{j=1}^m B_j^{(0)} 
    \leq
    2 (n \vee p)^{2 - c_{\mathrm{bal}} C_b / C_\lambda^2}
    \leq
    2 (n \vee p)^{-c} ,
    \label{eq-lemma-prop-1}
\end{align}
where $C_\lambda$ is the constant in Assumption~\ref{ass:latent} and $c > 0$ is the fixed probability exponent in this lemma;
the second inequality used our choice of block size choice in Theorem~\ref{thm:fed-upper}, which is also used by the adaptive procedure;
the third used $r_k \leq p$ and $\sum_{j=1}^m B_j^{(0)} \leq n$;
and the last holds provided
\[
C_b \geq
\frac{(c + 2) C_\lambda^2}{c_{\mathrm{bal}}} .
\]

We now work on the event $\mathcal{E}_{\mathrm{bal}}^{(0)}$.
Fix $j\in[m]$, condition on $\mathcal{P}_j^{(0)}$,
and write $B :=B_j^{(0)}$.

By \eqref{eq:adaptive-grid} and $r_k\leq p$, we have
\begin{equation}
\label{eq:adaptive-grid-cardinality}
S_k+1 \leq C_\gamma \log(n\vee p),
\end{equation}
where $C_\gamma >0$ is an absolute constant.
Also, note that
\begin{equation}\label{eq:adaptive-direct-block-siz-1}
B \geq
\frac{n_j^{(0)}}{b_j} -1
\geq
\frac{n_j}{2b_j} -2.
\end{equation}
Using $b_j \leq (C_b+1)\{r_k + \log(n\vee p)\}$ and \eqref{eq:adaptive-direct-block-siz-1}, we obtain
\begin{equation}\label{eq:adaptive-direct-block-size-2}
B \varepsilon_j
\geq
\frac{n_j\varepsilon_j}{2(C_b +1) \{r_k + \log(n\vee p)\}}
- 2 \varepsilon_j.
\end{equation}
Since $\varepsilon_j\leq 1$, combining
\eqref{eq:adaptive-pilot-size} and \eqref{eq:adaptive-direct-block-size-2} gives 
\begin{equation}\label{eq:adaptive-direct-block-size}
B \varepsilon_j
\geq 
C_2 \log\{m(n\vee p)^{c+1}\},
\quad
B \geq C_2 \log\{m (n\vee p)^{c+1}\},
\end{equation}
where $0 < C_2 \leq C_1/\{2(C_b+1)\} - 2 /\log(2)$ is a large enough absolute constant. In particular, $B\geq 8$.

Define
\begin{align*}
    \mathcal{E}_{\mathrm{score}}^{(0)} &:=
    \bigcap_{j=1}^m
    \left\{
    \max_{0\leq s\leq S_k}
    \left| q_{j,k}^{(0)}(d_s) - B_j^{(0)}G_{j,k}(d_s) \right|
    \leq 
    \frac{B_j^{(0)}}{64}
    \right\}, \\
    \mathcal{E}_{\mathrm{noise}}^{(0)} &:=
    \bigcap_{j=1}^m
    \left\{
    \max\{ |Y_j|, |Y_{j,0}|, \ldots, |Y_{j,S_k}| \}
    \leq 
    \frac{B_j^{(0)}}{64}
    \right\}.
\end{align*}
By \eqref{eq:adaptive-pair-cdf},
\[
\E\{ q_{j,k}^{(0)}(d_s) \mid \mathcal{P}_j^{(0)} \}
=
B_j^{(0)} G_{j,k}(d_s) .
\]
Conditional on the pilot partitions, block projectors computed from distinct blocks are independent, and changing one block projector changes $q_{j,k}^{(0)}(d_s)$ by at most two. 
The bounded differences inequality
\citep[e.g., Corollary 2.21 in][]{wainwright2019high}
and a union bound over the grid give, for an absolute constant $c_1>0$,
\begin{equation}
\label{eq:adaptive-private-score-deviation}
\P\left\{
(\mathcal{E}_{\mathrm{score}}^{(0)})^c
\mathrel{\Big|}
\mathcal{P}_1^{(0)}, \ldots, \mathcal{P}_m^{(0)}
\right\}
\leq
2(S_k+1) \sum_{j=1}^m e^{-c_1 B_j^{(0)}}.
\end{equation}
The Laplace tail bound and a union bound over the grid give, for an absolute constant $c_2 >0$,
\begin{equation}
\label{eq:adaptive-private-noise-deviation}
\P\left\{
(\mathcal{E}_{\mathrm{noise}}^{(0)})^c
\mathrel{\Big|}
\mathcal{P}_1^{(0)}, \ldots, \mathcal{P}_m^{(0)}
\right\}
\leq
2(S_k+1) \sum_{j=1}^m
e^{-c_2 B_j^{(0)} \varepsilon_j}.
\end{equation}

Combining \eqref{eq:adaptive-grid-cardinality},
\eqref{eq:adaptive-direct-block-size}, 
\eqref{eq:adaptive-private-score-deviation},
and \eqref{eq:adaptive-private-noise-deviation},
we have 
\begin{align}\label{eq-lemma-prop-2}
\P\left\{
(\mathcal{E}_{\mathrm{score}}^{(0)}
\cap\mathcal E_{\mathrm{noise}}^{(0)})^c
\mathrel{\Big|}
\mathcal{P}_1^{(0)}, \ldots, \mathcal{P}_m^{(0)}
\right\}
\leq
\frac{4 C_\gamma \log(n\vee p)}
{(n\vee p)^{c+1}}
\leq 
4 C_\gamma (n\vee p)^{-c}.
\end{align}

Next, by the upper bound in \eqref{eq:adaptive-upper-conditions-1},
the definition \eqref{eq:adaptive-intrinsic-radius},
and $b_j \leq (C_b + 1)\{r_k + \log(n \vee p)\}$,
we have
\begin{align*}
    \varrho_{j,k} &\geq
    \sqrt{r_k}
    \frac{\sigma}{\sqrt{\lambda_{j,1}^{(k)}}}
    \sqrt{\frac{p_k}{b_j}} 
    \geq
    (n \vee p)^{-\gamma}
    \sqrt{\frac{r_k + \log(n \vee p)}{c_0 b_j}} 
    \geq
    \frac{(n \vee p)^{-\gamma}}{\sqrt{c_0 (C_b + 1)}} ,
\end{align*}
so that by the choice $c_0 \leq c_-^2 / (C_b + 1)$, we conclude $d_{\min} = (n \vee p)^{-\gamma} \leq c_- \varrho_{j,k}$.

We now work on the event $\mathcal{E}_{\mathrm{bal}}^{(0)} \cap \mathcal{E}_{\mathrm{score}}^{(0)} \cap \mathcal{E}_{\mathrm{noise}}^{(0)}$. If $d_s\leq c_- \varrho_{j,k}$,
monotonicity of $G_{j,k}$ and Lemma~\ref{lem:adaptive-model-calibration} give $G_{j,k}(d_s)\leq 3/4$.
On $\mathcal{E}_{\mathrm{score}}^{(0)}$,
\[
q_{j,k}^{(0)}(d_s)
\leq 
B G_{j,k}(d_s) + \frac{B}{64}
\leq
\frac{49 B}{64} .
\]
On $\mathcal{E}_{\mathrm{noise}}^{(0)}$, we have $Y_{j,s}\leq B/64$ and $Y_j \geq -B/64$.
Consequently,
\[
q_{j,k}^{(0)}(d_s) +Y_{j,s}
\leq
\frac{50B}{64}
<
\frac{55B}{64}
=
\frac{7B}{8} -\frac{B}{64}
\leq
\frac{7B}{8} +Y_j .
\]
Thus no grid point at or below $c_- \cdot \varrho_{j,k}$ crosses the noisy
threshold.

To bound the selected grid point from above, first suppose that $c_+ \cdot \varrho_{j,k}\leq d_{S_k}$, and let $d_{s_*}$ be the first grid point satisfying $d_{s_*}\geq c_+ \cdot \varrho_{j,k}$. 
By monotonicity of $G_{j,k}$ and Lemma~\ref{lem:adaptive-model-calibration}, we have
\[
G_{j,k}(d_{s_*}) \geq \frac{15}{16}
= \frac{60}{64} .
\]
On $\mathcal{E}_{\mathrm{score}}^{(0)}$, this gives
\[
q_{j,k}^{(0)}(d_{s_*})
\geq 
B G_{j,k}(d_{s_*}) -\frac{B}{64}
\geq
\frac{59 B}{64} .
\]
On $\mathcal{E}_{\mathrm{noise}}^{(0)}$, we have $Y_{j,s_*}\geq -B/64$ and $Y_j\leq B/64$.
Consequently,
\[
q_{j,k}^{(0)}(d_{s_*}) +Y_{j,s_*}
\geq \frac{58 B}{64}
> \frac{57B}{64}
= \frac{7B}{8} + \frac{B}{64}
\geq \frac{7B}{8} +Y_j.
\]
Thus $d_{s_*}$ crosses the noisy threshold, so the search stops no later than $s_*$.

Finally, $d_0=d_{\min} \leq c_- \cdot \varrho_{j,k} < c_+ \cdot \varrho_{j,k}$, so $s_*\geq 1$.
By the construction of $s_*$, $d_{s_*-1}< c_+ \cdot \varrho_{j,k}$. Since consecutive grid points differ by a factor of two, it holds that
\[
d_{\widehat{s}_{j,k}}
\leq  d_{s_*}
= 2 d_{s_*-1}
< 2 c_+ \cdot \varrho_{j,k}.
\]
If instead $d_{S_k} < c_+ \cdot \varrho_{j,k}$, then $d_{S_k} \geq \sqrt{2 r_k}$ by
\eqref{eq:adaptive-grid}. 
The Frobenius distance between any two rank-$r_k$ orthogonal projectors is
at most $\sqrt{2 r_k}$. Hence every ordered block pair is counted at $d_{S_k}$, and
$q_{j,k}^{(0)}(d_{S_k}) =B$.

On $\mathcal{E}_{\mathrm{noise}}^{(0)}$, we have $Y_{j,S_k}\geq -B/64$ and $Y_j\leq B/64$. Therefore,
\[
q_{j,k}^{(0)}(d_{S_k}) +Y_{j,S_k}
\geq
\frac{63 B}{64}
>
\frac{57 B}{64}
=
\frac{7B}{8} +\frac{B}{64}
\geq
\frac{7B}{8} +Y_j.
\]
The terminal grid point crosses the noisy threshold.
The search thus selects a grid point satisfying
\[
d_{\widehat{s}_{j,k}}
\leq d_{S_k}
< c_+ \cdot \varrho_{j,k}.
\]

We have already shown that no grid point at or below $c_- \cdot \varrho_{j,k}$ crosses the noisy threshold.
Combining this lower bound with the upper bound from either case, we obtain, on 
$\mathcal{E}_{\mathrm{bal}}^{(0)} \cap \mathcal{E}_{\mathrm{score}}^{(0)} \cap \mathcal{E}_{\mathrm{noise}}^{(0)}$,
that
\begin{equation}
\label{eq:adaptive-scale-sandwich}
c_- \cdot \varrho_{j,k}
<
d_{\widehat{s}_{j,k}}
\leq
2c_+ \cdot \varrho_{j,k},
\quad j\in[m].
\end{equation}
The conditional probability bound for
$\mathcal{E}_{\mathrm{score}}^{(0)} \cap \mathcal{E}_{\mathrm{noise}}^{(0)}$
obtained above is uniform over the pilot partitions.
Together with \eqref{eq-lemma-prop-1} and \eqref{eq-lemma-prop-2},
it shows that \eqref{eq:adaptive-scale-sandwich} holds with probability at least $1 - C(n\vee p)^{-c}$.

Finally, \eqref{eq:adaptive-radius-comparison} gives $R_{j,k}\asymp \varrho_{j,k}$. 
Since $\widehat{R}_{j,k} = C_A d_{\widehat{s}_{j,k}}$, choose $C_A$ sufficiently large that $C_A c_- \cdot \varrho_{j,k} \geq R_{j,k}$ for all $j$.
Whenever \eqref{eq:adaptive-scale-sandwich} holds,
\[
R_{j,k}
\leq \widehat{R}_{j,k}
\leq 2 C_A c_+ \varrho_{j,k}
\leq C' R_{j,k},
\quad j\in[m].
\]
This proves \eqref{eq:adaptive-certificate-sandwich} 
and thus completes the proof of this lemma.
\end{proof}

\subsection{Proofs of propositions and theorems}

\begin{proof}[Proof of Proposition~\ref{prop:adaptive-FDP}]
Fix $j\in[m]$ and a realization of the data-independent split and block partitions. 
Consider any two neighboring datasets of client~$j$, with all other clients' data held fixed. As the split is fixed, the record on which the two datasets differ belongs to either the pilot half or the main half.

\medskip
\noindent
\textbf{Case 1: the differing record belongs to the pilot half.}
Changing this record changes at most one pilot block projector. It hence changes each query in the definition of $q_{j,k}^{(0)}(d_s)$ (in Algorithm~\ref{alg:adaptive-fed-est}) by at most two. 
With the noises $Y_j$ and $Y_{j,s}$ (also in Algorithm~\ref{alg:adaptive-fed-est}), the output $\widehat{s}_{j,k}$ is $(\varepsilon_j, 0)$-CDP.

The main observations are identical in the two datasets.
Once $\widehat{s}_{j,k}$ is given, the quantities $\widehat{R}_{j,k}$, $\widehat\tau_{j,k}$, $\widehat\Psi_{j,k}$, and $\widehat\iota_{j,k}$ are fixed,
and the remaining output is generated from the same main observations under both datasets. 
Thus the complete local message is a randomized post-processing of $\widehat{s}_{j,k}$. It is $(\varepsilon_j, 0)$-CDP in this case.

\medskip
\noindent
\textbf{Case 2: the differing record belongs to the main half.}
The pilot observations are identical in the two datasets,
so $\widehat{s}_{j,k}$ has the same distribution under both.
Condition on any value of $\widehat{s}_{j,k}$. The calibration parameters, certificate, and flag are now fixed. 
If $\widehat\iota_{j,k}=1$, Proposition~\ref{prop:standard-FDP} gives $(\varepsilon_j,\delta_j)$-CDP for the augmented main output.  If $\widehat\iota_{j,k}=0$, the main output is drawn independently of the data. Hence for every fixed pilot output, the complete local message is $(\varepsilon_j,\delta_j)$-CDP with respect to the main observations.  Averaging over the common distribution of the pilot output preserves this guarantee.

These two cases cover every pair of neighboring datasets.
Thus, for each fixed split and block partition, the local message is $(\varepsilon_j,\delta_j)$-CDP. 
The split and block partitions are generated independently of the data, so averaging the privacy inequality over their randomness gives the same guarantee for Algorithm~\ref{alg:adaptive-fed-est} without conditioning on a particular partition.

Finally, the server weights and final singular value decomposition are post-processing of the client messages.
Applying the argument to each client proves $(\boldsymbol\varepsilon, \boldsymbol\delta)$-FDP.
\end{proof}

\begin{proof}[Proof of Theorem~\ref{thm:adaptive-upper}]
For notational convenience, write $P:=A_kA_k^\top$ throughout this proof.
Our proof proceeds in three steps, detailed as follows.

\medskip\noindent\textbf{Step 1. Pilot calibration and screening.}
Apply Lemma~\ref{lem:block-balance} to all pilot and main blocks. The same union-bound argument as in its proof and a sufficiently large $C_b$ show that simultaneous factor balance fails with probability at most $C(n\vee p)^{-c}$. 
On the factor-balance event, apply Lemma~\ref{lem:adaptive-private-calibration}.
Define
\[
\big( \tau_{j,k}^{(1)} \big)^2 :=
\frac{64 R_{j,k}^2}{(B_j^{(1)})^2 \varepsilon_j^2}
\log\left( \frac{2.5}{\delta_j} \right),
\qquad
\Psi_{j,k}^{(1)} :=
\frac{R_{j,k}}{\sqrt{r_k}}
\sqrt{\frac{\log(n \vee p)}{B_j^{(1)}}}
+ \tau_{j,k}^{(1)}
\sqrt{p_k + \log(n \vee p)}.
\]
Then simultaneously for all clients,we have
\begin{equation}
\label{eq:adaptive-radius-sandwich-proof}
R_{j,k}
\leq
\widehat{R}_{j,k}
\leq
C R_{j,k},
\quad
c \Psi_{j,k}^{(1)}
\leq
\widehat\Psi_{j,k}
\leq
C \Psi_{j,k}^{(1)}.
\quad j \in [m],
\end{equation}

By the lower bound in \eqref{eq:adaptive-upper-conditions-1}, the block-size choice,
and \eqref{eq:adaptive-radius-sandwich-proof},
$\widehat{R}_{j,k}$ is uniformly small when $C_0$ is sufficiently large.
The condition \eqref{eq:adaptive-upper-conditions-2} also ensures $\sqrt{r_k} \widehat\Psi_{j,k} \leq c_0$
and the requirement in \eqref{eq:adaptive-screen},
after increasing $C_0$ if necessary.
Hence we conclude $\widehat{\mathcal{G}}_k = [m]$ on the calibration event.

\medskip\noindent\textbf{Step 2. Conditional local ideal releases.}
Let $\mathcal{T}_k$ be the sigma-algebra generated by the data-independent split and partitions, the pilot observations, and the pilot mechanism. 
Conditional on $\mathcal{T}_k$, all adaptive radii, certificates, and server weights are fixed, 
while the main observation noises remain independent across clients and blocks.

For each main block,
repeat all the construction in
\eqref{eq:truncated-block-projector}--\eqref{eq:proof-only-projector}, replacing
\[
B_j ,\quad 
R_{j,k} ,\quad
\tau_{j,k} ,\quad
\Psi_{j,k}^{\mathrm{loc}}
\]
by
\[
B_j^{(1)} ,\quad
\widehat{R}_{j,k} ,\quad
\widehat\tau_{j,k} ,\quad
\widehat\Psi_{j,k} ,
\]
respectively. 
Denote the resulting projector by $\breve\Pi_{j,k}^{\mathrm{ad}}$. The first inequality in \eqref{eq:adaptive-radius-sandwich-proof} and the main-block version of Lemma~\ref{lem:block-projector} show that the truncation agrees with every empirical main-block projector except on an event of probability $C(n\vee p)^{-c}$.

The conditional Hilbert-space Hoeffding bound in \eqref{eq:Hilbert-Hoeffding-local} and the GOE bound in \eqref{eq:GOE-tail} remain valid after these replacements. 
Since $p_k\gtrsim \log(n\vee p)$, their right-hand side is bounded by $C \sqrt{r_k} \widehat\Psi_{j,k}$. 
The reflection argument in \eqref{eq:truncated-reflection-invariance}--\eqref{eq:ideal-expectation-cross-block} also remains valid
conditionally on $\mathcal{T}_k$. Consequently,
\begin{align}
    (I_{p_k}-P) \E(\breve\Pi_{j,k}^{\mathrm{ad}} \mid \mathcal{T}_k) P
    &= 
    0,
    \label{eq:adaptive-conditional-cross-block}\\
    \left\|
    \E(\breve\Pi_{j,k}^{\mathrm{ad}} \mid \mathcal{T}_k) -P
    \right\|
    &\leq 
    c,
    \label{eq:adaptive-conditional-mean-close}\\
    \left\|
    \breve\Pi_{j,k}^{\mathrm{ad}}
    - \E(\breve\Pi_{j,k}^{\mathrm{ad}} \mid \mathcal{T}_k)
    \right\|_{\mathrm{F}}
    &\leq
    C \sqrt{r_k} \widehat\Psi_{j,k}
    \quad\text{almost surely}.
    \label{eq:adaptive-conditional-centered-bound}
\end{align}
The projectors are conditionally independent across $j$.

It remains to verify the coupling. 
On the simultaneous block concentration event, every pair of main-block projectors is at distance at most $R_{j,k}/8 \leq \widehat{R}_{j,k}/8$. Hence, exactly as in \eqref{eq:Huber-equals-block-average}, the Huber aggregate equals the arithmetic block average.
Note that $B_j^{(1)} \geq n_j /(2 b_j) - 1$.
Hence \eqref{eq:adaptive-upper-conditions-2} and our choice of block size imply
\[
B_j^{(1)} \varepsilon_j \geq
C_{\mathrm{test}}
\log\left\{ \frac{n \vee p}{\delta_j} \right\},
\]
where $C_{\mathrm{test}} > 0$ is an absolute constant that can be made sufficiently large by increasing $C_0$.
Conditional on the simultaneous block concentration event, the main stability score equals $B_j^{(1)}$ for every client, so that by \eqref{eq:ideal-test-failure} and taking a union bound over $m \leq n$ clients, we conclude that every main stability test passes,
except on an event of conditional probability at most $C (n \vee p)^{-c}$.
Hence we have
\begin{equation}
\label{eq:adaptive-conditional-coupling}
\P\left\{
\breve\Pi_{j,k}^{\mathrm{ad}}
= \widetilde{U}_{j,k} \widetilde{U}_{j,k}^\top ,
\; \forall j\in[m]
\mathrel{\Big|}
\mathcal{T}_k
\right\}
\geq
1 - C(n\vee p)^{-c} .
\end{equation}

\medskip\noindent\textbf{Step 3. Adaptive server aggregation.}
Conditional on $\mathcal{T}_k$, define
\[
M_k^{\mathrm{ad}} :=
\sum_{j=1}^m \widehat{w}_{j,k}
\E(\breve\Pi_{j,k}^{\mathrm{ad}} \mid \mathcal{T}_k) , 
\quad
\breve{M}_k^{\mathrm{ad}} :=
\sum_{j=1}^m \widehat{w}_{j,k} \breve\Pi_{j,k}^{\mathrm{ad}} .
\]
Then \eqref{eq:adaptive-conditional-cross-block} and \eqref{eq:adaptive-conditional-mean-close} show that
$M_k^{\mathrm{ad}}$ is block diagonal with respect to $\mathrm{span}(A_k) \oplus \mathrm{span}(A_k)^\perp$, 
its leading $r_k$-dimensional eigenspace is $\mathrm{span}(A_k)$, 
and its eigengap is at least an absolute positive constant.

Conditional independence and \eqref{eq:adaptive-conditional-centered-bound}, followed by the Hilbert-space Hoeffding inequality, give, for every $u>0$,
\begin{equation}
\label{eq:adaptive-conditional-Hoeffding}
\P\left\{
\|\breve{M}_k^{\mathrm{ad}} -M_k^{\mathrm{ad}}\|_{\mathrm{F}}^2
> C r_k u  \sum_{j=1}^m \widehat{w}_{j,k}^2 \widehat\Psi_{j,k}^2
\mathrel{\Big|} \mathcal{T}_k
\right\}
\leq
2 e^{-u} .
\end{equation}
By the definition of $\widehat{w}_{j,k}$ in Algorithm~\ref{alg:adaptive-fed-est}, we have
\begin{equation}
\label{eq:adaptive-weight-identity}
\sum_{j=1}^m \widehat{w}_{j,k}^2 \widehat\Psi_{j,k}^2
=
\frac{1}{\sum_{j=1}^m \widehat\Psi_{j,k}^{-2}} .
\end{equation}
Take $u=\log(n\vee p)$ in \eqref{eq:adaptive-conditional-Hoeffding}. 
By the Davis--Kahan theorem \citep[Theorem~2]{yu2015useful}, 
\eqref{eq:adaptive-conditional-coupling}, 
and \eqref{eq:adaptive-weight-identity}, we have
\[
\big\| \widehat{A}_k^{\mathrm{ad}}
(\widehat{A}_k^{\mathrm{ad}})^\top - P \big\|_{\mathrm{F}}^2
\leq
\frac{C r_k \log(n\vee p)}{\sum_{j=1}^m \widehat\Psi_{j,k}^{-2}}
\wedge (2 r_k) .
\]

By \eqref{eq:adaptive-radius-sandwich-proof}, we have
\[
\sum_{j=1}^m \widehat\Psi_{j,k}^{-2}
\geq
\frac{1}{C}
\sum_{j=1}^m (\Psi_{j,k}^{(1)})^{-2} .
\]
Moreover, the lower bound on $B_j^{(1)}$ implies $B_j^{(1)}\asymp B_j$, as the main sample contains a fixed fraction of the observations. 
Together with $p_k \gtrsim \log(n\vee p)$, this gives
\[
\Psi_{j,k}^{(1)} \leq
C \Psi_{j,k}^{\mathrm{loc}} .
\]
Altogether, we conclude the desired upper bound,
thus concluding the stated probability guarantee and this proof.
\end{proof}

\end{document}